\documentclass[11pt,a4,onesided]{amsart}
 \usepackage{epsfig}
\usepackage{amssymb}
\usepackage{color}
\usepackage{mathrsfs}
\usepackage{hyperref}
\newcommand\N{{\mathbb N}}

\newcommand\R{{\mathbb R}}

\newcommand\C{{\mathbb C}}

\newcommand\Pp{{\mathbb P}}

\def\AA{{\mathcal A}}
\def\BB{{\mathcal B}}

\def\DD{{\mathcal D}}
\def\EE{{\mathcal E}}
\def\FF{{\mathcal F}}

\def\JJ{{\mathcal J}}
\def\KK{{\mathcal K}}
\def\LL{{\mathcal L}}

\def\RR{{\mathcal R}}

\def\TT{{\mathcal T}}
\def\UU{{\mathcal U}}
\def\VV{{\mathcal V}}
\def\WW{{\mathcal W}}

\def\BBB{{\mathscr B}}
\def\CCC{{\mathscr C}}

\def\GGG{{\mathscr G}}

\def\KKK{{\mathscr K}}

\def\UUU{{\mathscr U}}
\def\VVV{{\mathscr V}}

\DeclareMathOperator{\sign}{sign}

\newcommand{\dd}{{\, \mathrm d}}

\let\oldmarginpar\marginpar
\renewcommand\marginpar[1]{\-\oldmarginpar[\raggedleft\footnotesize #1]%
{\raggedright\footnotesize #1}}

\def\eps{{\varepsilon}}

\newtheorem{theo}{Theorem}
\newtheorem{prop}[theo]{Proposition}

\newtheorem{lem}[theo]{Lemma}
\newtheorem{cor}[theo]{Corollary}

\theoremstyle{definition}
\newtheorem{defin}[theo]{Definition}

\theoremstyle{remark}
\newtheorem{rem}[theo]{Remark}

\newcommand{\beqn}{\begin{equation}}
\newcommand{\eeqn}{\end{equation}}
\newcommand{\bear}{\begin{eqnarray}}
\newcommand{\eear}{\end{eqnarray}}
\newcommand{\bean}{\begin{eqnarray*}}
\newcommand{\eean}{\end{eqnarray*}}

\def\Nt{|\hskip-0.04cm|\hskip-0.04cm|}
\def\signsm{\bigskip \begin{center} {\sc St\'ephane
      Mischler\par\vspace{3mm} 
      Universit\'e Paris-Dauphine \& IUF\par
      CEREMADE, UMR CNRS 7534\par
      Place du Mar\'echal de Lattre de Tassigny \par
      75775 Paris Cedex 16 FRANCE\par\vspace{3mm} e-mail:}
    \tt{mischler@ceremade.dauphine.fr} \end{center}}

\def\signjs{\bigskip \begin{center} {\sc Justine Scher
\par\vspace{3mm}
  Universit\'e Paris-Dauphine \par
      CEREMADE, UMR CNRS 7534\par
      Place du Mar\'echal de Lattre de Tassigny \par
      75775 Paris Cedex 16 FRANCE}
 \end{center}}

\begin{document}

\title[Spectral analysis and growth-fragmentation]{Spectral analysis of semigroups and growth-fragmentation equations}
%\title[Spectral analysis and growth-fragmentation]{Spectral analysis of semigroups in Banach spaces and growth-fragmentation equation}

%{Semigroup spectral Analysis illustrated on the growth-fragmentation equation}

\author{S. Mischler, J. Scher}

% \footnotetext[1]{Department of Mathematics, The University of Texas
%   at Austin, 1 University Station C1200, Texas 78712, USA.  E-mail:
%   \texttt{gualdani@math.utexas.edu}}

% \footnotetext[2]{CEREMADE, Universit\'e Paris IX-Dauphine, Place du
%   Mar\'echal de Lattre de Tassigny, 75775 Paris, France.  E-mail:
%   \texttt{mischler@ceremade.dauphine.fr}}

% \footnotetext[3]{CEREMADE, Universit\'e Paris IX-Dauphine, Place du
%   Mar\'echal de Lattre de Tassigny, 75775 Paris, France.  E-mail:
%   \texttt{cmouhot@ceremade.dauphine.fr}}

\maketitle

\begin{center} {\bf %Preliminary v
Version of \today}
\end{center}

\begin{abstract} The aim of this paper is twofold:

(1) On the one hand, the paper revisits the spectral analysis of semigroups in a general Banach
space setting. It presents some new and more general versions, and provides comprehensible proofs, 
of classical results such as the spectral mapping theorem, 
some (quantified) Weyl's Theorems and the Krein-Rutman Theorem.
Motivated by evolution PDE  applications, the results apply to a wide and natural class of generators which split as a dissipative part plus a more regular part, 
without assuming any symmetric structure on the operators nor Hilbert structure on the space, and give some growth estimates and spectral gap estimates for the associated semigroup.
The  approach relies on some factorization and summation arguments reminiscent of the Dyson-Phillips  series in the spirit of those used in  \cite{Mcmp,MMcmp,GMM,MM*}. 

(2) On the other hand, we present the semigroup spectral analysis for three important classes of ``growth-fragmentation" equations, namely  the cell division  equation, the self-similar fragmentation equation and the McKendrick-Von Foerster age structured population equation.
By showing that these models lie in the class of equations for which our general semigroup analysis theory applies, 
we prove the exponential rate of convergence of the solutions to the associated  first eigenfunction or self-similar profile for a very large and natural class of fragmentation rates.
Our results generalize similar estimates  obtained in  \cite{MR2114128,MR2536450} for the cell division model with (almost) constant total fragmentation rate and in \cite{MR2832638,MR2821681} for the self-similar fragmentation equation and the cell division equation restricted to 
 smooth and  positive fragmentation rate and  total fragmentation rate which does not increase more
rapidly than quadratically.  It also improves the convergence results without rate obtained in \cite{MMP,EMRR} which have been established 
under similar assumptions to those made in the present work. 

%under strong assumptions of  smoothness and  positivity on the fragmentation rate and the restriction for total fragmentation rate which does not increase more
%rapidly than quadratically.   
  \end{abstract}

\bigskip

Partial differential equations, Asymptotic behavior of solutions
 [see also 45C05, 45K05,  35410],
92D25 
Population dynamics [see also 92C37, 82D60]

\bigskip


%\newpage

\tableofcontents

%%%%%%%%%%%%%%%%%%%% Introduction %%%%%%%%%%%%%%%%%%%%%%%%%%%%%%%%%%%

\bigskip
\section{Introduction}
\label{sec:intro}
\setcounter{equation}{0}
\setcounter{theo}{0}

\smallskip
This paper deals with the study of decay properties for 
$C_0$-semigroup  of bounded and linear operators and their link with spectral properties of  their generator in a Banach framework as well as some applications to the long-time asymptotic of growth-fragmentation equations. 
 
\subsection{Spectral analysis of semigroups}
The study of spectral property of (unbounded) operators and of $C_0$-semigroups of operators has a long history which goes back 
(at least) to the  formalization of functional analysis by D.~Hilbert \cite{MR1511713,MR0056184} at the beginning of the 20th century for the first issue
and surely before the modern theoretical formalization of  $C_0$-semigroups of operators in general Banach spaces
impulsed by  E. Hille and K. Yosida \cite{MR0006610,MR0028537,MR0025077} in the 1940's  for the second issue. For a given operator $\Lambda$
on a Banach space $X$ which generates a $C_0$-semigroup $S_\Lambda(t) = e^{t\Lambda}$ of bounded operators, the two following issues are of 
major importance:

\smallskip
$\bullet$ describe its spectrum $\Sigma(\Lambda)$,  the set of its eigenvalues  and the associated eigenspaces;

\smallskip
$\bullet$ prove the spectral mapping theorem 
%$\Sigma(e^{t\Lambda}) = e^{t\Sigma(\Lambda)}$ 
\beqn\label{eq:SpectralMapping}
\Sigma(e^{t\Lambda}) \backslash \{ 0 \} =  e^{t \Sigma(\Lambda)}, 
\eeqn
and deduce the
asymptotical behaviour of trajectories associated to the semigroup. 

\smallskip
Although it is well-known that the first issue can be a complicated task and the second issue is false in general (see \cite{MR0089373} or \cite[Section IV.3.a]{EngelNagel} for some counterexamples), there exists 
some particular classes of operators (among which is the class of  self-adjoint operators  with compact resolvent  in a Hilbert space) for
which these problems can be completely solved. In the present paper, motivated by evolution partial differential equations  applications and 
inspired by the recent paper \cite{GMM} (see also \cite{Mcmp,MMcmp,MM*}),  we identify a class of operators which 
split as
\beqn\label{eq:L=A+B}
\Lambda = \AA + \BB,
\eeqn 
where $\AA$ is ``much more regular than $\BB$" and $\BB$ has some dissipative property (and then a good localization of its spectrum) for which a positive answer can be given. The dissipative property assumption  we make can be formulated in terms of the 
 time indexed family of iterated time convolution operators $(\AA S_\BB)^{(*k)}(t)$ in the following way
%(so that we may show that the spectrum of $\Lambda$ is not ``too different" from the spectrum of $\BB$). 
%More precisely, we assume that $\BB$ is dissipative and
% satisfies the growth 
% and that 
 \begin{itemize}
 \item[{\bf (H1)}]    for some $a^* \in \R$ and  for any $a > a^*$, $\ell\ge 0$,   there exists a constant $C_{a,\ell}  \in (0,\infty)$  
  such that the following growth estimate holds
\beqn\label{eq:hypH1}
\,\, \forall \, t \ge 0, \qquad 
\| S_\BB *  (\AA S_\BB)^{(*\ell)}(t) \|_{\BBB(X)}  \le  C_{a,\ell} \, e^{at}. 
\eeqn
    \end{itemize}
On the other hand, we make the key assumption that some iterated enough time convolution enjoys 
the growth and regularizing estimate:  
  \begin{itemize}
   \item[{\bf (H2-3)}]  there exist an integer $n \ge 1$ such that for any $a > a^*$, there holds 
\beqn\label{eq:hypH3-2}
\,\, \forall \, t \ge 0, \qquad 
\|  (\AA S_\BB)^{(*n)}(t) \|_{\BBB(X,Y)}  \le  C_{a,n,Y} \, e^{at},
\eeqn
or 
\beqn\label{eq:hypH2-3}
\,\, \forall \, t \ge 0, \qquad 
\|  S_\BB * (\AA S_\BB)^{(*n)}(t) \|_{\BBB(X,Y)}  \le  C_{a,n,Y} \, e^{at},
\eeqn
for some suitable subspace $Y \subset X$ and a constant $C_{a,n,Y}  \in (0,\infty)$.    
   \end{itemize}

\noindent
In assumption {\bf (H2-3)}
 we will typically assume that $Y \subset D(\Lambda^\zeta)$, $\zeta >0$,  and/or
$Y \subset X$ with compact embedding. 

\smallskip
Roughly speaking, for such a class of operators, we will obtain the following set of results:

\smallskip
$\bullet$ {\it Spectral mapping theorem.}  We prove  a partial, but principal, spectral mapping theorem which asserts that 
\beqn\label{eq:PartialSpectralMapping}
\Sigma(e^{t\Lambda}) \cap \Delta_{e^{at}}  = e^{t\Sigma(\Lambda) \cap \Delta_a},
\qquad\forall\, t\ge0, \,\, \forall \, a > a^*, 
\eeqn
where we define the half-plane $\Delta_a := \{ \xi \in \C; \,\, \Re e\xi > a \}$ for any $a \in \R$.
Although \eqref{eq:PartialSpectralMapping} is less accurate than \eqref{eq:SpectralMapping}, 
it is strong enough to describe the semigroup evolution at first order in many situations. 
In particular, it implies that the spectral bound $s(\Lambda)$ and the growth bound $\omega(\Lambda)$ coincide if they are at the right hand side of $a^*$, or in other words  
\beqn\label{eq:s=omega}
\max(s(\Lambda),a^*) = \max(\omega(\Lambda),a^*),
\eeqn
and it gives even more accurate asymptotic information on the semigroup  whenever $\Sigma(\Lambda) \cap \Delta_{a^*} \not=\emptyset$.

\smallskip
$\bullet$ {\it Weyl's Theorem. } We prove some (quantified) version of the Weyl's Theorem which asserts that the part of the spectrum $\Sigma(\Lambda) \cap \Delta_{a^*}$ consists only of {\it discrete eigenvalues} and we get some information on the localization and number of eigenvalues as well as some estimates on the total dimension of the associated sum of eigenspaces.

\smallskip
$\bullet$  {\it Krein-Rutman Theorem. }  We prove some (possibly quantified) version of the Krein-Rutman Theorem under some additional (strict) positivity hypothesizes on the generator $\Lambda$ and the semigroup $S_\Lambda$. 

%\smallskip
%$\bullet$ we extend these results to the case where $\BB$ is only {\it sub-dissipative}. 

\medskip
Let us  describe our approach in order to get the above mentioned {\it ``spectral mapping theorem"}, this one being the key result in order to get our versions of Weyl's Theorem and  Krein-Rutman Theorem. Following \cite{Mcmp,MMcmp,GMM} 
(and many authors before!), the spectral analysis of the operator $\Lambda$ with splitting structure \eqref{eq:L=A+B} is performed by writing the resolvent factorization identity (with our definition of the resolvent in \eqref{def:Resolv})
\beqn\label{eq:RLRB-RBARL}
\RR_\Lambda(z) = \RR_\BB(z) - \RR_\BB(z) \AA \RR_\Lambda(z) 
\eeqn
as well as 
\beqn\label{eq:RLRB-RLARB}
\RR_\Lambda(z) = \RR_\BB(z) - \RR_\Lambda(z) \AA \RR_\BB(z) , 
\eeqn
or an iterative version of \eqref{eq:RLRB-RLARB}, and by exploiting the information that one can deduce from \eqref{eq:hypH1}, \eqref{eq:hypH3-2}
and \eqref{eq:hypH2-3} at the level of the resolvent operators. 
 
 At the level of the semigroup, \eqref{eq:RLRB-RBARL}  yields  the Duhamel formula
 $$
S_\Lambda =S_\BB + S_\BB *  (\AA S_\Lambda),
$$ 
in its classic form, and  \eqref{eq:RLRB-RLARB} yields the Duhamel formula 
\beqn\label{eq:DuhamelSL=SB+SLASB}
S_\Lambda =S_\BB + S_\Lambda *  (\AA S_\BB), 
\eeqn 
in a maybe less standard form (but also reminiscent of perturbation arguments in semigroup theory). 

 \smallskip
 On the other hand, iterating infinitely one of the above identities, it is well-known since the seminal articles by Dyson and Phillips \cite{MR0028203,MR0054167}, that $S_\Lambda$ can be expended  as the Dyson-Phillips series
\beqn\label{eq:DPseries}
S_\Lambda = \sum_{\ell=0}^\infty S_\BB *  (\AA S_\BB)^{(*\ell)},
\eeqn
as soon as the right hand side series converges, and that matter has not an easy answer in general. The summability issue of the Dyson-Phillips  series can be  circumvented by considering the finite iteration of the  Duhamel formula \eqref{eq:DuhamelSL=SB+SLASB}%stopped Dyson-Phillips  series 
\beqn\label{eq:stopDyson}
S_\Lambda = \sum_{\ell=0}^{n-1} S_\BB *  (\AA S_\BB)^{(*\ell)} +  S_\Lambda * (\AA S_\BB)^{(*n)}.
\eeqn
It is a usual trick in order to establish eventual norm continuity (see \cite{MR1753014} 
and the references therein) and it has been also recently used in \cite{GMM,MM*}  
in order to enlarge  or to shrink the functional space where the semigroup $S_\Lambda$ satisfies  some spectral gap estimate. 
 
In the present work, and in the case $\Delta_{a^*} \cap \Sigma(\Lambda) = \{ \lambda \}$, where $\lambda \in \C$ is a semisimple eigenvalue and   $a^* < a< \Re e \lambda$  in order to make the discussion simpler, 
our spectral mapping theorem simply follows by using the classical representation of the semigroup by means of the inverse Laplace transform of the resolvent, as already established by Hille in \cite{MR0006610}. More precisely, we may write 
%\bean
%S_\Lambda = \Bigl( \, \sum_{k=0}^{n}  (-1)^k S_\BB *  (\AA S_\BB)^{(*k)} \Bigr) \Bigl( \sum_{\ell=0}^\infty(-1)^{(n+1)\ell} (\AA S_\BB)^{(*(n+1) \ell)}  \Bigr)
%\eean
\bear\label{eq:stopDysonBIS}
S_\Lambda (t)  
&=& \Pi_{\Lambda,\lambda} \, e^{\lambda t} + \sum_{\ell=0}^{N-1}  \Pi_{\Lambda,\lambda}^\perp S_\BB *  (\AA S_\BB)^{(*\ell)} (t)
\\ \nonumber
&&+  {i \over 2\pi}  \,  \int_{a-i\infty}^{a+i\infty}  e^{zt} \,   (-1)^N \, 
\Pi_{\Lambda,\lambda}^\perp \RR_\Lambda(z) \, (\AA \RR_\BB(z))^{N}  \, dz, 
%\\ \nonumber
%&&+  {i \over 2\pi} \int_{i\R \backslash[-iM,iM]}   e^{zt} \,   \sum_{\ell=N+1}^\infty(-1)^{\ell} 
%\RR_\BB(z) \, (\AA \RR_\BB(z))^{\ell}  \, dz
\eear
for $N$ large enough, in such a way that each term is  appropriately bounded. Here  $\Pi_{\Lambda,\lambda}$ stands for the projector
on the (finite dimensional) eigenspace associated to $\lambda$ and $\Pi_{\Lambda,\lambda}^\perp:= I - \Pi_{\Lambda,\lambda}$. From that formula, we will deduce the spectral mapping theorem \eqref{eq:PartialSpectralMapping}, and more importantly for us, we will generalize the Liapunov result \cite{MR0208093}  about the asymptotic behaviour of trajectories to that class of equations. It is worth emphasizing that our result characterizes the class of operators with ``separable spectrum" for which the partial spectral mapping theorem holds (in the sense that we exhibit a condition which is not only sufficient but also necessary!)
and then, in some sense,  we prove for {\it general semigroup} the spectral mapping theorem known for {\it analytic semigroup} or more generally {\it eventually continuous semigroup} \cite{MR0044737,MR0089373} and for {\it general Banach space}  the partial spectral mapping theorem in a {\it Hilbert space} framework that one can deduce from the Gearhart-Pr\"uss Theorem  \cite{MR0461206,MR743749}. 

\smallskip
With such a representation formula at hand, the precise analysis of the semigroup $S_\Lambda$ reduces to the analysis of the spectrum of $\Lambda$ at the right hand side of $a$, for any $a > a^*$. In other words, the next fundamental issue consists in describing the part of the spectrum $\Sigma(\Lambda) \cap \Delta_a$ in order to take advantage of the information given by \eqref{eq:PartialSpectralMapping}. The simplest situation is when $\Sigma(\Lambda) \cap \Delta_a$ only contains eigenvalues with finite (algebraical) multiplicity which is the situation one gets when one can apply Weyl's Theorem \cite{MR1511560} (see also \cite[Theorem IV.5.35]{Kato}). In our second main result in an abstract setting, we recover Voigt's version \cite{Voigt80} of Weyl's Theorem (for which we give a comprehensive and elementary proof) and we deduce a characterization of semigroups in a general Banach space for which the partial spectral mapping theorem holds with finite and discrete eigenspectrum in $\Delta_a$. We must emphasize that our proof is very simple 
(it exclusively uses  the Fredholm alternative  \cite{MR1554993} in its most basic form) and in particular does not use the essential spectrum set nor the Fredholm operators theory.
Moreover we are able to formulate a quantified version of the Weyl's theorem in the sense that we exhibit a bound on the total dimension of the eigenspaces associated to the discrete eigenvalues which lie in  $\Sigma(\Lambda) \cap \Delta_a$.

\smallskip
In order to describe in a more accurate way the  spectrum $\Sigma(\Lambda)$, one of the most popular techniques is to use a self-adjointness argument for the operator $\Lambda$ as an infinite dimensional generalization of the symmetric structure of matrix. That  implies  $\Sigma(\Lambda) \subset \R$ and, together with Weyl's Theorem, leads to a completely satisfactory description of the  operator's spectrum and the dynamics of the associated semigroup. 
One of the most famous application of that strategy is due to Carleman \cite{Carleman}  who carried on with the study of the linearized space homogeneous Boltzmann equation initiated by Hilbert and who obtained the spectral gap for the associated operator by combining Weyl's Theorem together with the symmetry  of the operator \cite{MR1511713} and the regularity of the gain term \cite{MR0056184}  
(see also the work  by Grad \cite{Grad1,Grad2} on the same issue, and the work by Ukai \cite{MR0363332} on the more complicated space non-homogeneous setting). It is worth emphasizing that this kind of hilbertian arguments have been recently extended to a class of operators, named as {\it ``hypocoercive operators"}, which split as a self-adjoint partially coercive operator plus an anti-adjoint operator. For such an operator one can exhibit an equivalent Hilbert norm which is also a Liapunov function for the associated evolution dynamics and then provides a spectral gap between the first eigenvalue and the remainder of the spectrum. We refer the interested reader to \cite{MR2339441} for a pedagogical introduction as well as to  \cite{MR2034753,MR2130405,MNeu,MR2215889,MR2294477,MR2562709,DMScras,DMS} for some of the original articles.

\smallskip 
In the seminal work \cite{Mcmp}, C. Mouhot started an abstract theory of  ``{\it enlargement} of the functional space of spectral analysis of operators"  
which aims to carry on 
the spectral knowledge on an operator $\Lambda$ and its associated semigroup $S_\Lambda$  in some space $E$ (typically a ``small Hilbert space" in which some self-adjointness structure can be exploited) to another larger general Banach space $\EE \supset E$.  It is worth emphasizing that the ``{\it enlargement} of functional spaces" trick for spectral analysis is reminiscent of several earlier papers  on Boltzmann equations and on Fokker-Planck type equations where, however,  the arguments are intermingled with some nonlinear stability arguments  \cite{MR946973,MR900501,MR1338453}  and/or reduced to some particular evolution PDE in some situation where explicit eigenbasis can be exhibited \cite{Bobylev1975,MR1912106}. 
While \cite{Mcmp} was focused on the linearized space homogeneous Boltzmann operator and the results applied to  sectorial operators, the ``{\it extension}  theory" (we mean here {\it ``enlargement"} or {\it ``shrinkage"}  of the functional space) has been next developed  in the series of papers  \cite{MMcmp,GMM,MM*} in order to deal with non-sectorial operators. 
A typical result of the theory is that the set  $\Sigma(\Lambda) \cap \Delta_a$ does not change when the functional space on which $\Lambda$ is considered changes. 
We do not consider the {\it extension}  issue, which is however strongly connected to our approach, in the present work and we refer to the above mentioned articles for recent developments on that direction.

\smallskip
 Let us also mention that one expects that the spectrum $\Sigma(\Lambda)$ of $\Lambda$ is close to the spectrum $\Sigma(\BB)$ of $\BB$ if $\AA$ is ``small".
Such a  {\it ``small perturbation method" }  is an efficient tool in order to get some information on the spectrum of an operator $\Lambda$ ``in a perturbation 
regime". It 
has been developed  in  \cite{MR0054167,MR0089373,MR0445335,Kato}, and more recently in  \cite{MR1133342,MR1753014}. 
Again, we do not consider that  {\it ``small perturbation" }  issue here, but we refer to  \cite{MMcmp,MR2476685,AGGMMS,TristaniInelastic} 
where that kind of method is investigated in the same framework as the one of the abstract results of the present paper. 
 
\smallskip
Last, we are concerned with a positive operator $\Lambda$ defined on a Banach lattice $X$ and the associated semigroup $S_\Lambda$ as introduced by R.S. Phillips in \cite{MR0146675}. For a finite dimensional Banach space and a strictly positive matrix the Perron-Frobenius Theorem \cite{MR1511438,Frobenius} states that the eigenvalue with largest real part is unique, real and simple. 
Or in other words, there exists $a^{**} \in \R$ such that $\Sigma(\Lambda) \cap \Delta_{a^{**}} = \{ \lambda \}$ with $\lambda \in \R$ a simple eigenvalue. 
In an infinite dimensional Banach space the Krein-Rutman  Theorem \cite{MR0027128} establishes the same  result for a class of  Banach lattices and under
convenient strict positivity and compactness assumptions on $\Lambda$. The Krein-Rutman  Theorem is then extended to broader classes of Banach lattices and broader classes of operator in many subsequent articles, see e.g. \cite{MR618205,MR617977,AGLMNS}. We present a very natural and general version of the abstract Krein-Rutman Theorem on a Banach lattice assuming that, additionally to the above splitting structure, the operator $\Lambda$ satisfies a weak and a strong maximum principle.  Our result improves the known versions of Krein-Rutman  Theorem in particular because  we weaken the required compactness assumption made on (the resolvent of) $\Lambda$. Moreover, our result is quite elementary and self-contained.

\smallskip
 Let us stress again that our approach is very similar to the ``{\it extension } of the functional space of spectral analysis of operators"  and that
 our starting point is the work by Mouhot \cite{Mcmp} where \eqref{eq:RLRB-RLARB}
% a stopped Dyson series  \eqref{eq:stopDyson} or \eqref{eq:stopDysonBIS} with rank $n=N=1$
 is used in order to prove an enlargement of the operator spectral gap for the 
space homogeneous linearized Boltzmann equation. 
Because of the self-adjointness structure of the space homogeneous linearized Boltzmann equation one can conclude thanks to a classical
argument (namely the operator $\BB$ is sectorial and the last term  in \eqref{eq:stopDysonBIS} with the choice $N=1$ converges,  
 see \cite{Pazy} or \cite[Corollary IV.3.2 \& Lemma V.1.9]{EngelNagel}).  
Our approach is  in fact reminiscent of the huge number of works on the spectral analysis of operators which  take advantage of a  splitting structure 
\eqref{eq:L=A+B}
% $\Lambda= \AA + \BB$ 
and then consider $\Lambda$ as a (compact, small) perturbation of $\BB$. 
 One of the main differences with the classical approach introduced by Hilbert and Weyl is that when one usually makes the
decomposition
\beqn\label{eq:SplitA0B0}
\Lambda = \AA_0 + \BB_0
\eeqn
where $\BB_0$ is dissipative and $\AA_0$ is $\BB_0$-compact, we make the additional splitting  $\AA_0 = \AA + \AA^c$ with $\AA$ ``smooth" and $\AA^c$ small  (it is the usual way to prove that $\AA_0$ is $\BB_0$-compact) and we write  
\beqn\label{eq:SplitAB}
\Lambda = \AA + \BB, \quad \BB := \AA^c + \BB_0.
\eeqn
In such a way, we get quite better ``smoothing properties" on $\AA$ (with respect to $\AA_0$) without losing  to much of the ``dissipative property" of $\BB$ (with respect to $\BB_0$). That splitting makes possible to get quantitative estimates (with constructive constants) in some situations. 
Of course, the drawback of the method is that one has to find an appropriate splitting \eqref{eq:L=A+B} for the operator as well as an appropriate space $X$ for which one is able to prove the estimates \eqref{eq:hypH1},  \eqref{eq:hypH3-2} and \eqref{eq:hypH2-3}. However, the efficiency of the method is  attested by the fact that it has been successfully used for several evolution PDEs such as the 
space homogenous and space nonhomogeneous elastic Boltzmann equations in \cite{Mcmp,GMM}, 
the space homogenous and space nonhomogeneous inelastic Boltzmann equations in \cite{MMcmp,TristaniInelastic}, 
some Fokker-Planck type equations in \cite{GMM,AGGMMS,MM*,EM-KS1}, the Landau equation in   \cite{KleberLandau}
and the growth-fragmentation equation in \cite{MR2832638,MR2821681} and in the present paper.   

\medskip

\subsection{Growth-fragmentation equations}  
The second aim of the paper is to establish the long-time asymptotics for 
the solutions of some growth-fragmentation equations  as a motivation, or an illustration, of the abstract theory developed in parallel. 
We then consider the  growth-fragmentation equation 
\beqn\label{eq:croisfrag}
\partial_t f   = \Lambda f := \DD f +  \FF f  \quad \hbox{in} \quad (0,\infty) \times (0,\infty). 
\eeqn
%for cell division/self-similar fragmentation equation
Here  $f = f (t,x) \ge 0$ stands for the  {\it number density of cells }(or {\it particles, polymers, organisms, individuals}), where  $t \ge 0$ is the {\it time }variable
and  $x \in (0,\infty)$ is the {\it size} (or {\it mass, age}) variable. In equation \eqref{eq:croisfrag} the {\it growth} operator is given by  
 \beqn\label{eq:OpDD}
  (\DD f)(x) := -  \partial_x (\tau(x) f (x) )- \nu(x) \, f (x)
 \eeqn 
 where the (continuous) function $\tau : [0,\infty) \to \R$  is the drift speed (or {\it growth rate}),
  and we will choose $\tau (x) = 1$ or $\tau (x) = x$ in the sequel, 
% (see \cite{MR2652618,MR3030711} and the references therein where more general rates are considered) 
and  the function $\nu : [0,\infty) \to [0,\infty)$ is a damping rate.
The drift and damping term $\DD$ models the growth (for particles and cells) or the aging 
(for individuals) and the death which can be  represented by the scheme
$$
\{ x \} \,\,\, \mathop{\longrightarrow}^{e^{-\nu(x)}} \,\,\, \{ x + \tau(x) \, dx \}.
 $$
On the other hand, the  {\it fragmentation} operator $\FF$  is defined by
 \beqn\label{eq:OpFF}
(\FF f)(x) := \int_x^\infty k(y,x) f(y) dy - K(x) f(x)
\eeqn
and the fragmentation kernel $k$ is related to the total rate of fragmentation $K$ by
\beqn\label{eq:K=intk} 
K(x) = \int_0^x k(x,y) \, {y \over x} \, dy.
\eeqn
The fragmentation operator $\FF$ models the division (breakage) of a single {\it mother particle} of size $x$ into two or more pieces ({\it daughter particles, offspring}) of size $x_i \geq 0$, or in other words, models  the event
\beqn\label{xtoxk}
\{ x \} \,\,\, \mathop{\longrightarrow}^k \,\,\, \{x_1 \} + .... + \{ x_i \} + ... \, ,
\eeqn
in such a way that the mass is conserved 
$$
x = \sum_i x_i, \quad 0 \le x_i \le x.
$$
It is worth emphasizing that the above mass conservation at the microscopic level is rendered by  equation \eqref{eq:K=intk} at the statistical level. 

In order to simplify the presentation we will only consider situations where the size repartition of offspring is invariant by the size scaling of the
mother particle, or more precisely that there exists a function (or abusing notation, a measure) $\wp : [0,1] \to \R_+$ such that 
\beqn\label{eq:k=Kwp}
k(x,y) = K(x) \, \kappa(x,y), \quad \kappa(x,y) = \wp(y/x)/x,
\eeqn
as well as, in order that the compatibility relation \eqref{eq:K=intk} holds, 
\beqn\label{eq:intwp=1}
\int_0^1 z \, \wp (dz) = 1. 
\eeqn
Here $\kappa(x,.)$ represents the probability density of the distribution of daughter particles resulting of the breakage of a mother particle of
size $x \in (0,\infty)$ and the assumption \eqref{eq:k=Kwp} means that this probability density is invariant by scaling of the size. 
As a first example, when a mother particle of size $x$ breaks into two pieces of exact size $\sigma \, x$ and $(1-\sigma) x$, $\sigma \in (0,1)$, the associated kernel is given by  
\beqn \label{ex:celldiv2Gal}
\kappa(x,y) =   \delta_{y=\sigma x}(dy) + \delta_{y=(1-\sigma) x} (dy)
=  {1 \over \sigma} \delta_{x=\frac y\sigma } (dx)+ \frac 1{1-\sigma} \delta_{x=\frac  y{1-\sigma}} (dx), 
\eeqn
or equivalently $\wp = \delta_\sigma + \delta_{1-\sigma}$. In the sequel, we will also consider the case when $\wp$ is a smooth function. 
In any cases, we define 
 \beqn\label{def:z0}
z_0 := \inf \hbox{supp} \, \wp \in [0,1).
\eeqn

The evolution equation  \eqref{eq:croisfrag} is complemented with an initial condition 
\beqn\label{eq:GFt=0}
f(0,.) = f_0    \quad \hbox{in} \quad (0,\infty), 
\eeqn
and a boundary condition that we will discuss for each example presented below.

 Instead of trying to  analyze the most general growth-fragmentation equation, we will focus our study on 
 some particular but relevant classes of models, namely the cell division equation with constant growth rate, the self-similar fragmentation equation 
 and the age structured population equation.

\subsubsection{\bf Example  1. Equal mitosis equation. }  
%\medskip\noindent
%{\bf Example  1. Equal mitosis equation. }
%(or cell division) equation. } 
We consider a population of cells which divide through a binary fragmentation mechanism with equal size offspring, 
%(for the sake of simplicity) which  
grow at constant rate and  are not damped. 
The resulting evolution equation is the {\it equal mitosis equation} which is associated to the operator $\Lambda = \DD + \FF$, where 
$\DD$ is defined by   \eqref{eq:OpDD} with the choice   $\tau=1$ and $\nu=0$, and where $\FF$ is  the equal mitosis operator defined by   \eqref{eq:OpFF} with the following   choice of fragmentation kernel  
\bear \label{ex:celldiv2}
  k(x,y) &=&  2 \, K(x) \,   \delta_{x/2}   (dy)
=  4 \,  K(x) \,  \delta_{2y } (dx).
\eear
Equivalently, $k$ is given by \eqref{eq:k=Kwp} with $\wp (dz) := 2 \, \delta_{z=1/2}$. 
The associated equal mitosis equation takes the form
\beqn \label{eq:mitose}
\frac{\partial}{\partial t} f(t,x) + \frac{\partial}{\partial x} f(t,x) + K(x) f(t,x) =4 K(2x) \, f(t,2x),
\eeqn
and it is complemented with the boundary condition 
\beqn\label{eq:BoundaryCondition}
 f(t,0) = 0.
\eeqn
As its name suggests, such an equation appears in the modeling of cell division when mitosis occurs   (see \cite{BellAnderson,SinkoS1971,MR860959,MR2024498} and the references therein for linear models as well as 
\cite{MR1072714,MR2024498} for more recent nonlinear models for tumor growth) but also appears in telecommunication systems to describe some internet protocols \cite{BaccelliMcDR}.  

We assume that the total fragmentation rate $K$ is a nonnegative $C^1$ function defined on $[0,\infty)$ which satisfies the positivity assumption
\beqn\label{eq:hypKmitose1}
\exists \, x_0 \ge 0, \quad K(x) = 0 \,\,\, \forall x < x_0, \quad K(x) > 0 \,\,\, \forall \, x > x_0,
\eeqn
as well as the growth assumption
\beqn\label{eq:hypKmitose2}
K_0 \, x^\gamma \, {\bf 1}_{x \ge x_1} \le K(x) \le K_1 \, \max(1,x^\gamma) ,
\eeqn
for  $\gamma \ge 0$, $x_0 \le x_1 < \infty$, $0 < K_0 \le K_1 < \infty$.

There is no conservation law for the equal mitosis equation. However, by solving the dual eigenvalue-eigenfunction  problem
\beqn\label{eq:VP*pb}
\DD^* \phi + \FF^* \phi = \lambda \phi, \quad \lambda \in \R, \quad \phi \ge 0, \,\, \phi\not\equiv0, 
\eeqn
 one immediately observes that any solution $f$ to the equal mitosis equation \eqref{eq:mitose} satisfies
$$
\int_0^\infty f(t,x) \, \phi(x) \, dx = e^{\lambda t} \int_0^\infty f_0 \, \phi(x) \, dx.
$$
The first eigenvalue $\lambda$ corresponds to an exponential growth rate of (some average quantity of) the solution. 
In an ecology context $\lambda$ is often called the {\it Malthus parameter } or the {\it fitness } of the cells/organisms population. 
In order to go further in the analysis of the dynamics, one can solve the primal eigenvalue-eigenfunction problem
\beqn\label{eq:VPpb}
\DD f_\infty + \FF f_\infty = \lambda f_\infty, \quad f_\infty \ge 0, \,\, f_\infty\not\equiv0, 
\eeqn
which provides a Malthusian profile $f_\infty$, and then define the  remarkable (in the sense  that it is a separated variables function) solution  $f(t,x) := e^{\lambda t} \, f_\infty(x)$  
to the equal mitosis equation \eqref{eq:mitose}. It is expected that \eqref{eq:VPpb} captures the main features of the model  or, 
more precisely, that 
\beqn\label{eq:GFtofinfty}
f(t,x) e^{-\lambda t} = f_\infty(x) + o(1) \quad\hbox{as}\quad t\to\infty. 
\eeqn
We refer to  \cite{MR745853,MR860959,MR2114128,MMP,MR2270822,MR2536450} (and the references therein) for   results about the existence of solutions to the primal and dual eigenvalue problems \eqref{eq:VP*pb} and \eqref{eq:VPpb} as well as for results on the asymptotic convergence \eqref{eq:GFtofinfty} without rate  or with exponential rate.  

 \subsubsection{\bf  A variant of example  1. Smooth cell-division equation. } We modify  the previous equal mitosis model by considering the case of a general fragmentation operator \eqref{eq:OpFF} where the total fragmentation rate $K$ still satisfies \eqref{eq:hypKmitose1} and \eqref{eq:hypKmitose2} 
and where, however, we restrict ourself to the case of a smooth size offspring distribution $\wp$. More precisely, we assume
\beqn\label{eq:wpsmooth}
\wp'_0 := \int_0^1 |\wp'(z)| \, dz < \infty.
\eeqn
We will sometimes make the additional  strong positivity assumption 
 \beqn\label{eq:wp>0}
 \wp(z) \ge \wp_* \,\, \, \forall \, z \in (0,1), \quad \wp_* > 0, 
 \eeqn
or assume the additional monotony condition and the constant number of offspring condition
 \beqn\label{eq:intkappa<0}
{\partial \over \partial y} \int_0^x \kappa (y,z) \, dz \le 0, 
\quad
n_F := \int_0^x \kappa (x,z) \, dz > 1, 
\quad \forall \, x,y \in \R_+.
\eeqn
Observe that this monotony condition is fulfilled for the equal and unequal mitosis kernel \eqref{ex:celldiv2Gal} and for the smooth distribution of offspring functions
$\wp(z) := c_\theta \, z^\theta$, $\theta > -1$, see \cite{MR2536450}. 

The smooth cell-division equation reads 
\beqn \label{eq:celldiv}
\frac{\partial}{\partial t} f(t,x) + \frac{\partial}{\partial x} f(t,x) = (\FF f(t,.))(x) ,
\eeqn
and it is  complemented again  with the boundary condition \eqref{eq:BoundaryCondition} and the initial condition \eqref{eq:GFt=0}. 
We call the resulting model as the {\it smooth cell-division equation}.

The general fragmentation operator is used in physics in order to model the dynamics of cluster breakage and it is often associated in that context to the coagulation operator which models 
the opposite agglomeration mechanism (see \cite{MR0087880}).  
It appears later associated with the drift operator $\DD = - \tau(x) \partial_x + \tau'(x)$ under the name of ``cell population balance model" (see  \cite{Fredrickson,MR764484}) in a chemical or biological context.  The general fragmentation operator is used in order to take into account unequal cell-division according to experimental evidence \cite{KochS1962,PainterM1968}.
% {\Red ICI population balance equation for the modeling of chemostat by : \cite{Fredrickson} modelling  of chemostat by cell population balance model} 
In recent years, the above smooth cell-division equation \eqref{eq:celldiv} has also appeared in many articles on the modeling of proteins \cite{MR2279568,MR2172205,MR2489341,MR2605707}.  

Concerning the mathematical analysis of the {\it smooth cell-division equation}, and in particular the long-time behaviour of solutions, a similar picture as for the 
equal mitosis equation is expected and some (at least partial) results have been obtained in  \cite{MR2065377,MMP,MR2250122,MR2536450,MR2652618,MR2832638,MR2935368,MR3030711}. 

With the above notation and for later references, for both equal mitosis and smooth cell-division equations, we introduce the critical exponent $\alpha^* \ge 1$ uniquely implicitly defined by the equation
\beqn\label{eq:defalpha*}
\wp_{\alpha^*} = K_0/K_1 \in (0,1], \quad  \wp_\alpha := \int_0^1 z^{\alpha} \, \wp(dz).
\eeqn

\subsubsection{\bf Example  2. Self-similar fragmentation equation. } We consider now a fragmentation rate associated to a power law total fragmentation rate %and then to a fragmentation rate
\beqn \label{eq:hypSSfrag}
 k(x,y) = K(x) \, \wp(y/x)/x, \quad K(x) = x^\gamma, \quad \gamma > 0,
\eeqn
where $\wp$ is a continuous  function satisfying \eqref{eq:wpsmooth}. 
The pure fragmentation model is then obtained for $\tau=\nu  \equiv 0$ in \eqref{eq:croisfrag} and therefore reads
$$
\partial_t g = \tilde \FF \, g , 
$$
where $\tilde \FF$ is a temporary notation for the fragmentation operator associated to the kernel $\tilde k := k/\gamma$.  By analogy with the probabilistic name for the associated Markov process, see e.g. \cite{MR2017852}, we call that model the {\it ``self-similar fragmentation equation"}. 
This equation arises in physics to describe fragmentation processes.  We refer to \cite{Filippov,MR814641,BCP} for the first study and the physics motivations, to 
\cite{MR2017852} and the references therein for a probabilistic approach. 

For this equation the only steady states are the Dirac masses, namely $x \, g(t,x) = \rho \, \delta_{x=0}$.
On the other hand, if $g$ is a solution to the pure fragmentation equation, we may introduce the rescaled density $f$ defined by
\beqn
\label{ftogF}
f(t, x)=e^{-2t}g\left(e^{\gamma\, t}-1, xe^{-t}\right), \eeqn
which is a solution to the fragmentation equation in self-similar variables (see for instance \cite{EMRR}) \beqn
\label{eq:fragSSV}
\frac{\partial}{\partial t} f =  \Lambda f = \FF \, f - x \,  \frac{\partial}{\partial x}f - 2 f. 
\eeqn
This is a mass preserving equation with no detailed balance condition. However, defining the associated adjoint operator $\Lambda^*$, one can show
that 
$$
\Lambda f_\infty = 0, \quad \Lambda^* \phi = 0, 
$$
for some positive function $f_\infty$  and for $\phi(x) = x$. 
Equation \eqref{eq:fragSSV} falls into the class of growth-fragmentation equation with first eigenvalue $\lambda = 0$ (because of the mass conservation). As for the cell-division equation one expects that the remarkable (self-similar profile)  solution $f_\infty$ is attractive and that \eqref{eq:GFtofinfty} holds again. 
Existence of the self-similar profile $f_\infty$ and convergence (without rate) to this one have been established in  \cite{EMRR,MMP}, while 
a rate of convergence for $\gamma \in (0,2)$ has been proved in \cite{MR2832638,MR2821681}.

%%%

\subsubsection{\bf  Example  3.  Age structured population equation. } 
We consider an age structured population of individuals which age, die and give birth, and which is described by the density $f(t,x)$ of individuals with age $x \ge 0$ at time $t \ge 0$.  The very popular associated evolution PDE is  
\beqn
\label{eq:renewal}
\frac{\partial}{\partial t} f  + \frac{\partial}{\partial x} (\tau(x) f )  = - \nu(x) \, f , \quad f(t,x=0)=\int_0^\infty K(y) f(t,y) dy,
\eeqn
and it is commonly attributed to A. McKendrik~\cite{McKendrik} and H. von Foerster~\cite{vonFoerster} (although the dynamics of age structured population 
has been anteriorly developed by  A. Lotka and F. Sharpe~\cite{SharpeLotka} and before by L. Euler, as well as by P.H. Leslie \cite{MR0015760} in a discrete time and age framework). In the sequel we call that model as the {\it age structured population equation}. Let us also mention that nonlinear versions of equation  \eqref{eq:renewal} which take into account possibly overcrowding effect  can be found in the work of 
Gurtin and MacCamy \cite{MR527566}.

In equation \eqref{eq:renewal} the function $K$ corresponds to the birth rate, the function $\nu$ to the death rate and the function $\tau$ to the aging rate (so that $\tau \equiv 1$). Notice that the age structured population equation can be seen as a particular example of the growth-fragmentation equation (\ref{eq:croisfrag})--(\ref{eq:OpFF}) making the following 
choice for $k$ :
\beqn \label{ex:celldiv1}
k(x,y)= K(x) \; [ \delta(y=x) +\delta(y=0)], 
\eeqn
which corresponds to the limit case $\sigma = 0$ in  \eqref{ex:celldiv2Gal}. 
In order to simplify the presentation we make the assumptions 
\beqn \label{eq:renewal1}
 \tau = \nu = 1, \quad 0 \le K \in C^1_b(\R_+) \cap L^1(\R_+), \quad \| K \|_{L^1} > 1.
\eeqn
The expected long-time behaviour of solutions to the age structured population equation is the same as the one described for a general growth-fragmentation equation, in particular the long-time convergence \eqref{eq:GFtofinfty} with exponential rate is known to hold. Here, the mathematical analysis is greatly simplified by remarking that the offspring number satisfies a (Volterra) delay equation (the so-called  renewal equation) which in turn can be handled through the direct Laplace transform,  as first shown by W. Feller \cite{MR0005419} and developed later in \cite{MR0210154,MR860959,MR772205,Iannelli}.  Let us also mention that the long-time convergence \eqref{eq:GFtofinfty} can also be obtained by   entropy method \cite{MR1946722,MMP}. 

%Our approach here is more abstract and can be seen as some Laplace transform method at the abstract level (...).

%\cite{MR860960,MR860962}

\subsubsection{\bf Main result. }\label{sss:MainResult} Let us introduce the functional spaces in which we will work. 
For any exponent $p \in [1,\infty]$ and any nonnegative weight function $\xi$, we denote by $L^p(\xi)$ the Lebesgue space $L^p(\R_+;\xi \, dx)$ or $L^p(\R;\xi \, dx)$ associated to the norm 
$$
\|u\|_{L^p(\xi)} :=  \|u \, \xi\|_{L^p},
$$
and we simply use the shorthand $L^p_\alpha$ for the choice $\xi(x) := \langle x \rangle^\alpha$, $\alpha \in \R$, $\langle x \rangle^2 := 1 + |x|^2$, 
% for any $x \in \R$ 
as well as the shorthand $\dot L^p_\alpha$ for the choice $\xi(x) :=|x|^\alpha$, $\alpha \in \R$.

\begin{theo}\label{theo:Frag} Consider the growth-fragmentation equations with the corresponding structure assumption and boundedness of coefficients as presented in the previous sections and define the functional space $X$ as follows: 

\begin{itemize}
 
\item[{\bf (1)}] {\it Cell-division equation: } take $X=L^1_\alpha$, $\alpha > \alpha^*$, where $\alpha^* \ge 1$ is defined in \eqref{eq:defalpha*}; 
 
\item[{\bf (2)}] Self-similar fragmentation equation: take  $X = \dot L^1_\alpha \cap \dot L^1_\beta$, $0 \le \alpha < 1 < \beta$; 

\item[{\bf (3)}] Age structured population equation: take  $X = L^1$. 
 
\end{itemize}
 
\smallskip
There exists a unique couple $(\lambda,f_\infty)$, with $\lambda \in \R$ and $f_\infty \in X$,  
solution to the stationary equation 
\beqn\label{eq:PrimalEVpb}
\FF f_\infty -   \partial_x (\tau f ) - \nu f_\infty = \lambda \, f_\infty, \quad f_\infty \ge 0, 
\quad \| f_\infty \|_{X} = 1.
\eeqn

There exists $a^{**} < \lambda$ and for any $a > a^{**}$ there exists $C_a$ such that for any $f_0 \in X$,  the associated solution $f(t) = e^{\Lambda t} f_0$  satisfies 
\beqn\label{eq:RateGrFr}
\| f(t) - e^{\lambda t} \, \Pi_{\Lambda,\lambda} f  \|_{X} \le C_a \, e^{at} \, \| f_0 -   \Pi_{\Lambda,\lambda} f_0  \|_{X} ,
\eeqn
where $\Pi_{\Lambda,\lambda}$ is the projection on the one dimensional space spanned by the remarkable solution $f_\infty$. It is   given by 
$$
\Pi_{\Lambda,\lambda} h  = \langle \phi, h \rangle \, f_\infty
$$
where $\phi \in X'$ is the unique positive and normalized solution to the dual first eigenvalue problem
\beqn\label{eq:DualEVpb}
\FF^*\phi +  \tau \partial_x \phi- \nu \phi = \lambda \, \phi , \quad \phi \ge 0, \quad \langle \phi, f_\infty \rangle = 1.
\eeqn
Moreover, an explicit bound on the spectral gap $\lambda-a^{**}$ is available for
\begin{itemize}
 
\item[{\bf (i)}] the cell-division equation with constant total fragmentation rate $K \equiv K_0$ on $(0,\infty)$, $K_0 > 0$, 
and a fragmentation kernel which satisfies the monotony condition and constant number of offspring condition \eqref{eq:intkappa<0};

\item[{\bf (ii)}] the self-similar fragmentation equation with smooth and positive offspring distribution in the sense that \eqref{eq:wpsmooth} and \eqref{eq:wp>0} hold.
 
\end{itemize}

\end{theo}

Let us make some comments about the above result. 

\smallskip
Theorem~\ref{theo:Frag} generalizes, improves and unifies the results on the long-time asymptotic convergence with exponential rate  which were known only for particular cases of growth fragmentation equation, namely for the cell division model with (almost) constant total fragmentation rate and monotonous offspring distribution  in  \cite{MR2114128,MR2536450} and for the self-similar fragmentation equation and the cell division equation restricted to 
 smooth and  positive fragmentation rate and  total fragmentation rate which increases at most  quadratically in \cite{MR2832638,MR2821681}. The rate of convergence \eqref{eq:RateGrFr} is proved under similar hypothesizes as those made in  \cite{MMP,EMRR}, but in \cite{MMP,EMRR} the  convergence is  established without rate. 
 It has been established  in   \cite{MR2114128,MR2536450} a similar  $L^1$-norm decay as in \eqref{eq:RateGrFr},  when however the initial datum is bounded in the sense of a  (strange and) stronger norm than the $L^1$-norm. It was also conjectured that the additional strong boundedness assumption on the initial datum is necessary in order to get an exponential rate of convergence.  Theorem~\ref{theo:Frag} disproves that conjecture in the sense that the norm involved in both sides of estimate  \eqref{eq:RateGrFr}  is the same. %, and it is  only a (weighted) $L^1$-norm. 
%shows that it is not the case and that the additional stronger norm boundedness assumption made in   \cite{MR2114128,MR2536450} can be removed
%(although the norm used in the present paper is different  from (and stronger than) those used in    \cite{MR2114128,MR2536450} the fundamental point is %that the norm measuring the distance between the solution and its limit is the same for positive time and at initial time).  
Let us emphasize that we do not claim that Theorem~\ref{theo:Frag} is new for the age structured population equation. However, we want to stress here that our proof of the convergence  \eqref{eq:RateGrFr} is similar for all these growth-fragmentation equations while the previous available proofs (of convergence results with rate) were
very different for the three  subclasses of models. It is also likely that our approach can be generalized to larger classes of growth operator and of fragmentation kernel such as considered in \cite{MR2250122,MR2652618,MR2935368,MR3030711}. However, for the sake of simplicity, we have not followed that line of research here.  
It is finally worth noticing that our result excludes the two ``degenerate equations" which are the equal self-similar fragmentation equation associated to equal mitosis offspring distribution and the age structured population equation associated to deterministic birth rate $K(z) := K \, \delta_{z=L}$. 
 
 \smallskip

Let us make some comments about the different methods of proof which may be based on linear tools (Laplace transform, Eigenvalue problem, suitable weak distance, semigroup theory) and nonlinear tools (existence of self-similar profile by fixed point theorems, GRE and E-DE methods). 

\smallskip

$\bullet$ {\sl Direct Laplace transform method. } For the  age structured population equation a direct Laplace transform  analysis can be performed at the level of the associated renewal equation and leads to an exact representation formula which in turn implies the rate of convergence \eqref{eq:RateGrFr} (see \cite{MR0005419,MR0210154,Iannelli}).

$\bullet$ {\sl PDE approach via compactness and GRE methods.} Convergence results (without rate) have been proved in \cite{MR1946722,EMRR,MMP} for a general class of 
growth-fragmentation equation which is basically the same class as considered in the present paper thanks to the use of the
so called {\it general relative entropy method}. More precisely, once the existence of a solution $(\lambda,f_\infty,\phi)$ to the primal eigenvalue problem \eqref{eq:PrimalEVpb} and dual eigenvalue problem \eqref{eq:DualEVpb} is established, we refer to \cite{MR2114128,MR2935368,MR2250122,MR2652618} where that issue is tackled, one can easily compute the evolution of the 
generalized relative entropy  $\JJ$ defined by  
$$
\JJ(f) := \int_0^\infty j(f/f_\infty) \, f_\infty \, \phi \, dx 
$$
for a  convex and non negative function $j : \R \to \R$ and for a generic solution $f=f(t)$ to the growth-fragmentation equation.
One can then show  the (at least formal) identity
\beqn\label{eq:GRE}
\JJ(f(t)) + \int_0^t \DD_\JJ(f(s)) \, ds = \JJ(f(0)) \qquad \forall \, t \ge 0, 
\eeqn
where $ \DD_\JJ \ge 0$ is the  associated generalized dissipation of relative entropy defined by 
$$
\DD_\JJ (f) := \int_0^\infty\!\! \int_0^\infty k(x_*,x) {\bf 1}_{x_*\ge x} \, [j(u) - j(u_*) - j'(u_*) \, (u-u_*) ]\, f_\infty \, \phi_* \, dxdx_*,
$$
with the notation $u := f/f_\infty$ and $h_* = h(x_*)$. Identity \eqref{eq:GRE} clearly implies that  $\JJ$ is a Liapunov functional for the growth fragmentation dynamics which in turn implies  the long-time convergence   $f(t) \to f_\infty$ as $t\to\infty$ (without rate) under positivity assumption on the kernel $k$ or the associated semigroup.

$\bullet$ {\sl Suitable weak distance. } In order to circumvent the possible weak information given by $\DD_\JJ$ in the case of cell-division models when $k$ is not positive 
(for the equal mitosis model for instance) an exponential rate of convergence for an (almost) constant total rate of fragmentation 
has been established by Perthame and co-authors in   \cite{MR2114128,MR2536450} by exhibiting a suitable weak distance.

$\bullet$  {\sl Entropy - dissipation of entropy (E-DE) method. }  In the case of the strong positivity assumption \eqref{eq:wp>0} an entropy - dissipation of entropy method has been implemented in  \cite{MR2832638,MR3030711} where the inequality 
$\DD_\JJ  \ge c \, \JJ$ is proved for $j(s) = (s-1)^2$. The Gronwall lemma then straightforwardly implies a rate of convergence in a weighted Lebesgue $L^2$ framework. 

$\bullet$ {\sl Extension semigroup method. } For a general fragmentation kernel (including the smooth cell-division equation and the self-similar fragmentation equation) the  {\it enlargement of semigroup spectral analysis} has been used in  \cite{MR2832638,MR2821681} in order to extend to a $L^1$ framework some of the convergence (with rate) results proved in  \cite{MR2832638}  in a narrow 
%strongly confinent 
weighted $L^2$ framework thanks to the above E-DE method.

$\bullet$ {\sl Markov semigroup method. } Markov semigroups techniques have been widely used since the first papers by Metz, Diekmann and Webb on the age structured equation and the equal mitosis equation in a bounded size setting \cite{MR772205,MR860959}, see also \cite{EngelNagel} and the references therein for more recent papers. That approach has also been applied to more general growth-fragmentation equations (still in a bounded interval setting) in  \cite{RudnickiP2000} (see also \cite{MR2915575} and the references therein for recent developments) and to other linear Boltzmann models (related to neuron transport theory) by 
Mokhtar-Kharroubi  in \cite{MR1151539,MR1612403,MR2237499}. The Markov semigroup method consists in proving that the spectral bound is an algebraic simple eigenvalue associated to a positive eigenfunction and a positive dual eigenfunction and the only eigenvalue with maximum real part. It fundamentally uses the positivity structure of the equation and some compactness argument (at the level of the resolvent of the operator or at the level of the terms involved in the Dyson-Phillips series). 
This approach classically provides a convergence result toward the first eigenfunction, but does not provide any rate of convergence in general (see however the recent work \cite{MR3192457} where a rate of convergence is established).

Let us conclude that our method in the present paper is clearly a Markov semigroup approach. Our approach is then completely linear, very accurate and still very general. The drawback is the use of the abstract semigroup framework
and some complex analysis tool  (and in particular the use of the Laplace transform and the inverse Laplace transform) at the level of the abstract (functional space) associated evolution equation. 
 The main novelty comes from the fact that we are able to prove that the considered growth-fragmentation operator falls in the class of operators with splitting structure as described in the first part of the introduction and for which our abstract Krein-Rutman Theorem applies. In particular, in order to deduce the growth estimates on the semigroup from the spectral knowledge on the generator, we use the iterated Duhamel formula (instead of the Dyson-Phillips series) as a consequence of our more suitable splitting \eqref{eq:SplitAB} (instead of \eqref{eq:SplitA0B0}). That more suitable splitting and our abstract semigroup spectral analysis theory make possible to establish some rate of convergence where the usual arguments lead to mere convergence results (without rate).

\medskip
The outline of the paper is as follows. In the next section, we establish the partial spectral mapping theorem in an abstract framework.
In Section~\ref{sec:Weyl}, we establish two versions of the Weyl's Theorem in an abstract framework  and we then  verify that they apply to the growth-fragmentation equations in Section~\ref{sec:SGWgrowthfrag}. Section~\ref{sec:KRthAbstract} is devoted to the statement and proof of an abstract version of the Krein-Rutman  Theorem. In the last section we apply the  Krein-Rutman theory  to the growth-fragmentation equation which ends the proof of Theorem~\ref{theo:Frag}.

\medskip
{\bf Acknowledgment. }  
The first author would like to thank O. Kavian for his encouragement  to  write the abstract part of the present paper as well as 
for his comments on a first proof  of Theorem~\ref{theo:SpectralMapping} which had led to many simplifications on the arguments. 
The first author also would like to thank C.~Mouhot for the many discussions of the spectral analysis of semigroups issue. 
We also would like to acknowledge %M. Doumic and 
O. Diekmann for having pointing out some interesting references related to our work.  The authors thank  the referees for helpful remarks and suggestions.
The authors gratefully acknowledge the support of the STAB ANR project (ANR-12-BS01-0019).

\bigskip
\section{Spectral mapping  for semigroup generators}
\label{sec:SpectralMap}
\setcounter{equation}{0}
\setcounter{theo}{0}

 %--------------------------------------------------------------------------------------------------------------------------------------------------------------

\subsection{Notation and definitions}
\label{sec:notation-definitions}

There are many textbooks addressing the theory of semigroups since the seminal books by Hille and Philipps \cite{MR0025077,MR0089373} 
among them the ones by Kato \cite{Kato}, Davies \cite{ MR591851}, Pazy \cite{Pazy},  Arendt et al \cite{AGLMNS} and more recently by Engel and Nagel  \cite{EngelNagel} to which we refer for more details. In this section we summarize some basic definitions and facts on the analysis of operators in a abstract and
general Banach space picked up from above mentioned books as well as from the recent articles \cite{GMM,moiED}. It is worth mentioning that we adopt the sign convention of Kato \cite{Kato} on the resolvent operator (which is maybe opposite to the most widespread convention), see \eqref{def:Resolv}. 
For a given real number $a \in \R$, we define the half complex plane
$$
\Delta_a :=  \{ z \in \C, \, \Re e \, z > a \}.
$$
 For some given Banach spaces $({X_1},\|\cdot \|_{X_1})$ and $({X_2},\| \cdot
\|_{X_2})$ we denote by $\mathscr{B}({X_1}, {X_2})$ the space of bounded linear
operators from ${X_1}$ to ${X_2}$ and we denote by   $\| \cdot
\|_{\mathscr{B}({X_1},{X_2})}$ or $\| \cdot \|_{{X_1} \to {X_2}}$ the associated norm
operator. We write $\mathscr{B}(X_1) = \mathscr{B}(X_1,X_1)$.
We denote by $\mathscr{C}(X_1,{X_2})$ the space of closed unbounded (and thus possibly bounded) linear
operators from $X_1$ to ${X_2}$ with dense domain, and $\mathscr{C}(X_1)=
\mathscr{C}(X_1,X_1)$. 
We denote by $\KKK(X_1,{X_2})$ the space of compact linear
operators from $X_1$ to ${X_2}$ and again $\mathscr{K}(X_1)=
\mathscr{K}(X_1,X_1)$. 

For a given Banach space  $(X,\|\cdot \|_X)$, we denote by $\mathscr{G}(X)$ the space of generators of a $C_0$-semigroup. 
For $\Lambda \in \GGG(X)$ we denote by $S_\Lambda(t) = e^{t \, \Lambda }$, $t \ge
0$, the associated semigroup, by $\textrm{D}(\Lambda)$ its domain, by
$\textrm{N}(\Lambda)$ its null space, by
$$
M(\Lambda) = \bigcup_{\alpha \in \N^*} N(\Lambda^\alpha)
$$
its algebraic null space, and by $\mbox{R}(\Lambda)$ its range. For any given integer $k \ge 1$, we define  $\textrm{D}(\Lambda^k)$
 the  Banach space associated with the norm
$$
\|f \|_{\textrm{D}(\Lambda^k)} =  \sum_{j=0}^k \| \Lambda^j f \|_X.
$$
For two operators $\AA, \, \BB \in \CCC(X)$, we say that $\AA$ is $\BB$-bounded if $\AA \in \BBB(D(\BB),X)$ or, in other words, if
there exists a constant $C \in (0,\infty)$ such that 
$$
\forall \, f \in X, \qquad \| \AA f \|_X \le C (\| f \|_X + \| \BB f \|_X).
$$
Of course a bounded operator $\AA$ is always $\BB$-bounded (whatever is $\BB \in \CCC(X)$). 
%and then for any real number $\zeta > 0$, we also denote $\textrm{D}(\Lambda^\zeta)$ the space defined 
%by interpolation from the family of spaces $\textrm{D}(\Lambda^k)$, $k \in \N$. 

 For $\Lambda \in \GGG(X)$, we denote by $\Sigma(\Lambda)$ its spectrum, so that for any $z \in \C
\backslash \Sigma(\Lambda)$ the operator $\Lambda - z$ is invertible
and the resolvent operator
\beqn\label{def:Resolv}
\RR_\Lambda(z) := (\Lambda -z)^{-1}
\eeqn
is well-defined, belongs to $\mathscr{B}(X)$ and has range equal to
$D(\Lambda)$. We then define the spectral bound $s(\Lambda) \in  \R \cup \{-\infty \}$ by 
$$
s(\Lambda) := \sup \{ \Re e \, \xi; \,\, \xi \in \Sigma(\Lambda) \} 
$$
and the growth bound $\omega(\Lambda) \in \R \cup \{-\infty \}$ by 
$$
\omega(\Lambda) := \inf \{ b \in \R; \,\, \exists \, M_b \,\, \hbox{s.t.} \,\,  \| S_\Lambda(t) \|_{\BBB(X)} \le M_b \, e^{bt} \ \forall \, t \ge 0 \} , 
$$
and we recall that $s(\Lambda) \le \omega(\Lambda)$ as a consequence of Hille's identity \cite{MR0006610}:  for any $\xi \in \Delta_{\omega(\Lambda)}$ there holds
\beqn\label{eq:RL=intSt}
- R_\Lambda(\xi) = \int_0^\infty S_\Lambda(t) \, e^{-\xi t} \, dt,
\eeqn
where the RHS integral normally converges.

  We say that $\Lambda$ is \emph{$a$-hypo-dissipative} on $X$
   if there exists some norm $\Nt  \cdot \Nt_X$ on $X$   equivalent to the initial norm $\| \cdot \|_X$ such that
\begin{equation}\label{eq:def-dissipative} 
  \forall \, f \in  D(\Lambda) , \,\, \exists \, \varphi \in F(f) \ \mbox{
    s.t. } \  \Re e
  \, \langle \varphi, (\Lambda-a) \, f \rangle \le 0, 
\end{equation} 
where $\langle \cdot , \cdot \rangle$ is the duality bracket for the duality in $X$ and $X'$
and $F(f)  \subset X'$ is the dual set of $f$
defined by
$$
F(f) = F_{\Nt \cdot \Nt} (f) := \left\{\varphi \in X'; \,\, \langle \varphi , f \rangle = \Nt
  f \Nt_X^2 = \Nt \varphi \Nt^2_{X'} \right\}.
$$
We just say that $\Lambda$ is {\it hypo-dissipative} if $\Lambda$ is $a$-hypo-dissipative for some $a \in \R$. 
 From the Hille-Yosida Theorem it is clear that any generator $\Lambda \in \GGG(X)$ is an hypo-dissipative operator
 and that 
$$
\omega (\Lambda) := \inf \{ b \in \R; \,\, \Lambda  \,\, \hbox{is } b\hbox{-hypo-dissipative} \}.
$$ 

\smallskip
A spectral value $\xi \in \Sigma(\Lambda)$ is said to be isolated if
\[
\Sigma(\Lambda) \cap \left\{ z \in \C, \,\, |z - \xi| \le r \right\} =
\{ \xi \} \ \mbox{ for some } r >0.
\]
In the case when $\xi$ is an isolated spectral value we may define
 the spectral projector $\Pi_{\Lambda,\xi} \in \mathscr{B}(X)$ by the Dunford integral
\begin{equation}\label{def:SpectralProjection} 
\Pi_{\Lambda,\xi} := {i \over
  2\pi} \int_{ |z - \xi| = r' }\RR_\Lambda(z) \, dz, 
\end{equation}
with $0<r'<r$. Note that this definition is independent of the value
of $r'$ as the resolvent
\[
\C \setminus \Sigma(\Lambda) \to \mathscr{B}(X), \quad z \to
\RR_{\Lambda}(z)\] is holomorphic. It is well-known  that $\Pi_{\Lambda,\xi}^2=\Pi_{\Lambda,\xi}$, 
so that $\Pi_{\Lambda,\xi}$ is indeed a  projector, and its range 
$R(\Pi_{\Lambda,\xi})=\overline{M(\Lambda -\xi)}$ is the closure of
the algebraic eigenspace associated to $\xi$. 
 More generally, for any compact part of 
the spectrum of the form $\Gamma = \Delta_a \cap \Sigma(\Lambda)$ we may define the 
associated spectral projector $\Pi_{\Lambda,\Gamma}$ by 
\begin{equation}\label{def:SpectralProjectionGamma} 
\Pi_{\Lambda,\Gamma} := {i \over
  2\pi} \int_{\gamma }\RR_\Lambda(z) \, dz, 
\end{equation}
for any closed path $\gamma : [0,1] \to \Delta_a$ which makes one direct turn around $\Gamma$. 

%We can also define the
%family of operators corresponding to the projected semigroup:
%\[
%S_{\Lambda,\xi}(t) := -\frac{1}{2i \pi} \int_{|z-\xi|=r'} e^{zt}
%\RR_\Lambda (z) \, dz = \Pi_{\Lambda,\xi}  S_{\Lambda}(t) =
%S_\Lambda(t)  \Pi_{\Lambda,\xi}, \quad t >0.
%\]

We recall that $\xi \in \Sigma(\Lambda)$ is said to be an eigenvalue
if $N(\Lambda - \xi) \neq \{ 0 \}$.  
The range of the spectral projector is finite-dimensional if and only
if there exists $\alpha_0 \in \N^*$ such that
\[
N(\Lambda -\xi)^\alpha = N(\Lambda -\xi)^{\alpha_0}  \neq \{ 0 \} \ \mbox{ for
  any } \ \alpha \ge \alpha_0,
\]
and in such a case
\[
\overline{M(\Lambda-\xi)} = M(\Lambda - \xi) = N((\Lambda
-\xi)^{\alpha_0}) \quad \hbox{and} \quad N(\Lambda - \xi) \neq \{ 0 \}. 
\]
In that case, we say that $\xi$ is a discrete eigenvalue, written as
$\xi \in \Sigma_d(\Lambda)$.  For any $a \in \R$ such that
\[
\Sigma(\Lambda) \cap \Delta_{a } = \left\{ \xi_1, \dots, \xi_J\right\}
\]
where $\xi_1, \dots, \xi_J$ are distinct discrete eigenvalues, we define
without any risk of ambiguity 
\[
\Pi_{\Lambda,a} := \Pi_{\Lambda,\xi_1} + \dots \Pi_{\Lambda,\xi_J}.
\]
 
 For some given Banach spaces $X_1$, $X_2$, $X_3$ and some given
  functions
  \[
  S_1 \in L^1(\R_+; \BBB(X_1,X_2)) \ \mbox{ and } \ S_2 \in L^1(\R_+;
  \BBB(X_2,X_3)),
  \]
  we define the convolution $S_2 \ast S_1 \in L^1(\R_+; \BBB(X_1,X_3))$ by
  $$
  \forall \, t \ge 0, \quad (S_2 * S_1)(t) := \int_0^t S_2(s) \, S_1 (t-s) d s. 
  $$
  When $S_1=S_2$ and $X_1=X_2=X_3$, we define recursively $S^{(*1)} =
  S$ and $S^{(*\ell)} = S * S^{(*(\ell-1))} $ for any $\ell \ge 2$.

\smallskip
For a generator $L$ of a semigroup such that $\omega (L) < 0$ we define the fractional powers $L^{-\eta}$ and 
 $L^\eta$ for $\eta \in (0,1)$ by Dunford formulas \cite{MR0201985,MR0206716}, see also \cite[section II.5.c]{EngelNagel}, 
\beqn\label{eq:FracPower}
L^{-\eta} := 
%{e^{-2i\pi\eta} - 1 \over 2\pi i} 
c_{-\eta} \int_0^\infty \lambda^{-\eta} \, \RR_L(\lambda) \, d\lambda,
\quad
L^\eta :=  
%{\Gamma(1) \over \Gamma(\eta) \Gamma(1-\eta)} 
c_\eta\int_0^\infty \lambda^{\eta-1} \, L \, \RR_L(\lambda) \, d\lambda,
\eeqn 
for some constants $c_\eta,c_{-\eta} \in \C^*$. The operator $L^{-\eta}$ belongs to $\BBB(X)$ and, denoting $X_\eta := R(L^{-\eta})$, the operator $L^\eta$ is the unbounded operator with domain $D(L^\eta) = X_\eta$ and defined by $L^\eta = (L^{-\eta})^{-1}$. We also denote $X_0 = X$ and $X_1 = D(L)$. Moreover, introducing the $J$-method interpolation norm
$$
\| f \|_{\tilde X_\eta} := \inf \Bigl\{ \sup_{t > 0} \| t^{-\eta} \, J(t,g(t)) \|_X; \,\, g \, \hbox{such that} \, f = \int_0^\infty g(t) \, {dt \over t} \Bigr\}
$$
with $J(t,g) := \max(\| g \|_X, t \| L g \|_X )$ and the associated Banach space $\tilde X_\eta$ which corresponds to the interpolation space  $S(\infty,-\eta,X_0;\infty,1-\eta,X_1)$ 
of Lions and Peetre \cite{MR0165343} defined by 
$$
\tilde X_\eta := \{ f \in X; \,\, \| f \|_{\tilde X_\eta} < \infty \},
$$
the following inclusions
\beqn\label{eq:XalphaSubsetTildeXalpha}
X_\eta \subset \tilde X_\eta \subset X_{\eta'}
\eeqn
hold with continuous embedding for any $0 <  {\eta'} <  \eta <1$. Let us emphasize that the first inclusion follows from the second inclusion in \cite[Proposition 2.8]{MR0206716} and  \cite[Theorem 3.1]{MR0206716} while the  second inclusion in \eqref{eq:XalphaSubsetTildeXalpha} is a consequence of  the first  inclusion in
\cite[Proposition 2.8]{MR0206716} together with \cite[Theorem 3.1]{MR0206716} and the classical embedding $S(\infty,-\theta,X_0;\infty,1-\eta,X_1) \subset S(1,-\theta',X_0;\infty,1-\theta',X_1)$ whenever $X_1 \subset X_0$ and $0 < \theta' < \theta < 1$.

 %--------------------------------------------------------------------------------------------------------------------------------------------------------------

\subsection{An abstract spectral mapping theorem}

We present in this section a {\it ``principal spectral mapping theorem"} for a class of semigroup generators $\Lambda$ on a Banach space which split as a
hypodissipative part $\BB$ and a ``more regular part" $\AA$,  as presented in the introduction. In order to do so, we introduce a more accurate version of 
the growth and regularizing  assumption  {\bf (H2-3)}, namely

  \begin{itemize}
   \item[{\bf (H2)}] there exist $\zeta \in (0,1]$ and $\zeta' \in [0,\zeta)$ such that $\AA$ is $\BB^{\zeta'}$-bounded and 
    there exists an integer $n \ge 1$ such that for any $a > a^*$, there holds 
%\beqn\label{eq:hypH2bis}
%\,\, \forall \, t \ge 0, \qquad 
%\|  (\AA S_\BB)^{(*n)}(t) \|_{\BBB(X,D(\Lambda^\zeta))}  \le  C_{a,n,\zeta} \, e^{at}
%\eeqn
%or 
\beqn\label{eq:hypH2}
\,\, \forall \, t \ge 0, \qquad 
\|  S_\BB * (\AA S_\BB)^{(*n)}(t) \|_{\BBB(X,D(\Lambda^\zeta))}  \le  C_{a,n,\zeta} \, e^{at}
\eeqn
for a constant $C_{a,n,\zeta}  \in (0,\infty)$.    
\end{itemize}

\begin{theo} \label{theo:SpectralMapping} 
Consider  a Banach space $X$,  the  generator $\Lambda$ of a  semigroup $S_\Lambda(t) = e^{t\Lambda} $ on $X$, two real numbers $a^*, a' \in \R$, $a^* < a'$, 
and assume that the spectrum $\Sigma(\Lambda)$ of $\Lambda$ satisfies the following separation condition
\beqn\label{eq:SepSpectrLambda}
\Sigma (\Lambda) \cap \Delta_{a^*} \subset   \Delta_{a'}.
\eeqn
The  following   growth estimate on the semigroup

 \begin{itemize}
 
 \item[{\bf (1)}] there exists a projector  $\Pi  \in \BBB(X)$  satisfying 
    $ \Lambda \Pi = \Pi \Lambda$, $\Lambda_1 :=  \Lambda_{|X_1} \in \BBB(X_1)$, $X_1 := R \Pi$,   $\Sigma(\Lambda_1) \subset \Delta_{a^*}$ and for any real number
    $a  > a^*$ there exists a constant $C_{a}$ such that 
  \begin{equation}\label{bddSlambda1st} 
     \forall \, t \ge 0, \quad 
     \Bigl\| e^{t \, \Lambda} (I-\Pi) 
     \Bigr\|_{\mathscr{B}(X)} 
     \le C_{a} \, e^{a \, t}  , 
   \end{equation}

  \end{itemize}
  \smallskip\noindent
is equivalent to the following splitting structure  hypothesis  

\smallskip
 \begin{itemize}
 \item[{\bf (2)}] there exist two operators $\ \AA, \BB \in  \mathscr{C}(X)$,  such that 
  $\Lambda =  \AA + \BB$ and  hypothesizes {\bf (H1)} 
  and {\bf (H2)} are satisfied. 
 \end{itemize}

Moreover, under assumption {\bf (2)},  for any $a > a^*$ there exists an explicitly computable constant $M = M(a,\AA,\BB)$   such that 
\beqn\label{eq:localizSpectrL} 
\Sigma (\Lambda) \cap \Delta_{a}   \subset B(0,M) := \{ z \in \C; \, |z| < M \}.  
\eeqn

\end{theo}

\begin{rem} \label{rem:theoEquivSGX}
     
 \begin{itemize}
   \item[{\bf (a)}] Theorem~\ref{theo:SpectralMapping}  gives a characterization (and thus a criterium with the conditions {\bf (H1)} and {\bf (H2)}) for an operator $\Lambda$ to satisfy a partial (but principal)  spectral mapping theorem under the only additional assumption that the spectrum satisfies a separation condition. 
   
  \item[{\bf (b)}] Hypothesis  {\bf (H1)} holds for any $\ell \in \N$ if it is true for $\ell=0$ (that is 
 $(\BB-a)$ is  hypodissipative in $X$ for any real number $a > a^*$) and $\AA \in \BBB(X)$.  
  
  \item[{\bf (c)}]  
 The implication {\bf (1)} $\Rightarrow$ {\bf (2)} is just straightforward by taking $\AA :=  \Lambda \Pi$ and $\BB:= \Lambda (I-\Pi)$. With such a choice
  we have $\AA  \in \BBB(X)$ by assumption,  next $\AA \in \BBB(X,D(\Lambda))$, because $\Lambda \AA = \AA^2  \in \BBB(X)$, and $\BB-a$ is hypodissipative in $X$ for any real number $a > a^*$ so that  hypothesis  {\bf (H1)}  is satisfied as well as  hypothesis {\bf (H2)} with $n=1$, $\zeta = 1$ and $\zeta' = 0$.

  \item[{\bf (d)}]   We believe that the implication {\bf (2)} $\Rightarrow$ {\bf (1)} is  new. It can be seen as a condition under which a
   {\it ``spectral mapping theorem for the  principal part of the spectrum holds"} in the sense that  \eqref{eq:PartialSpectralMapping} 
   holds. Indeed, defining $\Lambda_0 := \Lambda (I-\Pi)$, for any $a > a^*$ there holds $\Sigma(\Lambda)  = \Sigma(\Lambda_0)  \cup \Sigma(\Lambda_1)$,
   $\Sigma(e^{t\Lambda})  = \Sigma(e^{t\Lambda_0})  \cup \Sigma(e^{t\Lambda_1})$,   $\Sigma(e^{t\Lambda_1}) = e^{t\Sigma(\Lambda_1)}$ (because $\Lambda_1 \in \BBB(X_1)$), $\Sigma( \Lambda_0) \subset \Delta_{a}$ (by hypothesis) and  $\Sigma(e^{t\Lambda_0}) \subset \Delta_{e^{at}}$ (from the conclusion   \eqref{bddSlambda1st}). 
In particular, under assumption {\bf (2)} the spectral bound $s(\Lambda)$ and the growth bound $\omega(\Lambda)$ coincide if they are at the right hand side of $a^*$, or in other words \eqref{eq:s=omega} holds.  

    \item[{\bf (e)}]  When $a^* < 0$ the above result  gives a characterization of uniformly exponentially stable semigroup  (see e.g. \cite[Definition V.1.1]{EngelNagel})  in a Banach space framework. 
    More precisely, if assumption {\bf(2)} holds with $a^* < 0$, then 
    $$
    S_\Lambda \hbox{ is uniformly exponentially stable iff } s(\Lambda) < 0.
    $$
   That last assertion has to be compared to the Gearhart-Pr\"uss Theorem which gives another 
   %(and in some sense weaker) 
   characterization of uniform exponential stability  in a Hilbert framework, see  \cite{MR0461206,MR743749,AGLMNS} as well as \cite[Theorem V.1.11]{EngelNagel} for a comprehensive proof.    
  
\item[{\bf (f)}]  Although the splitting condition in {\bf (2)} may seem to be strange, it is in fact quite natural for many partial differential operators, including numerous cases of operators which have not any self-adjointness property, as that can be seen in the many examples studied in  \cite{Mcmp,MMcmp,MR2476685,GMM,AGGMMS,MR2832638,MR2821681,KleberLandau,EM-KS1,TristaniInelastic,CEM-KS2}. 

 \item[{\bf (g)}] For a sectorial operator $\BB$,  hypothesis  {\bf (H1)}  holds for any operator $\AA$ which is suitably  bounded with respect to $\BB$.
 More precisely, in a Hilbert space framework, for a ``hypo-elliptic operator $\BB$ of order $\zeta$" and for a  $\BB^{\zeta'}$-bounded operator $\AA$, with $\zeta' \in [0,\zeta)$, in the sense that 
for any $f \in D(\BB)$ 
$$
\qquad\quad (\BB f,f) \le - a \, \| \BB^{\zeta} f \|^2 + C \| f \|^2 
\quad\hbox{and}\quad
\| \AA f \| \le C \, ( \| \BB^{\zeta'} f \| + \| f \|), 
$$
 then {\bf (H1)} holds and {\bf (H2)} holds with $n=0$. A typical example is  $\BB = - \Delta$ and $\AA = a(x) \cdot  \nabla$ in the space $L^2(\R^d)$ with $a \in L^\infty(\R^d)$.
 In that case,   Theorem~\ref{theo:SpectralMapping} is nothing but the classical spectral mapping theorem
 which is known to hold in such a sectorial framework. We refer to  \cite[section 2.5]{Pazy} and \cite[Section II.4.a]{EngelNagel} for an introduction to sectorial operators (and analytic semigroups) as well as \cite[Section IV.3.10]{EngelNagel} for a proof of the spectral mapping theorem in that framework. 
 
 \item[{\bf (h)}] Condition {\bf (H2)} is fulfilled for the same value of $n$  if $\AA$ and $\BB$ satisfy {\bf (H1)} as well as 
 $$
t \mapsto \|  (\AA S_\BB)^{(*n)}(t ) \|_{\BBB(X,D(\Lambda^\zeta))} \, e^{-a t} \in L^p (0,\infty) 
$$
 for any $a > a^*$ and for some $p \in [1,\infty]$. 

\item[{\bf (i)}]    It is    worth emphasizing that our result does not require any kind of ``regularity" on the semigroup as it is usually the case for ``full" spectral mapping theorem. In particular,  we do not require that the semigroup is eventually norm continuous as in \cite{MR0089373} or \cite[Theorem IV.3.10]{EngelNagel}.  Let us also stress that some ``partial"  spectral mapping theorems have been obtained in several earlier papers as in  \cite{MR0103419,MR0259662,Voigt80,MR1107184,MR1767391} under stronger assumptions than {\bf (H1)-(H2)}. See also Remark~\ref{rem:theoWeyl}-{\bf (e)}.

 \end{itemize}

 \end{rem}

 \medskip
\noindent{\sl Proof of Theorem~\ref{theo:SpectralMapping}.} We only prove that {\bf (2)} implies {\bf (1)} since the reverse implication is 
clear (see Remark~\ref{rem:theoEquivSGX}-{\bf (c)}). 

\smallskip
The proof is split into four steps. From now on, let us fix $a \in (a^*,a')$.

 \smallskip
\noindent {\sl Step 1. }   We establish the cornerstone estimate 
\beqn\label{eq:Vestim1}
\forall \, z \in \Delta_a, \quad
\|  (\AA \RR_\BB(z))^{n+1}  \|_{\BBB(X)} \le C_{a} \, {1 \over \langle z \rangle^\alpha}, 
\eeqn
where  $n \ge 1$ is the integer given by assumption {\bf (H2)} and $\alpha := {(\zeta-\zeta')^2 / 8} \in (0,1)$.
We recall the notation $X_s = D(\Lambda^s)$ for $0 \le s \le 1$, with $X_0=X$ and $X_1 = D(\Lambda)$.  
On the one hand, from \eqref{eq:hypH2} we have 
$$
\forall \, z \in \Delta_a, \quad \| \RR_\BB (z) \, (\AA \RR_\BB (z))^n \|_{X \to X_\zeta} \le C'_{a,n}  , 
$$
for a constant $C_{a,n}' \in (0,\infty)$ only depending on $a$,  $a^*$ and $C_{a,n,\zeta}$. As a consequence, writing 
$$
\RR_\BB (z)^{1-\eps} \, (\AA \RR_\BB (z))^n  = (\BB -z)^\eps \, [ \RR_\BB (z) \, (\AA \RR_\BB (z))^n], 
$$
we get 
\beqn\label{eq:SGborne1}
\forall \, z \in \Delta_a, \quad \| \RR_\BB  (z)^{1-\eps}  \, (\AA \RR_\BB (z))^n \|_{X \to X_{\zeta-\eps} } \le C_{a,n} , 
\eeqn
for any $\eps \in [0,\zeta]$. 
 
 On the other hand, we claim that 
\beqn\label{eq:SGborne2}
\forall \, z \in \Delta_a, \quad
\| \RR_\BB  (z)^\eps \|_{\tilde X_{\zeta-\eps} \to \tilde X_{\zeta'}} \le C_{a,\eps} / \langle z \rangle^{\eps \, (\zeta-\zeta' - \eps)},
\eeqn
for any   $0 \le \eps < \zeta-\zeta'$  and a constant $C_{a,\eps}$. 
First, from the resolvent identity
$$
\RR_\BB(z) = z^{-1} \, (\RR_\BB(z) \BB - I),
$$
we have for $i,j = 0,1$, $i \ge j$,
$$
\forall \, z \in \Delta_a, \quad \| \RR_\BB(z) \|_{X_i \to X_j} \le C^{1,a}_{i,j}(z), 
$$
with $C^{1,a}_{0,0}(z) = C^{1,a}_{1,1}(z) = C^{1,a}$ and $C^{1,a}_{1,0}(z) = C^{1,a}_{1,0}/\langle z \rangle$, $C^{1,a}_{i,j} \in \R_+$. 
Next, thanks to the first representation formula in \eqref{eq:FracPower}, we
clearly have 
$$
\| \RR_\BB(z)^\eps \|_{X_1 \to X_0} \le C(a,\eps) \int_0^\infty {\lambda^{-\eps} \over \langle z \rangle + \lambda} \, d\lambda 
\le C^{\eps,a}_{1,0}(z) := {C^{\eps,a}_{1,0}  \over \langle z\rangle^\eps}.
$$
Last, thanks to the interpolation theorem \cite[Th\'eor\`eme 3.1]{MR0165343} we get with some $0 \le \theta < \theta' < 1$
$$
\| \RR_\BB(z)^\eps \|_{\tilde X_{\theta'} \to \tilde X_\theta} \le \| \RR_\BB(z)^\eps \|_{X_1 \to X_1}^{\theta} \, \| \RR_\BB(z)^\eps \|_{X_1 \to X_0}^{\theta'-\theta} \, \| \RR_\BB(z)^\eps \|_{X_0 \to X_0}^{1-\theta'} , 
$$
from which we conclude to \eqref{eq:SGborne2}.  
Writing now
$$
\RR_\BB  (z) \, (\AA \RR_\BB (z))^n = \RR_\BB (z) ^\eps \, [ \RR_\BB (z)^{1-\eps} \, (\AA \RR_\BB (z))^n],
$$
with the optimal choice $\eps = (\zeta-\zeta')/2$, 
we finally deduce from \eqref{eq:SGborne1}, \eqref{eq:SGborne2} and the inclusions \eqref{eq:XalphaSubsetTildeXalpha} that 
$$
\forall \, z \in \Delta_a, \quad
\| \RR_\BB  (z) \, (\AA \RR_\BB (z))^n \|_{X \to X_{\zeta'}}   \le {C_a  \over \langle z \rangle^\alpha}. 
$$
Because  $\AA$ is $\BB^{\zeta'}$-bounded we conclude to \eqref{eq:Vestim1}. 

\smallskip
\noindent {\sl Step 2. }  We prove \eqref{eq:localizSpectrL}. From the splitting $\Lambda = \AA + \BB$,  we have 
$$
\RR_\Lambda(z) = \RR_\BB(z) - \RR_\Lambda(z) \, \AA \,  \RR_\BB(z)
$$
on the open region of $\C$ where $\RR_\Lambda$ and $\RR_\BB$ are well defined functions (and thus analytic), 
and iterating the above formula, we get  
\beqn\label{eq:RL=ABN}
 \RR_\Lambda(z) =   \sum_{\ell=0}^{N-1}  (-1)^{\ell} \RR_\BB(z) \, (\AA \RR_\BB(z))^{\ell} + \RR_\Lambda(z)  (-1)^{N} (\AA \RR_\BB(z))^{N}
\eeqn
for any integer $N \ge 1$.  
We define 
\beqn\label{def:UU}
\UU(z) := \RR_\BB(z) - \dots + (-1)^{n} \RR_\BB(z) \, (\AA \RR_\BB(z))^{n}
\eeqn
and 
\beqn\label{def:VV}
\VV(z) := (-1)^{n+1} (\AA \RR_\BB(z))^{n+1}, 
\eeqn
where  $n \ge 1$ is again the integer given by assumption {\bf (H2)}. With $N=n+1$, we may rewrite the identity \eqref{eq:RL=ABN} as 
\beqn\label{eq:RLV=U}
  \RR_\Lambda(z) ( I - \VV(z)) = \UU(z). 
\eeqn
Observing that \eqref{eq:Vestim1}  rewrites as
\beqn\label{eq:Vestim1bis}
\forall \, z \in \Delta_a, \quad
\| \VV(z) \|_{\BBB(X)} \le C_{a} \, {1 \over \langle z \rangle^\alpha},
\eeqn
we have in particular, for $M$ large enough, 
\beqn\label{eq:Vestim2}
z \in \Delta_a, \, |z| \ge M \quad\Rightarrow \quad \| \VV(z) \|_{\BBB(X)} \le {1 \over 2}.
\eeqn
As a consequence, in the region $\Delta_a \backslash B(0,M)$ the operator $I-\VV$ is invertible and thus $  \RR_\Lambda$ is well-defined because of 
\eqref{eq:RLV=U},
which is nothing but \eqref{eq:localizSpectrL}.

\smallskip
\noindent {\sl Step 3. }  We write a representation formula for the semigroup which follows from some  classical complex analysis arguments  as in \cite[Proof of Theorem 2.13]{GMM} or \cite[Chapter 1]{Pazy}). 
First, the inclusion \eqref{eq:localizSpectrL}  together with the separation condition \eqref{eq:SepSpectrLambda} implies that 
$$
\Gamma := \Sigma(\Lambda) \cap \Delta_a \subset \Delta_{a'} \cap B(0,M),
$$
so that $\Gamma$ is a compact set of $\Delta_a$ and we may define the projection operator $\Pi$ on the 
corresponding eigenspace thanks to the Dunford integral  \eqref{def:SpectralProjectionGamma}. The projector $\Pi$ fulfils all the requirements stated in {\bf (1)}. 
In order to conclude we just have to prove \eqref{bddSlambda1st}.
In order to do so, for any integer $N \ge 0$, we write 
\bean
S_\Lambda \,  (I-\Pi) 
&=& (I-\Pi) \, S_\Lambda
\\
&=&  (I-\Pi)  \sum_{\ell=0}^{N-1}  (-1)^\ell S_\BB *  (\AA S_\BB)^{(*\ell)} +  ((I-\Pi) S_\Lambda) * (\AA S_\BB)^{(*N)}, 
\eean
where we have used the iteration of the  Duhamel formula \eqref{eq:stopDyson} in the second line.
For $b > \max(\omega(\Lambda),a)$,  we may use the inverse Laplace formula 
\bean
\TT(t) f 
&:=& ((I-\Pi) S_\Lambda) * (\AA S_\BB)^{(*N)} f  
\\
&=& \lim_{M' \to \infty}  {i \over
  2\pi}\int_{b-iM'}^{b+iM'} e^{zt} \,  (-1)^{N+1} \,  (I-\Pi)   \RR_\Lambda(z) \, (\AA \RR_\BB(z))^{N}  f \, dz
\eean
for any $f \in D(\Lambda)$  and $t \ge 0$, and we emphasize that the term $\TT(t) f$ might be only defined as a semi-convergent  integral. 
Because $z \mapsto  (I-\Pi)   \RR_\Lambda(z) \, (\AA \RR_\BB(z))^{N+1}$ is an analytic function on a neighborhood of $\bar \Delta_a$, 
we may move the segment on which the integral is performed, and we obtain the representation formula
\bear\label{eq:SL=2}
\quad &&S_\Lambda(t) (I-\Pi) f
=  \sum_{\ell=0}^{N-1}  (-1)^\ell (I-\Pi) S_\BB *  (\AA S_\BB)^{(*\ell)} (t)f
\\ \nonumber
&&\qquad +  \lim_{M' \to \infty}  {i \over
  2\pi}\int_{a-iM'}^{a+iM'} e^{zt} \,  (-1)^{N+1} \,  (I-\Pi)   \RR_\Lambda(z) \, (\AA \RR_\BB(z))^{N}  f \, dz, 
\eear
 for any $f \in D(\Lambda)$  and $t \ge 0$.

 \smallskip
\noindent {\sl Step 4. } In order to conclude we only have to explain why the last term in \eqref{eq:SL=2} is appropriately bounded for $N$ large enough. 
Thanks to \eqref{eq:RLV=U}, we have
\beqn\label{eq:RL=VU}
\WW(z) := \RR_\Lambda(z) \,  (\AA \RR_\BB(z))^{N}  = \UU(z) \, \sum_{\ell=0}^\infty \VV(z)^\ell \,  (\AA \RR_\BB(z))^{N} ,
\eeqn
provided the RHS series converges.   On the one hand, defining $N := ( [1/\alpha]+1) (n+1)$ and $\beta :=  ( [1/\alpha]+1) \alpha > 1$,  we get from \eqref{eq:Vestim1} that 
\beqn\label{eq:ARBN+1}
\|  (\AA \RR_\BB(z))^{N} \|_{\BBB(X)} \le \|  \VV(z) \|_{\BBB(X)}^{[1/\alpha]+1}
\le {C \over \langle y \rangle^\beta}
\eeqn
for any $z= a+ iy$, $|y| \ge M$. On the other hand,  we get from \eqref{eq:Vestim2} that the series term in \eqref{eq:RL=VU} is normally convergent uniformly in $z= a+ iy$, $|y| \ge M$. All together, we obtain
\bear\nonumber
\| \WW(z)\|_{\BBB(X)} 
&\le& \|  \UU(z) \|_{\BBB(X)} \, \Bigl\| \sum_{\ell=0}^\infty \VV(z)^\ell \Bigr\|_{\BBB(X)} \,  \| (\AA \RR_\BB(z))^{N} \|_{\BBB(X)}
\\ \label{eq:bdWW}
&\le& {C \over | y |^\beta}
\eear
for any $z= a+ iy$, $|y| \ge M$. 

\smallskip
In order to estimate the last term in \eqref{eq:SL=2}, we write 
\bean
&&\Bigl\| \lim_{M' \to \infty}  {i \over  2\pi}\int_{a-iM'}^{a+iM'}  e^{zt} \, (I-\Pi) \WW(z) \, dz \Bigr\|_{\BBB(X)} 
\\
&&\qquad\le {e^{at} \over 2\pi} \int_{a-iM}^{a+iM} \|  (I-\Pi) \RR_\Lambda(z) \|_{\BBB(X)} \,  \| (\AA \RR_\BB(z))^{N} \|_{\BBB(X)} \, dy
\\
&&\qquad\qquad+ {e^{at} \over 2\pi} \int_{\R \backslash [-M,M]} \| (I-\Pi) \WW(a+iy) \|_{\BBB(X)} \, dy,
\eean
where the first integral is finite thanks to  $ \Sigma(\Lambda) \cap [a-iM,a+iM] = \emptyset$ and \eqref{eq:ARBN+1}, and the second
integral is finite because of \eqref{eq:bdWW}.
\qed

\medskip
We give now two variants of Theorem~\ref{theo:SpectralMapping} which are sometimes easier to apply. First, for a generator $\Lambda$ of a semigroup $S_\Lambda$ with
the splitting structure \eqref{eq:L=A+B} we introduce an alternative growth assumption to {\bf (H1)}, namely  
 \begin{itemize}\item[{\bf (H1$'$)}] for some $a^* \in \R$ and for any $a > a^*$,  the operator $\BB-a$ is hypodissipative and  there exists a constant $C_{a} \in (0,\infty)$ such that 
\beqn\label{eq:H1bis}
  \int_0^\infty e^{-a \, s} \, \| \AA \, S_\BB(s) \|_{\BBB(X)} \, ds \le C_{a}. 
\eeqn
  \end{itemize}

\begin{rem}\label{rem:hypH1}  

\begin{itemize}  \item[] Estimate \eqref{eq:H1bis} is reminiscent of the usual condition under which for  a given generator $\BB$ of a semigroup $S_\BB$ the perturbed operator $\AA+\BB$ also generates a semigroup (see \cite{MR0044737}, \cite{MR0209900}, \cite[condition (1.1)]{MR0445335}).  Here however the condition $C_{a} < 1$ is not required  since we have already made the assumption that $S_{\Lambda}$ exists.  \end{itemize} \end{rem}

It is clear by  writing 
$$
S_\BB * (\AA S_\BB)^{(*\ell)} = [ S_\BB * (\AA S_\BB)^{(*\ell-1)}] * (\AA S_\BB)   
$$
and performing an iterative argument, that {\bf (H1$'$)} implies {\bf (H1)}. We then immediately deduce from Theorem~\ref{theo:SpectralMapping} 
a first variant: 

\begin{cor}\label{cor:ThSM1}  Consider a Banach space $X$, the  generator  $\Lambda$ of a  semigroup $S_\Lambda(t) = e^{t\Lambda} $ on $X$ and a real number $a^* \in \R$. Assume that the spectrum $\Sigma(\Lambda)$ of $\Lambda$ satisfies the  separation condition \eqref{eq:SepSpectrLambda} and that there exist two operators $\ \AA, \BB \in  \mathscr{C}(X)$ such that $\Lambda =  \AA + \BB$ and  hypothesises {\bf (H1$'$)}   and {\bf (H2)} are met. Then the conclusions   \eqref{eq:localizSpectrL} ans {\bf (1)}
 in Theorem~\ref{theo:SpectralMapping}  hold. 
\end{cor}

Next, for a generator $\Lambda$ of a semigroup $S_\Lambda$ with
the splitting structure \eqref{eq:L=A+B} we introduce the alternative growth and regularizing assumptions to {\bf (H1$'$)} and {\bf (H2)}, namely  
 \begin{itemize}
 \item[{\bf (H1$''$)}] 
 for some $a^* \in \R$, $\zeta \in (0,1]$ and for any $a > a^*$, $\zeta' \in [0,\zeta]$ 
 the operator $\BB-a$ is hypodissipative and assumption {\bf (H1) } or {\bf (H1$'$)} 
 also hold with $\BB(X)$ replaced by $\BB(D(\Lambda^{\zeta'}))$; 
 \item[{\bf (H2$''$)}] 
  there exist an integer $n \ge 1$ 	and a real number $ b > a^*$ such that 
$$
 \forall \, t \ge 0, \qquad 
\|  (\AA S_\BB)^{(*n)}(t) \|_{\BBB(X,D(\Lambda^\zeta))}  \le  C_{b,n,\zeta} \, e^{bt},
$$
for a constant $C_{b,n,\zeta}  \in (0,\infty)$. 
 \end{itemize}

Our second variant of Theorem~\ref{theo:SpectralMapping} is the following. 

\begin{cor}\label{cor:ThSM2}  Consider a Banach space $X$, the  generator  $\Lambda$ of a  semigroup $S_\Lambda(t) = e^{t\Lambda} $ on $X$ and a real number $a^* \in \R$. Assume that the spectrum $\Sigma(\Lambda)$ of $\Lambda$ satisfies the  separation condition \eqref{eq:SepSpectrLambda} and that there exist two operators $\ \AA, \BB \in  \mathscr{C}(X)$ such that $\Lambda =  \AA + \BB$ and  hypothesises {\bf (H1$''$)}   and {\bf (H2$''$)} are met. Then the conclusions   \eqref{eq:localizSpectrL} ans {\bf (1)} in Theorem~\ref{theo:SpectralMapping}  hold. 
\end{cor}

Corollary~\ref{cor:ThSM2} is an immediate consequence of Theorem~\ref{theo:SpectralMapping} or Corollary~\ref{cor:ThSM1} together with the following
simple variant of \cite[Lemma 2.15]{GMM} which makes possible to deduce the more accurate regularization and growth condition {\bf (H2)} 
from the rough  regularization and growth condition {\bf (H2$''$)}. 

 \begin{lem}  \label{lem:Tn} 
Consider two Banach spaces $E$ and $\EE$ such that
$E \subset \EE$ with dense and continuous embedding. 
Consider some operators $\LL$,  $\AA$ and $\BB$ on $\EE$ such that $\LL$ splits as $\LL = \AA + \BB$, some real number $a^* \in \R$
and some integer $m \in \N^*$.
Denoting with the same letter $\AA$, $\BB$ and $\LL$ the restriction of these
operators on $E$, we assume that the two following dissipativity conditions  are satisfied: 

\begin{itemize}

\item[{\bf (i)}]    for any $ a > a^*$,  $\ell \ge m$,  $X = E$ and  $X = \EE$, there holds
\[
\| (\AA \, S_\BB)^{(*\ell)} \|_{\BBB(X) }  \le    C_{a,\ell} \,   e^{a t}  ;
\]
  \item[{\bf (ii)}]  for some constants $b \in \R$, $b > a^*$,   and $C_b \ge 0$, there holds 
\[
\| (\AA \, S_\BB)^{(*m)} \|_{\BBB(\EE,E) }  \le    C_b \,   e^{b t}   .
\]
   \end{itemize}

\smallskip\noindent

Then for any $a > a^*$, there exist some constructive constants $n  = n(a)\in \N$, $C_{a} \ge 1$ such that
$$
\forall \, t \ge 0 \qquad 
\| (\AA S_\BB)^{(*n)}(t) \|_{\BBB(\EE,E)} \le C_{a}  \, e^{ a \, t}.
$$ 
\end{lem}

\medskip\noindent
{\sl Proof of Lemma~\ref{lem:Tn}. } We fix $a > a^*$ and $a' \in (a^*,a)$,  and we note $T := (\AA \, S_\BB)^{(*m)}$. For $n = p \, m$, $p \in \N^*$, we write 
\begin{eqnarray*}
 (\AA S_\BB)^{(*n)} (t) &=& T^{(*p)}(t) 
  \\
  &=& \int_0^t \dd t_{p-1} \, \int_0^{t_{p-1}} \dd t_{p-2} \, \dots \,
  \int_0^{t_2} \dd t_1 \, T (\delta_p) \, \dots \, T (\delta_1)
\end{eqnarray*} 
with 
\[
\delta_1 = t_1, \quad \delta_2 = t_2 - t_1, \ \dots, \  
\delta_{p-1} = t_{p-1} - t_{p-2} \mbox{ and } \ \delta_p = t-t_{p-1}.
\]
For any $p \ge 1$, there exist at least one increment $\delta_{r}$, $r \in  \{1, ..., p \}$,   such that $\delta_{r} \le t/p$. 
Using {\bf (ii)} in order to estimate $\|T (\delta_{r})
\|_{\BBB(\EE,E)}$ and   {\bf (i)}  in order to bound the
other terms $\| T (\delta_q) \|_{\BBB(X)}$ in the appropriate
space $X = \EE$ or $X=E$, we have  
\bean
&&\|T (\AA S_\BB)^{(*n)} (t)   \|_{\BBB(\EE,E)} \le \\
&&\quad  \le    
   \int_0^t \!\! dt_{p-1} \int_0^{t_{p-1}} \!\!\!\! 
   \dd t_{p-2} \, \dots \int_0^{t_2} \!\! \dd t_1 C'_{b}  \, e^{b \, \delta_r}   \, \prod_{q \not= r } C_{a',m}  \, e^{a' \, \delta_q} 
      \\
&&\quad  \le    
C'_{b}  \, C_{a',m}^{p-1} \,  e^{a' \, t}  \,  e^{ b \,  t /p }
   \int_0^t \!\! \dd t_{p-1} \int_0^{t_{p-1}} \!\!\!\! \dd t_{p-2} \,
   \dots \int_0^{t_2} \!\! \dd t_1   
     \\
&&\quad =  
 C \, t^{p-1}  \, e^{ (a' +  b/p ) \, t}  \le C' \, e^{at},
\eean
by taking $p$ (and then $n$) large enough so that $a'+b/p < a$. 
\qed

 %%%%%%%%%%%%%%%%%%%% Spectral mapping and Weyl's  %%%%%%%%%%%%%%%%%%%%%%%%%%%%%%%%%%%

\bigskip
\section{Weyl's Theorem for semigroup generators}
\label{sec:Weyl}
\setcounter{equation}{0}
\setcounter{theo}{0}

\subsection{An abstract semigroup Weyl's Theorem}
We present in this section a version of Weyl's Theorem about compact perturbation of dissipative generator in the spirit of 
the above spectral mapping theorem. For that purpose, we introduce the growth and compactness assumption
% for an operator $\Lambda$

 \begin{itemize}
 
 \item[{\bf (H3)}]  for the same integer $n \ge 1$ as in assumption  {\bf (H2)} and with the same notation, the time indexed family of operators $(\AA S_\BB)^{(*(n+1))}(t)$ satisfies the growth and compactness estimate  
\[
\forall \, a > a^*, \quad \int_0^\infty 
\|  (\AA S_\BB)^{(*(n+1))}(t) \|_{\BBB(X,Y)}  \, e^{-at} \, dt \le  C''_{n+1,a},
\]
for some constant $C''_{n+1,a} \ge 0$ and some (separable) Banach space  $Y$ such that $Y  \subset X$ with compact embedding.

   \end{itemize}

\begin{theo}  \label{theo:GalWeyl}
Consider a Banach space $X$, the  generator $\Lambda$ of a  semigroup $S_\Lambda(t) = e^{t\Lambda} $ on $X$ and a real number $a^* \in \R$. 
The  following quantitative growth estimate on the semigroup

 \begin{itemize}
 
 \item[{\bf (1)}] for any $a > a^*$ there exist  an integer $J \in \N$,  a finite family of  distinct complex  numbers $\xi_1, ..., \xi_J \in \bar \Delta_{a}$,
    % (with the convention $\{\xi_1, \dots, \xi_k  \}= \emptyset$ if $k=0$), 
    some finite rank projectors $\Pi_1, \dots \, \Pi_J \in \BBB(X)$ and 
    some operators $T_j \in \BBB(R \Pi_j)$, satisfying 
    $\Lambda \Pi_j = \Pi_j \Lambda = T_j \Pi_j$, $\Sigma(T_j) = \{ \xi_j \}$, in particular
    \[
\quad \Sigma(\Lambda) \cap \bar \Delta_a = \left\{ \xi_1, \dots,
    \xi_J \right\} \subset \Sigma_d(\Lambda),
\]
 and  a constant $C_{a}$ such that 
  \begin{equation}\label{bddSlambda} 
     \forall \, t \ge 0, \quad 
     \Bigl\| e^{t \, \Lambda} - \sum_{j=1}^J e^{t \, T_j} \, \Pi_j 
     \Bigr\|_{\mathscr{B}(X)} 
     \le C_{a} \, e^{a \, t}  , 
   \end{equation}

  \end{itemize}

\smallskip\noindent
is equivalent to the following splitting structure of the generator

\smallskip
 \begin{itemize}
 \item[{\bf (2)}] there exist two operators $\ \AA, \BB \in  \mathscr{C}(X)$ such that 
  $\Lambda =  \AA + \BB$ and  hypothesises {\bf (H1)}, {\bf (H2)} and {\bf (H3)}  are met.

\end{itemize}
 
\end{theo}

%
%
%  \begin{theo}
%  \label{theo:GalWeyl}
%In a Banach space $X$, consider two operators $\AA, \BB \in \CCC(X)$ such that 
%
%\begin{itemize} 
%
%\item[{\bf (i)}]  $\AA$ and $\BB$ satisfy conditions  {\bf (H1)} and {\bf (H2)} in  Theorem~\ref{theo:SpectralMapping};
%
%\item[{\bf (ii)}] $ \AA$ is $\BB$-compact in the sense that for the same integer $n \ge 1$ as in assumption  {\bf (H2)},  the time indexed family of operators $(\AA S_\BB)^{(*(n+1))}(t)$ satisfies the growth and compactness estimate  
%\[
%\forall \, a > a^*, \,\, \forall \, t \ge 0, \qquad 
%\|  (\AA S_\BB)^{(*(n+1))}(t) \|_{\BBB(X,Y)}  \le  C''_{n+1,a} \, e^{at},
%\]
%for some constant $C''_{n+1,a} \ge 0$ and some (separable) Banach space  $Y$ such that $Y  \subset X$ with compact embedding.
%\end{itemize}
%
%Then, for any $a > a^*$ there exist a computable constant $M$ and a finite family of  distinct complex  numbers $\xi_1, ..., \xi_J \in \Delta_{a}$ such that 
%$\Lambda = \AA + \BB$ satisfies
%\bean
%&& \Sigma (\Lambda) \cap \Delta_{a} = \{\xi_1, ... , \xi_J \} \subset \Sigma_d(\Lambda) \cap B(0,M),
%\eean and such that $S_\Lambda$ satisfies the growth estimate \eqref{bddSlambda}.
%%{\bf (1)} of Theorem~\ref{theo:SpectralMapping}. 
%  \end{theo}
%  
  
  \begin{rem} \label{rem:theoWeyl}
     
 \begin{itemize}
 
  \item[{\bf (a)}] When $\Sigma(\Lambda) \cap \Delta_{a^*} \not = \emptyset$ the above theorem gives a description of the principal asymptotic behaviour of the $C_0$-semigroup  $S_\Lambda$, namely it states that it is essentially compact  (see e.g.  \cite[Definition 2.1]{MR1297033}). As a matter of fact, if conditions {\bf (H1)}, {\bf (H2)} and {\bf (H3)}
 hold with  $a^* < 0$,  then the  $C_0$-semigroup  $S_\Lambda$ is either uniformly exponentially stable 
 or essentially compact.
  
  \item[{\bf (b)}]
  We adopt the convention 
  $\{ \xi_1, \dots, \xi_J \}= \emptyset$ if $J = 0$. 

  \item[{\bf (c)}] In the case when $\BB$ is sectorial and  $ \AA$ is $\BB$-compact, Theorem~\ref{theo:GalWeyl} is nothing but the classical  Weyl's Theorem on the spectrum 
   \cite{MR1511560} or 
   \cite[Theorem IV.5.35]{Kato} combined with the spectral mapping theorem, see \cite[Section IV.3.10]{EngelNagel} for instance. 
 
  \item[{\bf (d)}] Assumption {\bf (2) } in Theorem~\ref{theo:GalWeyl} is similar to the definition in  \cite{Voigt80} of the fact that {\it ``$\AA$ is $\BB$-power compact"}, but it is written  at the semigroup level rather than at the resolvent level. Under such a power compact hypothesis, Voigt establishes a generalisation of Weyl's Theorem. His proof uses the analytic property of the resolvent  function $\RR_\Lambda(z)$ obtained by Ridav and Vidav in \cite{MR0236741} as we present here. As a matter of fact,  Theorem~\ref{theo:GalWeyl}   is a simple consequence of \cite[Theorem 1.1]{Voigt80} together with Theorem~\ref{theo:SpectralMapping}. However, for the sake of completeness we give an elementary proof of Theorem~\ref{theo:GalWeyl} (which  consists essentially to prove again \cite[Theorem 1.1]{Voigt80} without the help of  \cite{MR0236741} and to apply  Theorem~\ref{theo:SpectralMapping}). 
 
 \item[{\bf (e)}]  It has been observed in \cite{Voigt80} (see also \cite{MR0103419,MR0259662,MR1107184,MR1767391})  that a similar conclusion as  \eqref{bddSlambda} holds under  a compactness condition on the remainder term
$$
\sum_{\ell = n-1}^\infty (\AA S_\BB)^{*\ell} = S_\Lambda * (\AA S_\BB)^{(*n)} = (S_\BB \AA )^{(*n)} * S_\Lambda
$$
in the Dyson-Phillips  series \eqref{eq:DPseries}. 
Hypothesis {\bf (H1)-(H2)-(H3)} provide some conditions on the only operators $\AA$ and $\BB$ such that the above remainder term enjoys such a nice property. 
   \end{itemize}
\end{rem}

 \noindent
 {\sl Proof of Theorem~\ref{theo:GalWeyl}. } Again, we only prove that {\bf (2)} implies {\bf (1)} since the reverse implication is 
clear (see Remark~\ref{rem:theoEquivSGX}-{\bf (c)}). 

 \smallskip
 The cornerstone of the proof is the use of the identity 
 $$
  \RR_\Lambda(z) ( I - \VV(z)) = \UU(z)
$$
established in \eqref{eq:RLV=U}, where we recall that the functions $\UU$ and $\VV$ defined by \eqref{def:UU} and  \eqref{def:VV} are analytic functions on $\Delta_{a^*}$ with values in $\BBB(X)$. Moreover, we have  
\beqn\label{eq:VVcompact}
\VV(z) \in \BBB(X,Y) \subset \KKK(X) \quad \forall \,  z \in \Delta_{a}
\eeqn
 because of assumption {\bf (H3)} and 
\beqn\label{eq:VVinf1/2}
\| \VV(z) \|_{\BBB(X)} \le 1/2  \quad \forall \,  z \in \Delta_a \cap B(0,M)^c
\eeqn
for some $M=M(a)$ thanks to  \eqref{eq:Vestim1}.

 \medskip\noindent
  {Step 1. } We prove that  $\Sigma(\Lambda) \cap \Delta_a$ is finite for any $a > a^*$. 
  
 Let us fix  $\xi \in \Delta_a$ and define $C(z) := I - \VV(z)$, $C_0 := C(\xi)$. Because of \eqref{eq:VVcompact} and thanks to the Fredholm alternative  \cite{MR1554993} (see also  \cite[Theor\`eme VI.6]{MR697382} for a modern and comprehensible statement and proof), there holds
 $$
 R(C_{0}) = N(C_{0}^*)^\perp, \quad \hbox{dim} N(C_{0}) = \hbox{dim} N(C_{0}^*) := N \in \N. 
 $$
If $N \ge 1$, we introduce $(f_{1},\dots,f_N)$ a basis of the null space $N(C_{0})$ and $(\varphi_{1},\dots,\varphi_{N})$ a family of independent linear forms in
 $X'$ such that $\varphi_{i}(f_{j})= \delta_{ij}$. 
 Similarly, we introduce $(\psi_{1},\dots,\psi_N)$  a basis of the null space $N(C_{0}^*)$ and $(g_{1},\dots,g_{N})$ a family of independent linear vectors in
 $X$ such that   $\psi_{i}(g_{j})= \delta_{ij}$. Then defining the projectors on $X$
 $$
 \pi_{0} := \sum_{i=1}^n f_{i} \, \varphi_{i}  \quad\hbox{and}\quad
  \pi_1 := \sum_{i=1}^n g_{i} \, \psi_{i},
 $$
 and $X_{1}:=  \pi_{0} (X)$, $X_{0}:= (I-\pi_{0} ) (X)$, $Y_{1}:= \pi_1(X)$, $Y_{0} := (I-\pi_1)(X)$,     we have 
  $$
 X = X_{0} \oplus X_1, \,\,\, X_1 = N(C_{0}) = \hbox{Vect}(f_{1},\dots,f_N), 
$$
and
$$
 X = Y_{0} \oplus Y_1, \,\,\,  Y_0 = R(C_{0}), \,\,\, Y_1  = \hbox{Vect}(g_{1},\dots,g_N).
 $$
 On the one hand, $C_{0} : X_{0} \to Y_{0}$ is  bijective by definition, and the family  of linear mappings  
\bean
D_0 (z)  := C_0^{-1}  \, (I-\pi_1) \,  C(z)_{|X_{0}}  : X_0 \to X_0,
\eean
is an analytic function with respect to the parameter $z \in \Delta_a$ and satisfies $D_{0}(\xi) = I_{X_{0}}$. We deduce that $D_{0}(z)$ is also invertible for $z$ belonging to a neighbourhood $B(\xi,r_0(\xi))$ of $\xi$, $r_0(\xi) > 0$, and then 
\beqn\label{eq:pi1perpCz}\quad
(I-\pi_1) \,  C(z)_{|X_{0}}  : X_0 \to Y_0 \,\,\, \hbox{is invertible for any} \,\, \,  z \in B(\xi,r_0(\xi)).
\eeqn
 On the other hand, the family  of linear mappings  
\bean
D_1 (z)  :=  \pi_1  \,  C(z)_{|Y_1}  : Y_1 \to Y_1,
\eean
is an analytic function with respect to the parameter $z \in \Delta_a$ and $D_{1}(z)$ is invertible for any $z \in \Delta_{a} \cap B(0,M')^c$ for $M'$ large enough,
because $\VV(z) \to 0$ as $\Im m z \to \infty$ from \eqref{eq:Vestim1}. Since $Y_1$ is 
finite dimensional we may define $z \mapsto \hbox{det}(D_1(z))$, which is an analytic function on  $\Delta_a$ and satisfies $\hbox{det}(D_1(z)) \not= 0$ for any $z \in \Delta_{a} \cap B(0,M')^c$. As a consequence, $z\mapsto \hbox{det}(D_1(z))$ has isolate zeros and since  $ \hbox{det}(D_1(z_{0})) = 0$ there exists a neighbourhood $B(\xi,r_1(\xi))$ of $\xi$, $r_1(\xi) > 0$,  such that  $\hbox{det}(D_1(z)) \not= 0$ for any $z \in  B(z_{0},r_1(\xi)) \backslash \{ \xi \}$ from which we deduce 
\beqn\label{eq:pi1Cz}\quad
\pi_1 \,  C(z)_{|Y_{1}}  : Y_1 \to Y_1 \,\,\, \hbox{is invertible for any} \,\, \,  z \in B(\xi,r_1(\xi))  \backslash \{ \xi \}.
\eeqn
Gathering \eqref{eq:pi1perpCz} and \eqref{eq:pi1Cz}, we get that $C(z) : X \to Y_0 \oplus Y_1 = X$ is surjective  for any $z \in B(\xi,r(\xi))  \backslash \{ \xi \}$,
$r(\xi) := \min (r_0(\xi), r_1(\xi) ) > 0$, and then bijective thanks to the Fredholm alternative. 

\smallskip
By compactness of the set $\bar\Delta_a \cap B(0,M')$, we may cover that set by a finite number of balls $B(\xi_j,r(\xi_j))$, $1 \le j \le J$, so that $C(z) = I - \VV(z)$ is invertible for any 
$z \in \Omega := \Delta_a \backslash \{ \xi_1, ... , \xi_J \}$. As a consequence, the function 
$$
\WW(z) := \UU(z) \, (I - \VV(z))^{-1}
$$
is well defined and analytic on $\Omega$ and coincide with $\RR_\Lambda(z)$ on the half plane $\Delta_b$, where $b = \omega(\Lambda)$ is the growth bound of $S_\Lambda$. These facts  immediately imply that $\RR_\Lambda = \WW$ on $\Omega$ and then
$$
\Sigma(\Lambda) \cap \Delta_a \subset \{ \xi_1, ... , \xi_J \} = \Delta_a \backslash \Omega.
$$
 
 \medskip\noindent
  {Step 2. } We prove that $\Sigma(\Lambda) \cap \Delta_{a^*} \subset \Sigma_d(\Lambda)$. 
  Indeed, for any $\xi \in \Sigma(\Lambda) \cap \Delta_{a^*}$, we may write with the notation of step 1
\bean
\Pi_{\Lambda,\xi} (I-\pi_1)  
&=& \int_{|z-\xi| = r(\xi)/2 }  \RR_\Lambda(z)    (I-\pi_1)  
\\
&=& \int_{|z-\xi| = r(\xi)/2 } \UU(z) \, ( I - \VV(z))^{-1}     (I-\pi_1)   =  0,
\eean
where in the last line we have used that $I-\VV(z) : X_0 \to Y_0$ is bijective for any $z \in B(\xi,r(\xi))$ so that 
$ ( I - \VV(z))^{-1}     (I-\pi_1) $ is well defined and analytic on $B(\xi,r(\xi))$. We deduce that 
  $$
\Pi_{\Lambda,\xi} = \Pi_{\Lambda,\xi} \pi_1,
$$  
and then 
$$
\hbox{dim} R(\Pi_{\Lambda,\xi} ) \le \hbox{dim} R( \pi_1 ) =N. 
$$
That precisely means $\xi \in  \Sigma_d(\Lambda)$. 

 \medskip\noindent
{\sl  Step 3. } The fact that $S_\Lambda$ satisfies the growth estimate \eqref{bddSlambda} is an immediate consequence of Theorem~\ref{theo:SpectralMapping}.
 \qed

%---------------------------------------------------------------------------------------------------------------------------------------------------------------------

\subsection{A quantified version of the Weyl's Theorem}

 We present now a quantified version of Weyl's Theorem~\ref{theo:GalWeyl}.

\begin{theo}  \label{theo:QuantWeyl}
Consider a Banach space $X$, the  generator $\Lambda$ of a  semigroup $S_\Lambda(t) = e^{t\Lambda} $ on $X$,
 a real number $a^* \in \R$ and
assume that there exist two operators $\ \AA, \BB \in  \mathscr{C}(X)$ such that 
  $\Lambda =  \AA + \BB$ and the  hypothesises {\bf (H1)}, {\bf (H2)} and {\bf (H3)} of Theorem~\ref{theo:GalWeyl} are met. 
Assume furthermore that 
\beqn\label{eq:hypQW}
\forall \, a > a^*, \quad
\int_0^\infty \| \AA  S_{\BB}(t) \|_{\BBB(X,Y)} \,  e^{-at} \, dt \le C_a,
\eeqn
for some constant $C_a \in \R_+$, 
and there exists a sequence of $N$ dimensional range increasing projectors $\pi_N$ and a sequence of positive real numbers $\eps_N \to 0$ 
such that 
\beqn\label{eq:YpiN}
\forall \, f \in Y \quad \| \pi^\perp_N f \|_X \le \eps_N \| f \|_{Y}. 
\eeqn
For any $a > a^*$, there exists an  integer $n^*$ (which depends on a constructive way on $a$, $(\pi_N)$, $(\eps_N)$ and the constants involved in the assumptions {\bf (H1)}, {\bf (H2)},  {\bf (H3)} and \eqref{eq:hypQW}) such that 
\beqn\label{eq:quantWeyl1}
\sharp (\Sigma(\Lambda) \cap \Delta_a) \le n^*, \quad \hbox{dim} \, R(\Pi_{\Lambda,a}) \le n^* .
\eeqn
We assume moreover {\bf (H1$''$)} and  $\| (S_\BB \AA)^{(*n)} \|_{\BBB(X,D(\Lambda^\zeta))} e^{-at} \in L^1(0,\infty)$ for any $a > a^*$. Then, for any $a > a^*$,  there exists 
a constant $C'_a$ such that for any Jordan  basis $(g_{i,j})$ associated to the eigenspace $R\Pi_{\Lambda,a}$
there holds
\beqn\label{eq:quantWeyl2}
\| g_{i,j} \|_X = 1, \quad \| g_{i,j} \|_{D(\Lambda^\zeta) \cap Y} \le C'_a.
\eeqn

\end{theo}

\noindent
{\sl Proof of Theorem~\ref{theo:QuantWeyl}. } {\sl Step 1. } Let us fix $a > a^*$ and let us define for any $z \in \Delta_a$ the compact perturbation of the identity 
$$
\Phi(z) := I + \AA \, R_{\BB}(z) : X \to X.
$$
On the one hand, because of \eqref{eq:hypQW}, we know that there exists a constant  $C_a$ such that 
$$
\forall \, z \in \Delta_a \qquad \| \AA \, \RR_\BB(z) \|_{\BBB(X,Y)} \le C_a,
$$
and then 
$$
\| \pi^\perp_N  \AA \, \RR_\BB(z) \|_{\BBB(X)} \le \eps_N \, C_a < 1
$$
for $N$ large enough. We deduce from the above smallness condition that 
$$
I + \pi^\perp_N  \AA \, \RR_\BB(z) : R\pi^\perp_N \to R\pi^\perp_N 
$$
is a an isomorphism for any $z \in \Delta_a$. On the other hand, thanks to the Fredholm alternative, it is clear that $\Phi(z)$ is invertible if, and only if, 
$\pi_N \Phi(z) $ has maximal rank $N$.

 \smallskip\noindent
 {\sl Step 2. } For a given basis $(g_1, ... , g_N)$ of $R\pi_N$ we denote by $\pi_{N,i}$ a projection on  $\C \, g_i$, $1 \le i \le N$,
and we define   $\Phi_{N,i}(z) := \pi_{N,i} \Phi(z)$. 
For a given $i \in \{1, ..., N \}$ and if $z_i \in \Delta_a$ satisfies  $\Phi_{N,i}(z_i) =0$ we have 
$$
\pi_{N,i} + \pi_{N,i} \, \AA \, \RR_\BB(z_i) = 0
$$
and
$$
 \Phi_{N,i}(z) = \pi_{N,i} +\pi_{N,i} \, \AA \sum_{n=0}^\infty \RR_\BB(z_i)^{n+1} \, (z-z_i)^n \quad \forall \, z \in \Delta_a \cap B(z_i,r),
$$
with $r := \|  \RR_\BB(z_i) \|^{-1} \ge C_a^{-1} > 0$. From the two last equations,  we deduce 
$$
  \Phi_{N,i}(z) = -   \pi_{N,i}   \sum_{n=1}^\infty \RR_\BB(z_i)^n \, (z-z_i)^n \quad \forall \, z \in \Delta_a \cap B(z_i,r_i)
$$
and we observe that  rank$\,   \Phi_{N,i}(z) = 1$ for any $z \in B(z_i,r)$, $z \not= z_i$. 
As a consequence, in any ball $B $  of radius $r$, rank$\, \pi_N   \Phi(z) = N$ for any $z \in B \cap \Delta_a$ except at most $N$ points $z_1$, ..., $z_N \in B$, and the total dimension of the ``defect of surjectivity"   $\hbox{codim} \,  R(\pi_N\Phi(z_i)) = \hbox{dim} \, \Pi_{\Lambda,z_i} $ is at most $N$. 
Covering the region $\bar\Delta_a \cap B(0,M)$ by $(1+2M/r)^2$ balls of radius $r \in (C_a^{-1}, M)$, we see that \eqref{eq:quantWeyl1} holds with 
$n^*  := (1+ 2MC_a)^2 N$. 

\smallskip\noindent
 {\sl Step 3. }  Consider a Jordan basis $(g_{j,m})$  associated to an eigenvalue $\xi \in \Sigma_d(\Lambda) \cap \Delta_a$, and then defined by 
$$
\Lambda g_{j,m} = \xi g_{j,m} + g_{j,m-1}.
$$
We write 
$$
g_{j,m} =  \RR_\BB(\xi) g_{j,m-1} - \RR_\BB(\xi) \AA g_{j,m} , 
$$
and iterating the formula
$$
g_{j,m} =  \sum_{\ell = 0}^{L-1} (-1)^\ell \RR_\BB(\xi) (\AA \, \RR_\BB(\xi) )^\ell  g_{j,m-1} +  (-1)^L (\RR_\BB(\xi) \AA)^n g_{j,m} , 
$$
from which we easily conclude that  \eqref{eq:quantWeyl2} 
holds.  \qed 
 
\smallskip

%------------------------------------------------------------------------------------------------------------------------------------------------------------

%\subsection{Application of  Weyl's and spectral mapping theorems to the growth-fragmentation equation}
\section{Semigroup Weyl's theorem for the growth-fragmentation equations}
\label{sec:SGWgrowthfrag}
%------------------------------------------------------------------------------------------------------------------------------------------------------------

%\subsection{Application of  Weyl's and spectral mapping theorems to the growth-fragmentation equation}
\subsection{Equal mitosis and smooth cell-division equations}
\label{subsec:CellDiv}
In this paragraph we are concerned with the equal mitosis equation \eqref{eq:mitose} and the cell division equation \eqref{eq:celldiv} with smooth offspring distribution, so that  we consider the operator
$$
\Lambda f (x) :=  -  \frac{\partial}{\partial x} f(x) - K(x) f(x) + (\FF^+ f)(x)
$$
where $K$ satisfies \eqref{eq:hypKmitose1} and \eqref{eq:hypKmitose2} and the gain part $\FF^+$ of the fragmentation operator is defined by
\beqn\label{eq:defFF+}
\FF^+ f (x) := \int_x^\infty k(y,x) \, f(y) \, dy
\eeqn
with $k$ satisfying \eqref{eq:k=Kwp} and $\wp = 2 \delta_{1/2}$ (in the equal mitosis case) or $\wp$ is a function which satisfies \eqref{eq:intwp=1} and \eqref{eq:wpsmooth} (in the case of the smooth cell-division equation).  We recall that 
%that except when $K$ vanishes on a neighbourhood of $0$ 
the operator is complemented with a boundary condition \eqref{eq:BoundaryCondition}. 
%\beqn\label{eq:condx0=0}
%f(x_0/2) = 0.
%\eeqn

\smallskip
We fix $\alpha > \alpha^* $, with $\alpha^* > 1$ defined thanks to Equation \eqref{eq:defalpha*}, and we set $K_0' := K_0 - \wp_\alpha K_1 >0$.
When $K(0) = 0$ or $\wp = 2\delta_{1/2}$, we then define the critical abscissa $a^* \in \R $  by 
\bean
&& a^*:=  - K'_0   \,\,\, \hbox{if} \,\,\, \gamma = 0,   
\qquad a^* :=  - \infty \,\,\, \hbox{if} \,\,\, \gamma >0. 
\eean
Observe that when $K(0) > 0$, the positivity conditions \eqref{eq:hypKmitose1} and  \eqref{eq:hypKmitose2} imply that there exists a constant $K_* > 0$
such that 
\beqn\label{eq:hypKmitose3}
K(x) \ge K_* \qquad \forall \, x \ge 0.
\eeqn
When $K(0) > 0$ and $\wp$ satisfies \eqref{eq:wpsmooth}, we then define the critical abscissa $a^* \in \R $  by
\bean
&& a^*:=  - \min(K'_0,K_*)   \,\,\, \hbox{if} \,\,\, \gamma = 0, 
\qquad a^* :=  - K_* \,\,\, \hbox{if} \,\,\, \gamma >0. 
\eean

We perform the spectral analysis of $\Lambda$ in the Banach space 
$X = L^1_\alpha$ defined at the beginning of Section~\ref{sss:MainResult}. We also define, for
later reference, the Sobolev spaces
$$
W_\alpha^{1,1} = \{ f \in L^1_\alpha, \,\, \partial_x f \in L^1_\alpha \}, \quad
\dot W^{1,1}_{\alpha} := \{ f \in L^1_{loc}; \,\, \partial_x f \in \dot L^1_{\alpha} \}, \quad \alpha \in \R.
$$

It is worth emphasizing that we classically have 
(and that is also a straightforward consequence of the lemmas which follow) 
$$
D(\Lambda) = \{ f \in L^1_{\alpha+\gamma}, \,\, \partial_x f \in L^1_\alpha, \,\, f(0) = 0 \}.
$$
As a first step in the proof of Theorem~\ref{theo:Frag} we have 

\begin{prop}\label{prop:GrFrSpectralMapping} Under the above assumptions and definitions of $\alpha^* > 1$ and $a^* < 0$, the conclusion {\bf (2)}, 
and then \eqref{eq:localizSpectrL}  and {\bf (1)}, 
of Theorem~\ref{theo:GalWeyl}  holds for the cell-division semigroup in $L^1_\alpha$ for any $\alpha > \alpha^*$ and for any $a > a^*$.  
Moreover the conclusions of Theorem~\ref{theo:QuantWeyl} hold under the additional assumptions that $K(0) = 0$ and $\wp$ satisfies the smoothness 
condition \eqref{eq:wpsmooth}.
\end{prop}

In order to establish Proposition~\ref{prop:GrFrSpectralMapping} we will introduce an adequate splitting $\Lambda = \AA + \BB$ and 
we prove that $\AA$ and $\BB$ satisfy  conditions {\bf (H1)},  {\bf (H2)} and {\bf (H3)} (or one of the ``prime" variants of them) for $n=1$ or $2$
as a consequence of the series of technical Lemmas~\ref{lem:BornesAA}, \ref{lem:BorneBB}, \ref{lem:OpLambda}, \ref{lem:GrFrH2} and \ref{lem:regMitose} below.

\smallskip
Taking a  real number $R \in [1,\infty)$ to be chosen later, we define 
$$
K_{R}  := K \, \chi_R 
\quad
K^c_R  := K \, \chi^c_R, 
\quad
k_R = K_R \, \kappa,
\quad
k_R^c = K_R \, \kappa
$$
where $\chi_R (x) = \chi(x/R)$, $\chi^c_R (x) = \chi^c(x/R)$,  $\chi$ being the Lipschitz function defined 
on $\R_+$ by $\chi(0) = 1$, $\chi' = - {\bf 1}_{[1,2]}$ and $\chi^c = 1 - \chi$,  
as well as 
$$
\FF^{+}_R = \int_x^\infty k_R(y,x) \, f(y) \, dy, \quad
\FF^{+,c}_R = \int_x^\infty k^c_R(y,x) \, f(y) \, dy, 
$$
and then 
\bean
\AA = \AA_R = \FF^{+}_R, \quad \BB = \BB_{R}=  -  \frac{\partial}{\partial x}  - K(x) + \FF^{+,c}_R
\eean
so that $\Lambda =  \AA + \BB$. 

\begin{lem}\label{lem:BornesAA} 

(1) For any $0 \le \alpha \le 1 \le \beta$ there holds $\AA  \in \BBB(\dot L^1_\alpha, \dot L^1_{\alpha} \cap \dot L^1_\beta)$. 

(2) For the mitosis operator, there holds $\AA  \in \BBB(W^{1,1} , W^{1,1}_\alpha)$ for any $\alpha \ge 0$. 
In particular  $\AA  \in \BBB (D(\Lambda^\eta))$ for   $\eta = 0,1$.

(3) Under assumption \eqref{eq:wpsmooth} on $\wp$  there holds $\AA \in \BBB(\dot L^1_{\beta-1}, \dot W^{1,1}_\beta)$ for any 
$\beta \ge 0$.

(4) Under assumption \eqref{eq:wpsmooth} on $\wp$  and the additional assumption  $K(0) = 0$,  there holds $\AA \in \BBB(L^1_\alpha, W^{1,1}_{\alpha})$ for any 
$\alpha \ge 0$. 

\end{lem}

\noindent
{\sl Proof of Lemma \ref{lem:BornesAA}.}
We split the proof into four steps. 

\noindent
{\sl Step 1. } Fix $0 \le f \in \dot L^1_\alpha$ as well as $\alpha' \ge \alpha$. Recalling from \eqref{eq:hypKmitose2} that $K(x) \le K_1 \langle x \rangle^\gamma$, $\langle x \rangle := (1+x^2)^{\gamma/2} $, we compute 
\bean 
\| \AA f \|_{\dot L^1_{\alpha'}} &=& \int_0^\infty f(x) \, \int_0^x k_R(x,y) \, y^{\alpha'} \, dy \, dx
\\
&=& \wp_{\alpha'} \int_0^\infty  f(x) \, K_R(x) \, x^{\alpha'} \, dx
\\
&\le&  \wp_{\alpha'} \,  K_1 \, \langle R \rangle^{\gamma+\alpha'-\alpha}   \int_0^{2R}  f(x)  \, x^{\alpha} \, dx,
\eean
so that $\AA \in \BBB(\dot L^1_\alpha, \dot L^1_{\alpha'})$.

\noindent
{\sl Step 2. } For the mitosis operator, we have 
$$
\partial_x \AA f = 8 \, (\partial_x K_R)(2x) \, f(2x) + 8 \, K_R(2x) \, (\partial_xf) (2x), 
$$
and similar estimates as above leads to the bound 
 $$
 \| \partial_x \AA_R \, f \|_{L^1_{\alpha}} \le 4 \, \langle R \rangle^{\alpha} \,  \| K  \|_{W^{1,\infty}(0,2R)} \, \| f \|_{W^{1,1}}
 $$
for any  $\alpha \ge 0$ and $f \in W^{1,1}$, from which (2) follows. 

\noindent
{\sl Step 3. } Under the regularity assumption \eqref{eq:wpsmooth}, there holds
$$
\partial_x \AA f =  \int_x^\infty K_R(y) \, \wp'(x/y) \, y^{-2} \, f(y) \, dy - K_R(x) \, \wp(1) \, f(x)/x,
$$
so that 
\bean
\| \partial_x \AA f  \|_{\dot L^1_{\beta}} \le  \bigl\{ \wp(1) + \wp'_{\beta} \bigr\} 
K_1 \, \langle R \rangle^{\gamma}   \int_0^{2R}  f(x)  \, x^{\beta-1} \, dx,
\eean
and we conclude with $\AA \in \BBB(\dot L^1_{\beta-1}, \dot W^{1,1}_\beta)$. \qed

\noindent
{\sl Step 4. } With the assumption of point (4) and recalling that $K$ is $C^1$, there holds $\| K(x)/x \|_{L^\infty(0,2R)}$ for any $R \in (0,\infty)$ and then 
from the above expression of $\partial_x \AA f$, we get 
\bean
\| \partial_x \AA f  \|_{L^1_{\alpha}}  &\le& \bigl\{ \wp(1) + \wp'_0 \bigr\} \, \| K(x)/x \|_{L^\infty(0,2R)}
\|   f  \|_{L^1_{\alpha}} 
\eean
for any $f \in L^1_\alpha$ and $\alpha \ge 0$. 
\qed

\medskip

\begin{lem}\label{lem:BorneBB} For any $a > a^*$ there exists $R^*(a) > 0$ such that the operator $\BB$ is $a$-hypodissipative in $L^1_\alpha$ for any $R \in (R^*(a),\infty)$. 

\end{lem}

\noindent
{\sl Proof of Lemma  \ref{lem:BorneBB}.}   We introduce the primitive functions
\beqn\label{def:primitivK} 
\KK(z) := \int_0^z K(u) \, du, \quad \KK(z_1,z_2) := \KK(z_2) - \KK(z_1)
\eeqn
and,  for any given $a \in (a^*,0]$, we define the space
$$
\EE := L^1(\phi), \quad \phi(x) :=  \phi_0(x) \, {\bf 1}_{x \le x_2} +  \phi_\infty(x) \, {\bf 1}_{x \ge x_2}, 
$$ 
where
$$
\phi_0 (x) := \frac{e^{\KK(x) + a \, x}}{e^{\KK(x_2) + a  \, x_2}},
\quad
\phi_\infty(x) := { x^\alpha \over x^\alpha_2}, 
$$
and where  
\bear\label{def:x1gamma=0}
x_2 &:=&  \max(1,x_1, 2\alpha/(a+K'_0)) \quad \hbox{if} \,\,\, \gamma = 0;
\\ \label{def:x1gamma>0}
x_2 &:=& \max(1,x_1,  [3\alpha/(x_1K'_0)]^{1/\gamma}, [-3a/ K'_0]^{1/\gamma} ) \quad \hbox{if} \,\,\, \gamma >  0.
\eear
We recall that $x_1$ is defined in \eqref{eq:hypKmitose2} and $K'_0$  is defined at the beginning of Section~\ref{subsec:CellDiv}.

\smallskip
Consider  $f \in D(\Lambda)$ a real-valued function and let us show that for any $a > a^*$ and any $R > R^*(a)$, $R^*(a)$ to be chosen later, 
\beqn\label{eq:B0MitoseDissip}
  \int_0^\infty  \sign(f(x)) \, \BB  f(x)\, \phi (x) \,dx
  \leq
  a \|f\|_{\EE}. 
\eeqn
Since the case of complex-valued functions can be handled in a similar way and since the norm  $\| \cdot \|_\EE$  is clearly equivalent to   $\| \cdot \|_{L^1_\alpha}$,  that will end the proof. 
 
\smallskip
On the one hand, we have 
\bean
\| \FF^{+,c} f \|_{L^1(\phi_0 \, {\bf 1}_{x \le x_2})} 
&=&  \int_R^\infty K^c_R(x) \,  f(x) \int_0^{x \wedge x_2} \wp(y/x) \, \phi_0 (y) \, {dy \over x} \, dx
\\
&\le& \eta(x_2/R)  \int_R^\infty K^c_R(x) \,  f(x) \, dx ,
\eean
with 
$$
 \eta(u) := \bigl( \sup_{[0,x_2]} \phi_0 \bigr) \Bigl( \int_0^{u} \wp(z) \, dz \Bigr),
 $$
so that, performing one integration by part, we calculate  
  \bear \nonumber
&&\int_0^{x_2}  \sign(f(x)) \, \BB  f(x)\, \phi_0 (x) \,dx =
 \\ \nonumber
&&\quad = 
   \int_0^{x_2}  \{   -K(x)|f(x)| - \partial_x |f(x)| \, \} \,  \phi_0 (x) \,dx + \int_0^{x_2}  \sign(f) \, (\FF^{+,c} f ) \, \phi_0   \,dx
 \\ \label{eq:B0MitoseDissip3}
&&\quad \le \ 
 -|f(x_2)| + a  \, \int_0^{x_2} |f(x)| \, \phi_0 (x) \,dx + \eta(x_2/R) \int_R^\infty K^c_R(x) \,  f(x) \, dx .
  \eear
On the other hand, performing one integration by part again, we compute 
    \bear  \nonumber
 &&\quad \int_{x_2}^\infty  \sign(f(x)) \, \BB_R f(x)\, \phi_\infty (x) \,dx = 
 \\ \nonumber
&&\quad \le \ 
 |f(x_2)| +   \int_{x_2}^\infty  |f(x)|  \,  \{ - K \phi_\infty +  \partial_x \phi_\infty  +  K^c_R(x) \, \int_0^x \wp(\frac yx) \, \phi_\infty(y) \, \frac {dy}x   \}  \,dx
 \\ \nonumber
&&\quad \le \ 
 |f(x_2)| +   \int_{x_2}^\infty  |f(x)|  \,  [ (\wp_\alpha K_1   - K_0) \, x^\gamma + \alpha/x ] \,  \phi_\infty (x)  \,dx
 \\ \label{eq:B0MitoseDissip4}
&&\quad\le \ 
 |f(x_2)| +  \int_{x_2}^\infty  |f(x)|  \, \Bigl\{ a + \theta \, x^\gamma \Bigr\}   \phi_\infty (x)  \,dx,
  \eear
  with $\theta := - (a+K'_0)/2$ if $\gamma = 0$ and $\theta := - K'_0/3$ if $\gamma > 0$  
thanks to the choice of $x_2$. Gathering \eqref{eq:B0MitoseDissip3} and \eqref{eq:B0MitoseDissip4}, and taking $R^*$ large enough in such a way that 
$\eta(x_2/R^*) \, K_1 + \theta \le 0$, we get that $\BB_R$ is $a$-dissipative in $\EE$.  \qed

\begin{lem}\label{lem:OpLambda}  
The operator  $\Lambda$ generates a $C_0$-semigroup on $L^1_\alpha$ . 

\end{lem}

\noindent
{\sl Proof of Lemma \ref{lem:OpLambda}.}  Thanks to Lemmas  \ref{lem:BornesAA} and \ref{lem:BorneBB},  we have that $\Lambda$ is $b$-dissipative in $\EE$ with $b := \| \AA_R \|_{\BBB(\EE)} + a$. 
On the other hand, one can show that $R(\Lambda - b) = X$ for $b$ large enough and conclude thanks to Lumer-Phillips Theorem (\cite{MR0132403} or \cite[Theorem I.4.3]{Pazy}) that $\Lambda$ generates a $C_0$-semigroup.  Equivalently,  one can argue as in \cite[Proof of Theorem 3.2]{EMRR} 
by introducing an approximation sequence of bounded total fragmentation rates $(K^n)$ and proving that for any fixed initial datum $f_0 \in X$ the associated sequence of solutions $(f^n)$ (constructed by a mere Banach fixed point Theorem in $C([0,T];\EE)$, $\forall \, T > 0$) is a Cauchy sequence.   That establishes that for any $f_0 \in \EE$ there exists a unique solution $f \in C(\R_+; \EE)$ to the Cauchy problem associated to the operator $\Lambda$ and then that $\Lambda$ generates a $C_0$-semigroup.  
\qed

\medskip 
We define 
$$
Y_r := \{ f \in W^{r,1}_\alpha \cap L^1_{\alpha+\gamma+1} (\R); \,\, \hbox{supp} f \subset [0,\infty) \}, \quad r \in [0,1],
$$
as a family of interpolating spaces between $Y_0 = L^1_{\alpha+\gamma+1} \subset X$ and $Y_1 \subset D(\Lambda)$. 

\begin{lem}\label{lem:GrFrH2} If  $\wp$ satisfies \eqref{eq:wpsmooth} and $K(0) = 0$,  for any $a> a^*$ there exists a constant $C_a$ such that 
\beqn\label{eq:U2wpK0=0}
\| \AA S_\BB  (t) \|_{\BBB(X,Y_1)} \le C_a \, e^{a t}. 
\eeqn
If  $\wp$ satisfies \eqref{eq:wpsmooth} and $K(0) \not= 0$,  for any $a> a^*$ and any $r \in [0,1)$ there exists a constant $C_{a,r}$ such that 
 \beqn\label{eq:U2wpK0not0}
\| (\AA S_\BB)^{(*2)} (t) \|_{\BBB(X,Y_r)} \le C_{a,r} \, e^{a t}. 
\eeqn
 \end{lem}

\noindent
{\sl Proof of Lemma \ref{lem:GrFrH2}.} 
Extending by $0$ a function $g \in L^1(\R_+)$ outside of $[0,\infty)$, we may identify $L^1(\R_+) = \{ g\in L^1(\R); \,\, \hbox{supp} \, g \subset [0,\infty) \} $
and extending $k=k(x,y)$ by $0$ outside of $\{ (x,y); \,\, 0 \le y \le x \}$, we may consider $S_\BB(t)$ and $\AA$ as operators acting on $L^1(\R)$ which preserve the support $[0,\infty)$ (they then also act on $L^1(\R_+)$).  

\medskip\noindent
 {\sl Step 1. } Assume first $K(0) = 0$. From Lemma~\ref{lem:BornesAA}-(4) and  Lemma~\ref{lem:BorneBB}, for any 
$f \in L^1_\alpha$, we get 
$$
\| \partial_x(\AA S_\BB(t) f ) \|_{L^1_\alpha} \le C \, \|  S_\BB(t) f   \|_{L^1_\alpha} \le C_a \, e^{at} \|  f   \|_{L^1_\alpha}. 
$$
 From Lemma~\ref{lem:BornesAA}-(1) and  Lemma~\ref{lem:BorneBB} we get a 
similar estimate on the quantity  $\| \AA S_\BB(t) \|_{\BBB(L^1_\alpha,L^1_{\alpha+\gamma+1})}$ and that ends the proof of \eqref{eq:U2wpK0=0}.

\smallskip\noindent
{\sl Step 2. } We assume now $K(0) \not=0$. We introduce the notation 
$$
\BB_0  f (x) := -  \frac{\partial}{\partial x} f(x) - K(x) f(x)
$$
and then the shorthands $\AA^c :=  \FF^{+,c}$, $U := \AA S_\BB$, $U^c := \AA^cS_\BB$,  $U_0 := \AA S_{\BB_0}$ and $U^c_0 := \AA^c  S_{\BB_0}$,
where $S_\BB$ (resp. $S_{\BB_0}$) is the semigroup associated to the generator $\BB$ (resp. $\BB_0$) complemented with the boundary condition \eqref{eq:BoundaryCondition}.  Thanks to the Duhamel formula
\beqn\label{eq:DuhamelF1}
S_\BB  =  S_{\BB_0}  + S_{\BB_0} *  \AA^c S_\BB,
\eeqn
we have 
\beqn\label{eq:U=U0+U0*Uc}
U  =  U_0  + U_0 * U^c. 
\eeqn
From the explicit representation formula
\beqn\label{eq:defSB0}
(S_{\BB_0}(t) f)(x) =  e^{- \KK(x-t,x)} f(x-t) \quad \forall \, t \ge 0, \,\, \forall \, x \in \R, 
\eeqn
where $\KK$ is the   primitive function defined in \eqref{def:primitivK}, we get  
$$
U_0 (t) f (x) =  \int_x^\infty k_R   (y,x) \, e^{- \KK(y-t,y)} f(y-t)  \, dy,
$$
so that 
\bean
\partial_x [ U_0 (t) f (x)] 
&=& \int_x^\infty K_R(y) \, \wp'(x/y) \, y^{-2} \, e^{- \KK(y-t,y)} f(y-t) \, dy
\\
&& - {K_R (x) \over x} \, \wp(1) \, e^{- \KK(x-t,x)} f(x-t) .
\eean
Using  that $\KK(y-t,y) \ge K_* \, t$ for any $y \ge t$, we have 
$$
 e^{- \KK(y-t,y)} |f(y-t)| \le e^{a t} \,  |f(y-t)|  \quad \forall \, t \ge 0, \,\, \forall \, y \in \R, 
 $$
 and we deduce 
\bean
\| \partial_x [ U_0 (t) f ] \|_{  L^1_\alpha} 
&\le&  \wp(1) \int_0^{2R} K_R(x) \, \langle x \rangle^\alpha {|f (x-t)| \over (x-t)+t} \, dx \, e^{at} 
\\
&&+ \int_0^{2R} K_R(x) \, \int_0^x |\wp'(y/x)| \,  \langle y \rangle^\alpha { dy \over x}  {|f (x-t)| \over (x-t)+t} \, dx  \, e^{at} 
\\
&\le&  \Bigl( \wp(1) +  \wp'_0 \Bigr) \, K_1 \, \langle 2R \rangle^{\gamma+\alpha}  {  e^{at}  \over t} \int_0^{2R} | f(x)|  dx.
\eean
Since clearly (from the explicit representation formula \eqref{eq:defSB0} for instance) the same estimate holds for 
$\|  U_0 (t)   \|_{ \BBB(L^1, L^1_{\alpha+\gamma+1})}$ we get by interpolation and for any $r \in [0,1]$
$$
\|   U_0 (t)   \|_{\BBB(L^1,Y_r)} \le C_a \, {e^{at} \over t^r}  .
$$
Thanks to the identity \eqref{eq:U=U0+U0*Uc}, the fact that $\AA^c \in \BBB(L^1_\gamma,L^1)$ and Lemma~\ref{lem:BorneBB}, we deduce for $r \in [0,1)$
\bean
\| U(t) f \|_{Y_r} 
&\le& C_a \, {e^{at} \over t^r}  \, \| f \|_{L^1} + \int_0^t C_a \, {e^{a(t-s)} \over (t-s)^r } \, \| \AA^c S_{\BB}(s) f \|_{L^1}
\\
&\le& C_{a,r} \, {e^{at} \over t^r }  \, \| f \|_{L^1_{\max(\alpha,\gamma)}} . 
\eean
On the other hand, from Lemma~\ref{lem:BornesAA}-(1) and Lemma~\ref{lem:BorneBB}, we know that 
$$
\| U (t) f \|_{L^1_{\alpha+\gamma}} 
\le \| \AA \|_{\BBB(L^1_\alpha,L^1_{\alpha+\gamma})} \, \| S_\BB(t) f \|_{L^1_\alpha} 
\le C_a e^{at} \| f \|_{L^1_\alpha} .
$$
These two estimates together imply \eqref{eq:U2wpK0not0}. \qed

\begin{lem}\label{lem:regMitose} For the equal mitosis equation,  there holds
\beqn\label{eq:regMitose}
\|(\AA \, S_\BB)^{(*2)} (t) \|_{\BBB(X, Y_1)} \le C \, t   \quad \forall \, t \ge 0
\eeqn
 for some constant $C$ which only depends on $K$ through its norm $\|K \|_{W^{1,\infty}(0,R)}$, 
where $R$ is defined in Lemma~\ref{lem:BorneBB}. 
 
\end{lem} 

\noindent
{\sl Proof of Lemma \ref{lem:regMitose}.}  
Thanks to the Duhamel formula \eqref{eq:DuhamelF1} and   the iterated Duhamel formula
$$
S_\BB  =  S_{\BB_0}  + S_{\BB_0} *  \AA^c S_{\BB_0}  + S_{\BB_0} *  \AA^c S_{\BB_0} *   \AA^c S_\BB  ,
$$
we have 
$$
U  =  U_0  + U_0 * U^c \quad\hbox{and}\quad U  =  U_0  + U_0 * U_0^c + U_0 * U^c_0 * U^c,
$$
from which we finally deduce 
\beqn\label{eq:U*2-1}
U^{*2} 
=  U_0^{*2} + U_0^{*2} * U^c +   U_0 * U_0^c * U  + U_0 * U^c_0 * U^c * U.
\eeqn
From the explicit representation formula \eqref{eq:defSB0}, we get  
\bean
(U_0(t) f)(x) &=& 4 \, K_R(2x) \, e^{- \KK(2x-t,2x)} f(2x-t) 
\eean
as well as 
\bean
(U^c_0(t) f)(x) &=& 4 \, K^c_R(2x) \, e^{- \KK(2x-t,2x)} f(2x-t) .
\eean
We then easily compute 
\bear \nonumber
&&(U_0^{*2} (t) f )(x) 
= \int_0^t (U_0(t-s)  U_0(s) f )(x) \, ds 
\\  \nonumber
&&= 4 \, K_R(2x) \int_0^t  e^{- \KK(2x-t+s,2x)} (U_0(s) f)(2 x - t + s)  \, ds 
\\  \nonumber
&&=  16 \, K_R(2x)  \int_0^t  K_R(2(2 x - t + s)) \, 
\\  \nonumber
&&\qquad e^{- \KK(2x-t+s,2x) - \KK(4x -2t+s,2(2 x - t + s))} f(4x -2t+s)  \, ds 
\\ \label{eq:U*2-2}
&&=  16 \, K_R(2x)  \int_{u_0}^{u_1}  K_R(2u-4x+2t)  \, e^{- \Theta(u)} f(u)  \, du 
%\\
%&&=  16 \, K_R(2x)  \int_{(4x-2t)\wedge(2x-t+R/2)}^{(4x-t)\wedge(2x-t+R/2)}  K(2u-4x+2t)  \, e^{- \KK(u-2x+t,2x) - \KK(u,2u-4x+2t)} f(u)  \, du 
\eear
with $u_0 := (4x-2t)\wedge(2x-t+R/2)$, $u_1 := (4x-t)\wedge(2x-t+R/2)$ and $\Theta (u) := \KK(u-2x+t,2x) + \KK(u,2u-4x+2t) $. 
Similarly, we have
\beqn\label{eq:U*2-3}
(U_0 * U_0^c (t) f )(x) 
=  16 \, K_R(2x)  \int_{u_0}^{u_1}  K^c_R(2u-4x+2t)  \, e^{- \Theta(u)} f(u)  \, du.
\eeqn
 The two last terms are clearly more regular (in $x$) than the initial function $f$. Moreover, extending by $0$ the function $f$ outside of $[0,\infty)$, we see thanks to \eqref{eq:defSB0} that $S_{\BB_0}(t) f$, and then $U_0^{*2} (t) f $ and  $U_0 * U_0^c (t) f $, are  well defined as functions in $L^1(\R)$.  
 Using the lower bound $\Theta \ge 0$ and performing some elementary computations, we easily get from \eqref{eq:U*2-1}, \eqref{eq:U*2-2} and \eqref{eq:U*2-3} the estimate 
 $$
 \| \partial_x [ (\AA S_B)^{(*2)}(t) f ] \|_{L^1(\R) } \le C_R \, t  \, \| f \|_{L^1}, \quad \forall \, t \ge 0,
 $$
 which ends the proof of \eqref{eq:regMitose} since then $ (\AA S_B)^{(*2)}(t) f \in C(\R)$ and \\ supp$ \, (\AA S_B)^{(*2)}(t) f \subset [0,\infty]$
 imply $(\AA S_B)^{(*2)}(t) f (0) = 0$. \qed
  
\medskip
With all the estimates established in the above lemmas, we are able to present the 

\medskip\noindent
{\sl Proof of Proposition~\ref{prop:GrFrSpectralMapping}. } We just have to explain why the hypothesis of Theorem~\ref{theo:GalWeyl} are 
satisfied in each cases. Hypothesis {\bf (H1)} is an immediate consequence of Lemma~\ref{lem:BornesAA}-(1) and Lemma~\ref{lem:BorneBB}
together with Remark~\ref{rem:theoEquivSGX}-(b). For that last claim we use that $Y_r \subset D(\Lambda^{r'})$ for any $0 \le r' < r < 1$ thanks to the classical interpolation theory, see \eqref{eq:XalphaSubsetTildeXalpha}, \cite{MR0206716} and \cite{MR0165343}. 

 \smallskip
 In the case when $\wp$ satisfies \eqref{eq:wpsmooth}, Lemma~\ref{lem:GrFrH2} 
 and Remark~\ref{rem:theoEquivSGX}-(h) imply that hypothesises  {\bf (H2)} and  {\bf (H3)} are met with $n=1$, $\zeta = 1$ and  $\zeta' = 0$ in the case $K(0) = 0$ and are met with $n=2$, 
 $\zeta \in (0,1)$ and $\zeta' = 0$ in the case $K(0) \not= 0$. Also notice that the additional assumptions of Theorem~\ref{theo:QuantWeyl} are met in the case that $K(0) = 0$.
 
 \smallskip
 For the equal mitosis equation,  Lemma~\ref{lem:BornesAA}-(2) and Lemma~\ref{lem:BorneBB} imply that assumption {\bf (H1)} also holds with $\BBB(X)$ replaced by $\BBB(D(\Lambda))$ and $\BBB(X,L^1_{\alpha+1})$, while  Lemma~\ref{lem:regMitose} implies that {\bf (H2$''$)} holds with $n=2$, $\zeta = 1$, $\zeta' = 0$ and any $b > 0$, so that we can apply Corollary~\ref{cor:ThSM2} (assumption {\bf (H2) } 
 holds with $n$ large enough thanks to Lemma~\ref{lem:Tn}). Finally, hypothesis {\bf (H3)}
  follows from {\bf (H2)} and the fact that $D(\Lambda) \cap L^1_{\alpha+1} \subset X$ compactly. \qed

%---------------------------------------------------------------------------------------------------------------------------------------------------------------------

\subsection{Quantified semigroup Weyl's Theorem for the self-similar fragmentation equation}

In this paragraph we are concerned with the self-similar fragmentation equation \eqref{eq:fragSSV}, so that  we consider the operator
$$
\Lambda f (x) :=  - x \,  \frac{\partial}{\partial x} f(x) - 2 \, f(x) - K(x) f(x) + (\FF^+ f)(x)
$$
where $K(x) = x^\gamma$ and   $\FF^+$ is defined in \eqref{eq:defFF+}
with $k$ satisfying \eqref{eq:k=Kwp},  \eqref{eq:intwp=1} and \eqref{eq:wpsmooth}.
 
We perform the spectral analysis of $\Lambda$ in the Banach space $X := \dot L^1_\alpha \cap \dot L^1_\beta$
for $0 \le \alpha < 1 < \beta < \infty$, where we recall that the homogeneous Lebesgue space $\dot L^1_\alpha$ has been defined at the beginning of Section~\ref{sss:MainResult},  and we set $a^* := \alpha-1 \in [-1,0)$.   It is worth emphasizing that  we classically have 
(that is again a consequence of the estimates established in the series of lemmas which  follow, moreover the inclusion $\supset$ is just straightforward and it is the only
inclusion that we will really use) 
$$
D(\Lambda) = \{ f \in \dot L^1_\alpha \cap \dot L^1_{\beta+\gamma}, \,\, \partial_x f \in \dot L^1_{\alpha+1} \cap \dot L^1_{\beta+1} \}. 
$$
As a first step in the proof of Theorem~\ref{theo:Frag} we have 

\begin{prop}\label{prop:FragSpectralMapping} Under the above assumptions and definitions,  %(of $\alpha,\beta$ and $a^*$), 
the conclusion {\bf (2)},  and then \eqref{eq:localizSpectrL}  and {\bf (1)},  of Theorem~\ref{theo:QuantWeyl} 
holds for the  self-similar fragmentation semigroup in $\dot L^1_\alpha \cap \dot L^1_\beta$ %for  %$0 \le \alpha < 1 < \beta < \infty$ and
 for any $a > a^*$.  
\end{prop}

In order to establish Proposition~\ref{prop:FragSpectralMapping} we introduce a suitable splitting $\Lambda = \AA + \BB$ and 
we prove that $\AA$ and $\BB$ satisfy  conditions {\bf (H1)},  {\bf (H2)}, {\bf (H3)} with $n=1$ as well as  \eqref{eq:hypQW} with $Y \subset D(\Lambda)$
as a consequence of the two technical Lemmas~\ref{lem:BornesAAFrag} and  \ref{lem:BorneBBFrag} below.

\smallskip
We introduce the following splitting inspired from \cite{MR2821681}.  With the notation of the previous section, we define
$$
k^s(x,y)  := k(x,y) \, \chi^c_\delta(x)  \, \chi^c_\eps(y)  \, \chi_{R} (y)
$$
for $0 < \eps \le \delta/2 \le 1$, $R \ge 2$ to be specified later, and then 
$$
k^c := k^{c,1} + k^{c,2} + k^{c,3}
$$
with $k^{c,1}(x,y) = k(x,y) \, \chi_\delta(x)$, $k^{c,2}(x,y) =  \chi^c_\delta(x)  \, \chi^c_R(y)$
and $k^{c,3}(x,y) =   k(x,y) \,$ $ \chi^c_\delta(x)     \, \chi _\eps(y)$.
We then define 
$$
\AA f (x) = \int_x^\infty k^s(y,x) \, f(y) \, dy, \quad
\AA^{c,i} f(x) = \int_x^\infty k^{c,i}(y,x) \, f(y) \, dy, 
$$
and then 
\bean
\AA^c = \AA^{c,1} + \AA^{c,2} + \AA^{c,3}, 
\quad
\BB_0 = - x \,  \frac{\partial}{\partial x}  - 2 - K(x),
\quad
\BB = \BB_0 + \AA^c 
\eean
so that $\Lambda =  \AA + \BB$.

\begin{lem}\label{lem:BornesAAFrag} 
For any $0 \le \alpha' \le 1 \le \beta'$, there holds $\AA  \in \BBB(\dot L^1_1, \dot W^{1,1}_{\alpha'} \cap \dot W^{1,1}_{\beta'})$
and $\AA^{c,i} \in \BBB(\dot L^1_{\alpha' + \gamma},\dot L^1_{\alpha'})$ for $i=1,2,3$, with 
\bean
 \| \AA^{c,1} f \|_{\dot L^1_{\alpha'}} &\le& \wp_{\alpha'} \, \delta^\gamma \,  \|   f \|_{\dot L^1_{\alpha' + \gamma}},
\\
 \| \AA^{c,2} f \|_{\dot L^1_{\alpha'}} &\le& \wp_{\alpha'} \, \int_R^\infty f(x) \, x^{\alpha'+\gamma} \, dx  ,
\\
 \| \AA^{c,3} f \|_{\dot L^1_{\alpha'}} &\le& \| \wp (z) \, z^{\alpha'} \|_{L^1(0,\eps/\delta)} \, 
\|   f \|_{\dot L^1_{\alpha' + \gamma}}. 
\eean
\end{lem}

\noindent
{\sl Proof of Lemma \ref{lem:BornesAAFrag}.} For $0 \le f \in \dot L^1_{\alpha' + \gamma}$, we compute 
\bean
\| \AA^{c,3} f \|_{\dot L^1_{\alpha'}} 
&=& \int_0^\infty f(x) \, x^\gamma \, \chi_\delta^c(x) \int_0^x \wp(y/x) \, y^{\alpha'}  \chi_\eps(y) \,  {dy \over x}  \, dx 
\\
&\le& \| \wp (z) \, z^{\alpha'} \|_{L^1(0,2\eps/\delta)} \| f \|_{\dot L^1_{\alpha' + \gamma}}, 
\eean
and that establishes the last claim. For the other claims we refer to Lemma~\ref{lem:BornesAA} and \cite[3.~Proof of the main theorem]{MR2821681}
where very similar estimates have been proved. 
\qed

\begin{lem}\label{lem:BorneBBFrag} For any $a > a^*$ there exist $R,\delta,\eps > 0$ such that the operator $\BB$ is $a$-hypodissipative
in $X$. 
\end{lem}

\noindent
{\sl Proof of Lemma  \ref{lem:BorneBBFrag}.}   We define the space
$$
\EE := L^1(\phi), \quad \phi(x) :=  x^\alpha + \eta \, x^\beta, \quad \eta > 0.
$$  
\smallskip
Consider  a real-valued function  $f \in D(\Lambda)$ and let us show that for any $a > a^*$ and for suitable $R \in (1,\infty)$, $\delta \in (0,1)$, $\eps \in (0,\delta/2)$ to be chosen later, 
\beqn\label{eq:BFragDissip}
  \int_0^\infty  \sign(f(x)) \, \BB  f(x)\, \phi (x) \,dx
  \leq
  a \|f\|_{\EE}. 
\eeqn
Since the case of complex-valued functions can be handled in a similar way and since the norm  $\| \cdot \|_\EE$  is clearly equivalent to   $\| \cdot \|_{X}$,  that will end the proof.

From the identity 
$$
 \int_0^\infty  \sign(f(x)) \, \BB_0  f(x)\, x^r \, dr = \int_0^\infty ((r-1) x^r - x^{r+\gamma}) \, |f| \, dx
 $$
 and the following  inequality which holds for $\eta > 0$ small enough
$$
\eta \, (\beta - 1 - a) \, x^\beta  \le {1+a - \alpha \over 2} x^\alpha + \eta \, { 1- \wp_\beta \over 2} \, x^{\beta + \gamma} \quad \forall \, x > 0, 
$$
 we readily deduce 
\begin{multline*}
 \int_0^\infty  \sign(f(x)) \, \BB_0  f(x)\, \phi(x) \, dx \\
 \le a \,  \int_0^\infty|f| \, \phi  \, dx + \int_0^\infty ({\alpha - a - 1 \over 2}  x^\alpha - \eta \, {\wp_\beta+1 \over 2} x^{\beta+\gamma}) \, |f| \, dx .
\end{multline*}
On the other hand, we know from Lemma~\ref{lem:BornesAAFrag}  that  
\bean
\| \AA^{c}  f \|_{\EE} 
&\le& \wp_\alpha \delta^\gamma \, \| f \|_\EE+ ({\wp_\alpha \over R^{\beta-\alpha}} + \eta \, \wp_\beta ) \int_R^\infty x^{\beta+\gamma} |f|
\\
&&+  \| \wp (z) \, z^{\alpha} \|_{L^1(0,2\eps/\delta)} \, \| f \|_\EE.
\eean
We then easily conclude to \eqref{eq:BFragDissip} putting together these two estimates and choosing $R$ large enough, $\delta$ small enough, 
and then $\eps/\delta$ small  enough.  \qed

\medskip\noindent
{\sl Proof of Proposition~\ref{prop:FragSpectralMapping}. } We just have to explain why the hypothesis of Theorem~\ref{theo:QuantWeyl} are 
satisfied. Hypothesis {\bf (H1)} is an immediate consequence of Lemma~\ref{lem:BornesAAFrag}-(1), Lemma~\ref{lem:BorneBBFrag} and
 Remark~\ref{rem:theoEquivSGX}-(b).  

 \smallskip
We define 
$$
Y := \{ f \in W^{1,1}(\R); \,\, \hbox{supp} \, f \subset [\eps,2R] \}
$$
endowed with the norm $\| \cdot \|_{W^{1,1}}$. We also define on $X$ the projection $\pi_N$ onto the $2N+1$ dimension subspace
$$
R \pi_N := \{ {\bf 1}_{[0,2R]} \, p, \,\, p \in \Pp_{2N}(\R) \},
$$
where $\Pp_{2N}$ stands for the set of polynomials of degree less that $2N$, by 
$$
(\pi_N f) (x) := \chi^c_{\eps}(x) \, \chi_{2R}(x) \, (p_{N,R} * f)(x), \quad p_{N,R} (x) := p(x/(4R)) / (4R), 
$$
 where $p_N$ stands for  the Bernstein polynomial $p_N (x) := \alpha_N \, (1-x^2)^N$ with $\alpha_N$ such that $\| p_N \|_{L^1(-1,1)} = 1$. 
By very classical approximation arguments, we have 
$$
\| f - \pi_N f \|_{X} \le {C_{ R}  \over \sqrt{N}} \, \| f \|_{W^{1,1}} \quad \forall \, f \in Y,
$$
so that  \eqref{eq:YpiN} is fulfilled.

\smallskip
Finally,  Lemma~\ref{lem:BornesAAFrag}-(1) and Lemma~\ref{lem:BorneBBFrag} imply that 
$$
\| \AA S_\BB(t) f \|_{Y} \le C_a \, e^{at} \, \| f \|_X
$$
for any $t \ge 0$, $f \in X$, $a > a^*$. We then deduce  {\bf (H2)} with $n=1$, $\zeta=1$ and $\zeta'=0$  thanks to Remark~\ref{rem:theoEquivSGX}-(h) and $Y \subset D(\Lambda)$, 
and also deduce {\bf (H3)} with $n=1$ and \eqref{eq:hypQW} from the fact that  $Y$ is  compactly embedded in $X$ .  \qed
 
 %--------------------------------------------------------------------------------------------------------------------------------------------------------------

\subsection{A remark on the age structured population equation. } The aim of this short section is to present a quantified version of the Weyl's Theorem 
for the age structured population equation \eqref{eq:renewal} in the simple case when $\tau=\nu = 1$ and $K \in C_b(\R_+) \cap L^1(\R_+)$. 
More precisely, we consider the evolution equation 
\beqn\label{eq:renewallBIS}
\partial_t f = \Lambda f = \AA f + \BB f
\eeqn
with $\AA$ and $\BB$ defined on $M^1(\R)$ by 
\bean
(\AA f ) (x) &:=& \delta_{x=0} \int_0^\infty K(y) \, f(y) \, dy
\\
(\BB f ) (x) 
%&:=& - \partial_x(\tau(x) f(x)) - \nu(x) \, f(x)
%\\
&:=& - \partial_x f(x)  -  f(x),
\eean
and we emphasize that the boundary condition in  \eqref{eq:renewal}  has been equivalently replaced by the term $\AA f$ involving a Dirac mass $(\delta_0)(x) = \delta_{x=0}$ in $x=0$. 
%We indifferently denote $(\delta_0)(x) = \delta_{x=0}$. 
We  perform the spectral analysis of $\Lambda$  in the space $L^1(\R_+)$ as well as in the space $X := M^1(\R)$ of bounded measures endowed with the total variation norm. In that last functional space, the domain $D(\Lambda)$ is  the space $BV(\R)$ of functions with  bounded variation.   

\begin{lem}\label{lem:Renewal1} In $X = M^1(\R)$, the operators $\AA$ and $\BB$ satisfy

(i) $\AA \in \BBB(X,Y)$ where $Y = \C \, \delta_0 \subset X$ with compact embedding;   

(ii) $S_\BB(t)$ is $-1$-dissipative; 

(iii) the family of operators $S_\BB * \AA S_\BB(t)$ satisfy 
$$
\| (S_\BB * \AA S_\BB) (t)   \|_{X \to D(\Lambda) } \le C \, e^{-t} \qquad \forall \, t \ge 0. 
$$

\end{lem}

\noindent{\sl Proof of Lemma~\ref{lem:Renewal1}. } 
We clearly have $\AA \in \BBB(X,Y)$ because $K \in C_b(\R_+)$ and the $-1$-dissipativity of $S_\BB$ follows
from the explicit formula 
$$
S_\BB(t) f (x) = f(x-t) \, e^{-t} \, {\bf 1}_{x-t \ge 0}. 
$$
We next prove (iii). We write
$$
\AA S_\BB (t) f = \delta_{x=0} \int_0^\infty K(y+t) \, f(y) \, dy \, e^{-t},
$$
next
\bean
(S_\BB * \AA S_\BB) (t) f 
&=&
\int_0^t S_\BB(s) \, \AA S_\BB(t-s) f \, ds 
\\
&=& \int_0^t  \delta_{x-s = 0} e^{-s} \, {\bf 1}_{x-s \ge 0}   \int_{0} ^\infty K(y+t-s) \, f(y) \, dy \, e^{-(t-s)} \, ds
\\
&=& e^{-t} \, {\bf 1}_{x \le t}    \int_{0} ^\infty K(y+t-x) \, f(y) \, dy,
\eean
and then 
\bean
\partial_x (S_\BB * \AA S_\BB) (t) f 
&=& - e^{-t} \, {\delta}_{x = t}    \int_{0} ^\infty K(y+t-x) \, f(y) \, dy
\\
&&- e^{-t} \, {\bf 1}_{x \le t}    \int_{0} ^\infty K'(y+t-x) \, f(y) \, dy.
\eean
As a consequence, we deduce 
$$
\| \partial_x (S_\BB * \AA S_\BB) (t) f  \|_{M^1} \le \| K \|_{W^{1,\infty}} \, e^{-t} \, \| f \|_{M^1} 
$$
and a similar estimate for $\| (S_\BB * \AA S_\BB) (t) f  \|_{M^1}$. \qed 

\medskip

As a first step in the proof of Theorem~\ref{theo:Frag} for the age structured population equation, we have 

\begin{prop}\label{prop:RenSpectralMapping} Under the above assumptions and notation, the conclusion {\bf (2)},  and then \eqref{eq:localizSpectrL}  and {\bf (1)},  of Theorem~\ref{theo:QuantWeyl} 
holds for the age structured population  semigroup $S_\Lambda$ in $M^1(\R)$ and in $L^1(\R_+)$ for any $a > a^* := - 1$. Moreover,  any eigenvalue $\xi \in \Sigma_d(\Lambda) \cap \Delta_{-1}$ is algebraically simple. 
\end{prop}

\noindent{\sl Proof of Proposition~\ref{prop:RenSpectralMapping}. } 
In order to prove the result in $X = M^1(\R_+)$, we just have to explain  again why the hypothesis of Theorem~\ref{theo:QuantWeyl} are 
satisfied. Conditions {\bf (H1)},  {\bf (H3)} and \eqref{eq:hypQW} are  immediate consequences of Lemma~\ref{lem:Renewal1}-(i) \& (ii) together with Remark~\ref{rem:theoEquivSGX}-(b). 
We refer to \cite{MR0354068,MR684412,MR860959} and the references therein for the existence  theory in $L^1(\R_+)$ (which extends without difficulty to $M^1(\R_+)$) for the semigroup $S_\Lambda$. 
Hypothesis {\bf (H2)} with $n=1$, $\zeta=1$ and $\zeta'=0$  is nothing but Lemma~\ref{lem:Renewal1}-(iii). 
 Finally, taking up again the proof of Theorem~\ref{theo:QuantWeyl} and using the additional fact that dim$Y=1$, we get the algebraic simplicity of the eigenvalues $\xi \in \Sigma(\Lambda) \cap \Delta_{-1}$.  
 
  \smallskip
We may then easily extend the  spectral analysis performed in $M^1(\R_+)$ to the functional space $ L^1(\R_+)$. Indeed, 
for a function $f \in L^1$, Theorem~\ref{theo:QuantWeyl} implies that 
$$
\| e^{t\Lambda} f - \sum_{j=1}^J e^{t\xi_j} \Pi_j f \|_{M^1} \le C_a \, e^{at} \, \| f \|_{M^1}
$$
for any $a > - 1$ and $t \ge 0$. 
Because  $S_\Lambda$ is well defined in $L^1$ and  the domain of $\Lambda$ as an operator in $M^1(\R_+)$ is $BV(\R_+) \subset L^1(\R_+)$,
all the terms involved in the above expression belong to $L^1$ and we can replace the norms $\| . \|_{M^1}$ by the norms $\| . \|_{L^1}$.   \qed

%%%%%%%%%%%%%%%%%%%% Spectral mapping and Weyl's  %%%%%%%%%%%%%%%%%%%%%%%%%%%%%%%%%%%

\bigskip
\section{The Krein-Rutman  Theorem in an abstract setting}
%\section{positive semigroup and Krein-Rutman Theorem}
\label{sec:KRthAbstract}
\setcounter{equation}{0}
\setcounter{theo}{0}

%\subsection{The Krein-Rutman  Theorem in an abstract setting}

In this section we consider a ``Banach lattice of functions" $X$. 
We recall that a Banach lattice is a Banach space endowed with an order denoted by $\ge $ (or $\le$) 
such that the following holds: 

- The set $ X_+ := \{ f \in X; \,\, f \ge 0 \}$ is a nonempty convex closed cone. 

-  For any  $f \in X$, there exist some unique (minimal) $f_\pm \in X_+$ such that $f = f_+ - f_-$, we then denote $|f| := f_+ + f_- \in X_+$. 

- For any $f,g \in X$, $0 \le f \le g$ implies  $\| f \| \le \| g \|$.  

\smallskip
We may define a dual order $\ge$ (or $\le$) on $X'$ by writing for $\psi \in X'$
$$
\psi \ge 0 \,\,\,  (\hbox{or } \psi \in X'_+) \quad\hbox{iff}\quad
\forall \, f \in X_+ \,\,\, \langle \psi,f \rangle \ge 0,
$$
so that $X'$ is also a Banach lattice.

\smallskip
We then restrict our analysis to the case when $X$ is a``space of functions". 
The examples of spaces we have in mind are the space of Lebesgue functions $X = L^p(\UUU)$, $1 \le p < \infty$,  $\UUU \subset \R^d$ borelian set, 
 the space $X = C(\UUU)$ of  continuous functions on a compact set $\UUU$ and the space $X = C_0(\UUU)$ of uniformly continuous functions defined on an open set $\UUU \subset \R^d$ and which tend to $0$ at the boundary of $\UUU$. For any element (function) $f$ in such a ``space of functions"  $X$, we may define without difficulty the composition functions $\theta(f)$ and $\theta'(f)$ for $\theta(s) = |s|$ and $\theta(s) = s_\pm$ as well as the support $\hbox{supp} f$ as a closed subset of $\UUU$.  
Although we believe that our results extend to a broader class of Banach lattices, by now on and in order to avoid technicality,  we will restrict ourself to these examples of  ``space of functions" without specifying anymore but just saying that we consider a ``Banach lattice of functions" (and we refer to the textbook \cite{AGLMNS} for possible generalisation). As a first consequence of that choice, we may then  obtain a nice and simple property on the generator of a positive semigroup. 
 
\begin{defin}\label{def:positiveSG} Let us consider a Banach lattice $X$ and a generator $\Lambda$ of a semigroup $S_\Lambda$ on $X$. 

(a) - We say that the semigroup $S_\Lambda$ is positive if $S_\Lambda(t) f \in X_+$ for any $f \in X_+$ and $t \ge 0$.

(b) - We say that a generator $\Lambda$ on $X$ satisfies Kato's inequalities if  the inequality 
\beqn\label{eq:KatoIneq}
\forall \, f \in D(\Lambda) \quad \Lambda \theta(f) \ge \theta'(f) \, \Lambda f
\eeqn
holds for $\theta(s) = |s|$ and $\theta(s) = s_+$. 

(c) - We say that  $-\Lambda$ satisfies a "weak maximum principle" if  for any $a > \omega(\Lambda)$ and $g \in X_+$ there holds
\beqn\label{eq:PM*}
f \in D(\Lambda) \hbox{ and } (- \Lambda + a) f = g \quad\hbox{imply}\quad f \ge 0. 
\eeqn

(d) - We say that the opposite of the resolvent is a positive operator if for any $a > \omega(\Lambda)$ and $g \in X_+$ there holds $- R_\Lambda(a) g \in X_+$. 

\end{defin} 

Here the correct way to understand Kato's inequalities is 
$$
\forall \, f \in D(\Lambda), \,\,\, \forall \, \psi \in D(\Lambda^*) \cap X'_+  \quad  \langle \theta(f), \Lambda^* \psi \rangle  \ge \langle \theta'(f) \, \Lambda f,\psi \rangle,
$$
where $\Lambda^*$ is the adjoint of $\Lambda$. 

\smallskip
It is well known (see \cite{NagelUhlig81} and \cite[Remark 3.10]{Arendt82} and the textbook \cite[Theorems  C.II.2.4, C.II.2.6 and Remark C-II.3.12]{AGLMNS} ) that the generator $\Lambda$ of a positive semigroup $S_\Lambda$ on one of our {\it Banach lattice space of functions} satisfies Kato's inequalities \eqref{eq:KatoIneq}. 
It is also immediate from the Hille's identity \eqref{eq:RL=intSt} that (a) implies (d) and then (c) in the general Banach lattice framework. 
For a broad class of spaces $X$  the properties  (a), (b), (c) and (d) are in fact equivalent and we refer again to the textbook \cite{AGLMNS} for more details on that topics.

\smallskip
Last, we need some strict positivity notion on $X$ and some strict positivity (or irreducibility) assumption on $S_\Lambda$ that we will formulate in term of {\it ``strong maximum principle"}. It is worth mentioning that we have not assumed that $X_+$ has nonempty interior, so that the strict positivity property cannot be defined 
using that interior set (as in Krein-Rutman's work \cite{MR0027128}, see also \cite{MR1064315}).
However, we may define the strict order $>$ (or $<$) on $X$ by writing for $f \in X$ 
$$
f >  0 \quad\hbox{iff}\quad
\forall \, \psi \in X'_+ \backslash\{0\}\,\,\, \langle \psi,f \rangle > 0 , 
$$
and similarly a strict order $>$ (or $<$) on $X'$ by writing for $\psi \in X'$ 
$$
\psi >  0 \quad\hbox{iff}\quad
\forall \, g \in X_+ \backslash\{0\}\,\,\, \langle \psi,g \rangle > 0 .
$$
It is worth emphasising that from the Hahn-Banach Theorem, for any $f \in X_+$ there exists $\psi \in X'_+$ such that $\| \psi \|_{X'} = 1$ and 
$\langle \psi,f \rangle = \| f \|_X$ from which we easily deduce that 
\beqn\label{eq:StrictOrderNorm}
\forall\, f,g \in X, \quad 0 \le f < g \quad \hbox{implies}\quad \| f \|_X < \| g \|_X.
\eeqn

 \begin{defin}\label{def:PMstrong} We say that $-\Lambda$ satisfies the ``strong maximum principle" if for any given  $f \in X$ and $ \mu \in \R$, 
 there holds
% $$
% f_+ \in D(\Lambda)  \backslash \{0\} \, \hbox{ and }\,  \Lambda f_+ = \mu f_+
%\quad\hbox{imply} \quad f > 0 
%$$
%and also 
$$
|f| \in D(\Lambda)  \backslash \{0\} \, \hbox{ and }\,  ( -  \Lambda + \mu ) |f| \ge 0
\quad\hbox{imply} \quad f > 0 \, \hbox{ or }\,  f < 0. 
$$
\end{defin}

We can now state the following version of the Krein-Rutman  Theorem in a general and abstract setting. 

%There exists $n_0 \in \N$ such that for $f \in X$, we have
%$$
%|f| \in D(\Lambda^{n_0}) \,\,\hbox{and}\,\,\, |f| > 0 \quad\hbox{imply} \quad |f| = f_+   \,\,\hbox{or}\,\,\,  |f| = f_-
%$$
%and
%$$
%f_+ \in D(\Lambda^{n_0}) \,\,\hbox{and}\,\,\, f_+ > 0 \quad\hbox{imply} \quad  f = f_+.
%$$
%and
%$$
%f \in D(\Lambda) \cap X_+ \backslash \{0\}, \,\, \mu \in \R, \,\, \Lambda f = \mu f
%\quad\hbox{imply} \quad f > 0; 
%$$
%
%defined in  the following strict positivity property holds
%$$
%f \in D(\Lambda) \cap X_+ \backslash \{0\}, \,\, \mu \in \R, \,\, \Lambda f = \mu f
%\quad\hbox{imply} \quad f > 0; 
%$$
%

\begin{theo}\label{theo:KR} We consider a  generator $\Lambda$ of a semigroup $S_\Lambda$ on a  Banach lattice of functions $X$, and we assume that 

{\bf (1)} $\Lambda$ satisfies the property {\bf (1)} of the semigroup Weyl's Theorem~\ref{theo:GalWeyl} for some $a^* \in \R$;

{\bf (2)}  there exist $b > a^*$ and $\psi \in D(\Lambda^*) \cap X'_+ \backslash \{0\}$ such that $\Lambda^* \psi \ge b \, \psi$;

{\bf (3)}  $S_\Lambda$ is positive (and $\Lambda$ satisfies Kato's inequalities);   
 
 {\bf (4)} $-\Lambda$ satisfies a strong maximum principle.
 % as presented in Definition~\ref{def:PMstrong}. 
 
\smallskip\noindent
Defining $\lambda := s(\Lambda)$, there holds
$$
a^* < \lambda = \omega(\Lambda)  \quad \hbox{and} \quad \lambda \in \Sigma_d(\Lambda), 
$$
and there exists $0 < f_\infty \in D(\Lambda)$ and $0 < \phi \in D(\Lambda^*)$ such that 
$$
\Lambda f_\infty = \lambda \, f_\infty, \quad \Lambda^* \phi = \lambda \, \phi, \quad 
R\Pi_{\Lambda,\lambda} = \hbox{Vect} (f_\infty),
$$
and then 
$$
\Pi_{\Lambda,\lambda} f = \langle f,\phi\rangle \, f_\infty \quad \forall \, f \in X.
$$
Moreover, there exists  $a^{**} \in (a^*,\lambda)$ and for any $a > a^{**}$ there exists $C_a > 0$ such that   for any $f_0 \in X$ 
$$
\|  S_\Lambda(t) f - e^{\lambda t} \,\Pi_{\Lambda,\lambda}  f_0 \|_X \le C \, e^{a t} \, \| f_0 - \,\Pi_{\Lambda,\lambda}  f_0 \|_X 
\qquad
\forall \, t \ge 0.
$$
 
 \end{theo}

\begin{rem}\label{rem:KR}
\begin{itemize}
\item[(a)] Theorem~\ref{theo:KR} generalises the Perron-Frobenius Theorem \cite{MR1511438,Frobenius} for strictly positive matrix  $\Lambda$, 
the Krein-Rutman  Theorem \cite{MR0027128} for irreducible, positive and compact semigroup on a Banach lattice with non empty interior cone and
the Krein-Rutman Theorem variant \cite[Corollary III-C.3.17]{AGLMNS} (see also \cite{MR618205}
 for the original proof)  for irreducible, positive, eventually norm continuous semigroup with compact resolvent in a general Banach lattice framework. 
We also refer to the book by Dautray and Lions \cite{MR1064315} for a clear and comprehensible version of the the Krein-Rutman  Theorem as well as
to the recent books \cite{MR1612403,MR2178970} and the references therein for more recent developments on the theory of positive operators.  
The main novelty here is that with assumptions {\bf (1)} and {\bf (2)} we do not ask for the semigroup to be eventually norm continuous 
and we only ask power compactness on the decomposition $\AA$ and $\BB$ of the operator $\Lambda$ instead of compactness on the resolvent $\RR_\Lambda$.  

\item[(b)] Condition {\bf (2) } is necessary because  the requirement {\bf (1) } only implies the needed compactness and regularity on the iterated  operator $(\AA \, \RR_\BB(z))^n$ for $z \in \Delta_{a^*}$. Condition {\bf (2)} can be removed (it is automatically verified) if condition {\bf (1)} holds for any $a^* \in \R$. 

\item[(c)] In a general Banach lattice framework and replacing the strong maximum principle hypothesis {\bf (4)} by the more classical irreducibility assumption on the semigroup (see for instance \cite[Definition C-III.3.1]{AGLMNS})  Theorem~\ref{theo:KR} is an immediate consequence of   Theorem~\ref{theo:GalWeyl} together with  \cite[Theorem C-III.3.12]{AGLMNS} (see also \cite{MR591851,MR618205,MR617977}). We do not know whether the strong maximum principle and the irreducibility are equivalent assumptions although both are related to strict positivity of the semigroup or the generator. Anyway the strong maximum principle for the operator $-\Lambda$ is a very natural notion and hypothesis from our PDE point of view and that is the reason why we have chosen to presented the statement of Theorem~\ref{theo:KR} in that way. Moreover, the proofs in \cite[part C]{AGLMNS} are presented in the general framework of Banach lattices and positive or reducible semigroups (no compactness assumption is required in the statement of \cite[Theorem C-III.3.12]{AGLMNS}) so that quite abstract 
%(and maybe obsur for no initiated people) 
arguments are used during the proof. We give below a short, elementary and somewhat self-contained proof of Theorem~\ref{theo:KR}, and thus do not use   \cite[Theorem C-III.3.12]{AGLMNS}. 

\end{itemize}
\end{rem}

\medskip

\noindent
{\sl Proof of Theorem~\ref{theo:KR}. } We split the proof into five steps. 

\smallskip\noindent
{\sl Step 1. } 
%We define 
%$$
%\lambda := \sup \{ \Re e \, \xi, \, \xi \in \Sigma(\Lambda) \}.
%$$
On the one hand, let us fix $0 \le f_0 \in D(\Lambda)$ such that $C :=  \langle f_0, \psi \rangle  > 0$ which exists by definition of $\psi$. Then denoting $f(t) := S_\Lambda(t) f_0$, we have 
$$
{d \over dt} \langle f(t) , \psi \rangle =  \langle \Lambda f(t) , \psi \rangle = \langle  f(t) , \Lambda^* \psi \rangle \ge 
b \, \langle f(t) , \psi \rangle, 
$$
which in turn implies
\beqn\label{eq:f(t)>ebt}
\langle f(t) , \psi \rangle \ge C \, e^{bt} \quad \forall \, t \ge 0. 
\eeqn
On the other hand, from   Theorem~\ref{theo:SpectralMapping} we know that $\omega(\Lambda) \le \max(a^*,\lambda)$. As a consequence,  if $\lambda  < b$ for any $a \in (\max(a^*,\lambda),b)$   there exists $C_a \in (0,\infty)$ such that  
$$
 \langle f(t) , \psi \rangle \le  \| \psi \|_{X'} \,  \| f(t) \|_X \le C_a \, \| \psi \|_{X'}  \, e^{a t} \, \| f_0 \|_X.
$$
That would be in contradiction with \eqref{eq:f(t)>ebt}. We conclude that $a^* < b \le \lambda = \omega(\Lambda)$. 
  
%$$
%s(\Lambda) \in \Sigma_d(\Lambda) \quad\hbox{and}\quad N(\Lambda-s(\Lambda)) \cap K^* \not = \emptyset.
%$$
\smallskip\noindent
{\sl Step 2. } We prove that there exists $f_\infty \in X$ such that 
\beqn\label{eq:existFinfty}
\| f_\infty \|=1, \,\,\, f_\infty > 0, \,\,\,  \Lambda f_\infty = \lambda f_\infty. 
\eeqn
Thanks to the Weyl's Theorem~\ref{theo:GalWeyl} we know that for some $a < 0$ 
$$
\Sigma (\Lambda) \cap \Delta_a = \cup_{j=1}^J \{ \xi_j \} \subset \Sigma_d(\Lambda), \quad \Re e \xi_j = \lambda,
$$
with $J \ge 1$.  We introduce the Jordan basis $\VVV := \{ g_{1,1}, ..., g_{J,L_J} \}$ of $\Lambda$ in the invariant subspace $R\Pi_{\Lambda,a}$ 
as the family of vectors
$$
g_{j,\ell} \not = 0, \quad \Lambda g_{j,\ell} = \xi_j g_{j,\ell} + g_{j,\ell-1} \quad \forall \, j \in \{1, ..., J \}, \,\, \forall \, \ell \in \{ 1, ..., L_j \}, 
$$
with the convention $g_{j,k} = 0$ if $k \le 0$ or $k \ge L_j+1$,  as well as the projectors (associated to the basis $\VVV$) 
$$
\Pi_{j,\ell} := \hbox{projection on} \, g_{j,\ell}, \quad \Pi_k := \hbox{projection on} \, \hbox{Vect}(g_{j,k}, \, 1 \le j \le J \}.
$$
 For any fix $g \in \VVV$ we write $g = g^1 - g^2 + i g^3 - i g^4$ with $g^\alpha \ge 0$ and we remark that there exists $\alpha \in \{1, ..., 4 \}$ and $k_0 \ge 1$ such that $\Pi_{k_0} g^\alpha = \Pi_{\Lambda,a} g^\alpha \not= 0$. We then define $k^* = k^*(g^\alpha) :=\max \{ k; \,\, \Pi_k g^\alpha \not= 0 \}$, that set being not empty
 (since it contains $k_0$).  
We may split the semigroup  as
\bean
e^{\Lambda t} g^\alpha &=& \sum_{j}\sum_\ell e^{\Lambda t} \, \Pi_{j,\ell} g^\alpha + e^{\Lambda t} \, (I-\Pi_{\Lambda,a}),
\eean
with
$\Pi_{j,\ell} g^\alpha = (\pi_{j,\ell} g^\alpha) \, g_{j,\ell} $, $\pi_{j,\ell} g^\alpha \in \C$, and 
$$
e^{\Lambda t} \,  g_{j,\ell} = e^{\xi_j t} g_{j,\ell} + ...  + t^{\ell-1} \, e^{\xi_j t} g_{j,1}. 
$$
Using the positivity assumption {\bf (3)}  and keeping only the leading order term in the above expressions, we have thanks to Theorem~\ref{theo:GalWeyl} 
$$
0 \le {1 \over t^{k^*-1}} e^{(\Lambda -\lambda)t} g^\alpha = \sum_{j=1}^J (\pi_{j,k^*} g^\alpha) e^{(\xi_j -\lambda)t} g_{j,1} + o(1). 
$$
There exist a sequence $(t_n)$ which tends to infinity and complex numbers $z_j \in \C$, $|z_j| = 1$, such that,  passing to the limit in the above expression, we get 
$$
0 \le \sum_{j=1}^J (\pi_{j,k^*} g^\alpha) \, z_j \,  g_{j,1} =: g_\infty. 
$$
Because of the choice of $\alpha$ and the fact that the vectors $g_{j,1}$, $1 \le j \le J$,  are independent, we have 
then $ g_\infty \in R\Pi_{\Lambda,a} \cap X_+ \backslash \{ 0 \}$. 
%$$
% g_\infty = \sum_{j=1}^J \theta_j \,  g_{j,1} \in R \Pi_1 \cap K \backslash \{ 0 \},
%$$
%for some (not all vanishing) complex numbers $\theta_j \in \C$. 
Applying again the semigroup, we get
$$
0 \le e^{\Lambda t}  g_\infty = \sum_{j=1}^J e^{\xi_j t}  [(\pi_{j,k^*} g^\alpha) \, z_j \,  g_{j,1} ] \quad \forall \, t \ge 0,
$$
which in particular implies 
$$
 \sum_{j=1}^J \Im m \bigl\{ e^{\xi_j t}  [(\pi_{j,k^*} g^\alpha) \, z_j \,  g_{j,1} ] \bigr\}  = 0  \quad \forall \, t \ge 0, 
$$
and then $\pi_{j,k^*} g^\alpha  = 0$ if $\Im m \, \xi_j \not = 0$. As a conclusion, we have proved that there exists $ g_\infty \in N(\Lambda-\lambda) \cap X_+ \backslash \{ 0 \}$. Together with the strong maximum principle we conclude that $f_\infty := g_\infty  / \| g_\infty \|$ satisfies \eqref{eq:existFinfty}. 

Moreover, the above argument for any $g = g_{j,L_j}$ associated to $\xi_j \not= \lambda$ and for any $\alpha \in \{1,..., 4 \}$ such that $\Pi_{j,L_j} g^\alpha \not = 0$
implies that $\pi_{j,k^*} g^\alpha  = 0$,  or in other words 
\beqn\label{eq:dimespacespropres}
\max_{\xi_j \not = \lambda} L_j < \max_{\xi_j = \lambda} L_j.
\eeqn
 
\smallskip\noindent
{\sl Step 3. } We prove that there exists $\phi \in X'$ such that 
\beqn\label{eq:existphi}
\phi > 0, \,\, \Lambda^* \phi = \lambda \phi. 
\eeqn
We define $S^*_\Lambda$ to be the dual semigroup associated to $S_\Lambda$ and we emphasize that it is not necessarily strongly  continuous (for the norm in $X'$) but only weakly continuous (for the weak topology $\sigma(X',X)$).  However, introducing the splitting 
$$
S^*_\Lambda = (S_\Lambda (I - \Pi_{\Lambda,a}))^* + (S_\Lambda  \Pi_{\Lambda,a})^*
$$
and observing that 
$$
\| (S_\Lambda (I - \Pi_{\Lambda,a}))^*\|_{\BBB(X)} \le C_a \, e^{at}
$$
for some $a < \lambda$, the same finite dimension argument as in Step 2 implies that there exists $\phi \in N(\Lambda^* - \lambda) \cap X' \backslash \{ 0 \}$. 

Let us prove the strict positivity property. For  $a > s(\Lambda)$  and $g \in X_+ \backslash \{ 0 \}$, thanks to the weak and strong maximum principles {\bf (3)} and {\bf (4)},   there exists $0 < f \in X$ such that 
$$
(-\Lambda + a) f = g.
$$
As a consequence, we have 
\bean 
\langle \phi, g \rangle &=& \langle \phi, (-\Lambda + a) f   \rangle 
\\
&=&  \langle  (a-\Lambda^*) \phi, f   \rangle = (a - \lambda) \, \langle  \phi, f   \rangle  > 0. 
\eean
Since $g \in X_+$ is arbitrary, we deduce that $\phi > 0$. That concludes the proof of  \eqref{eq:existphi}.

\smallskip\noindent
{\sl Step 4. } We prove that $N(\Lambda-\lambda) = \hbox{Vect}(f_\infty)$. Consider a normalized eigenfunction $f \in X^\R \backslash\{0\}$ associated to the eigenvalue  $\lambda$. 
First we observe that from Kato's inequality 
\bean
\lambda |f| =   \lambda f  \, \hbox{sign} (f)  
 = \Lambda f  \, \hbox{sign} (f) 
\le \Lambda |f|. 
\eean
That inequality is in fact an equality, otherwise we would have 
\bean
\lambda \langle |f|,\phi \rangle \not=   \langle \Lambda |f| ,\phi \rangle  = \langle  |f| , \Lambda^* \phi \rangle  = \lambda \langle |f|,\phi \rangle,
\eean
and a contradiction. As a consequence, $|f|$ is a solution to the eigenvalue problem $\lambda |f|  = \Lambda |f|$ so that the strong maximum principle assumption {\bf (4)} implies $f > 0$ or $f< 0$, and without loss of generality we may assume $f > 0$. Now, thanks to Kato's inequality again, we write 
\bean
\lambda  (f-f_\infty)_+ =  \Lambda (f-f_\infty) \, \hbox{sign}_+(f-f_\infty)  \le  \Lambda (f-f_\infty)_+ , 
\eean
and for the same reason  as above that last inequality is in fact an inequality. Since $(f-f_\infty)_+ = |(f-f_\infty)_+|$, the strong maximum principle implies that either $(f-f_\infty)_+ = 0$, or in
other words $f \le f_\infty$,  either $(f-f_\infty)_+ > 0$ or in other words $f > f_\infty$. Thanks to \eqref{eq:StrictOrderNorm} 
and  to the normalization hypothesis $\| f \| = \| f_\infty \| = 1$ the second case in the above alternative is not possible. 
Repeating the same argument with $(f_\infty-f)_+$ we get that $f_\infty \le f$ and we conclude with  $f = f_\infty$. For a general eigenfunction $f \in X^\C$ associated to the eigenvalue  $\lambda$ we may introduce the decomposition $f = f_r + i f_i$ and we immediately get that $f_\alpha \in X^\R$ is an eigenfunction associated to  $\lambda$ for $\alpha = r,i$. As a consequence of what we have just established, we have $f_\alpha = \theta_\alpha f_\infty$ for some $\theta_\alpha \in \R$ and we conclude that  $f = (\theta_r + i \theta_i) \, f_\infty  \in \hbox{Vect}(f_\infty)$ again. 

\smallskip\noindent
{\sl Step 5. } We first claim that $\lambda$ is algebraically simple. Indeed, if it was not the case, there would exist $f \in X^\R$ such that  $\Lambda f = \lambda f + f_\infty$ 
and then  
$$
 \lambda \langle f,  \phi \rangle = \langle f, \Lambda^* \phi \rangle = \langle \Lambda f, \phi \rangle = \langle \lambda f + f_\infty, \phi \rangle, 
 $$
 which in turn implies $\langle f_\infty, \phi \rangle = 0$ and a contradiction. With the notation of step 2 and thanks to \eqref{eq:dimespacespropres}, that implies that for any $\xi_j \not= \lambda$ there holds $L_j < 1 $  or in other words $\Sigma(\Lambda) \cap \bar\Delta_\lambda = \{ \lambda \}$.   We conclude the proof by using the semigroup Weyl's Theorem~\ref{theo:GalWeyl}, which in particular implies that $\Sigma(\Lambda) \cap \Delta_{a^{**}} = \{ \lambda \}$ for some $a^{**}  \in (a^*, \lambda)$ because  $\Sigma(\Lambda) \cap \Delta_{a}$ is finite for any $a > a^*$. 
\qed

%---------------------------------------------------------------------------------------------------------------------------------------------------------------------

\section{Krein-Rutman Theorem for the growth-fragmentation equations}\label{subsec:KRGrFr}
%---------------------------------------------------------------------------------------------------------------------------------------------------------------------

This section is devoted to the proof of  Theorem~\ref{theo:Frag}.  
\subsection{General growth-fragmentation equations}\label{subsec:KRGrFr}
We present below the proof of the part of   Theorem~\ref{theo:Frag} which holds in full generality. Namely, we prove \eqref{eq:PrimalEVpb}, \eqref{eq:RateGrFr}
and \eqref{eq:DualEVpb} for the three models of growth-fragmentation equations. 

\smallskip
We start with the cell-division equation for which we apply the Krein-Rutman Theorem~\ref{theo:KR} in the Banach lattice $L^1((x_0z_0,\infty); \langle x \rangle^\alpha dx)$, $\alpha > 1$, where $z_0$ is defined by \eqref{def:z0}, instead of $L^1_\alpha$ in order that the operator $\Lambda$ enjoys a strong maximum principle. 
We have proved {\bf (1)}   
in Proposition~\ref{prop:GrFrSpectralMapping}. We have $\Lambda^* 1 = K(x) \, (1- \wp_0) \ge 0$ so that {\bf (2)} holds
with $b = 0 > a^*$. The weak maximum principle {\bf (3)} is an immediate consequence of Kato's inequalities which in
turn follows from the fact that $\FF^+ \theta (f) \ge \theta'(f) \FF^+ f$ $a.e$  for any $f \in D(\Lambda)$, $\theta(s) = |s|$ and $\theta(s) = s_+$. 
The strong maximum principle {\bf (4)} follows from the fact that the equation
$$
|f| \in D(\Lambda)  \backslash \{0\} \, \hbox{ and }\,  ( -  \Lambda + \mu ) |f| \ge 0
$$
may be rewritten as
$$
- \partial_x |f| + (K(x) + \mu) |f|  \ge \FF^+ |f| \ge 0,
$$
and we conclude as in the \cite[proof of Theorem 3.1]{MMP} that the continuous function $|f|$ does not vanish on $(x_0z_0,\infty)$, so that 
$f > 0$ on $(x_0z_0,\infty)$ or $f < 0$ on $(x_0z_0,\infty)$. 

\smallskip
For the self-similar equation the proof is exactly the same by applying the Krein-Rutman Theorem~\ref{theo:KR} in the Banach lattice $L^1(\R_+; (x^\alpha + x^\beta) dx)$, $0 \le \alpha < 1 < \beta < \infty$. Let us just emphasize that condition {\bf (1)} has been proved in Proposition~\ref{prop:FragSpectralMapping} and that condition {\bf (2)} 
with $b=0 > a^*$ follows from the fact that
$\Lambda^* \phi = 0$ for the positive function  $\phi(x) = x$. We refer to \cite[Section 3]{EMRR} for the proof of the weak and strong maximum principles
 {\bf (3)} and {\bf (4)}. 
  
\smallskip
Finally, for the age structured population equation we apply the Krein-Rutman  Theorem~\ref{theo:KR} in the Banach lattice $L^1(\R_+)$.  As in \cite[Appendix]{MR1946722}, we observe that,  denoting by $\lambda > -1$ the real number such that 
$$
\int_0^\infty K(x) \, e^{-(1+\lambda) \, x} \, dx = 1, 
$$
which exists thanks to condition \eqref{eq:renewal1}, 
the function 
$$
\psi (x) := e^{(1+\lambda) x} \int_x^\infty K(y)  \, e^{-(1+\lambda) \, y} \, dx 
$$
is a solution to the dual eigenvalue problem 
$$
\Lambda^* \psi = \partial_x \psi - \psi + K(x) \psi(0) = \lambda \psi, \quad 0 \le \psi \in L^\infty(\R_+), 
 $$
 so that in particular {\bf(2)} holds with $\lambda > a^*= -1$.  Condition {\bf (1)} has been proved in Proposition~\ref{prop:RenSpectralMapping} and
 the proof of the weak and strong maximum principles {\bf (3)} and {\bf (4)} is classical.  

%---------------------------------------------------------------------------------------------------------------------------------------------------------------------

\subsection{Quantified spectral gap theorem for the cell division equation with constant total  fragmentation rate}
%{Growth-fragmentation equation with constant total rate}

We consider the particular case of the cell-division equation with constant total fragmentation rate and 
fragmentation kernel which furthermore fulfils condition \eqref{eq:intkappa<0} 
for which we can give an accurate long-time asymptotic behaviour (as formulated in point (i) of Theorem~\ref{theo:Frag}) and answer to a question formulated
in \cite{MR2114128,MR2536450}.  We then consider the equation
\beqn\label{eq:GrowthFragCst}
\partial_t f + \partial_x f +  K_0 f  = K_0 \int_x^\infty \kappa(y,x) \, f(y) \, dy
\eeqn
with vanishing boundary condition \eqref{eq:BoundaryCondition}, where $K_0 > 0$ is a constant and $\kappa$ satisfies  \eqref{eq:intkappa<0}. 
In such a situation, we have the following accurate description of the spectrum.  
 
\begin{prop}   \label{theo:K=Cste} The first eigenvalue is given by $\lambda = s(\Lambda):= (n_F-1) K_0$ with $n_F$ defined in  \eqref{eq:intkappa<0}.
On the other hand, for any $a^{**} \in (-  K_0, (n_F-1) K_0)$ and any $\alpha > \alpha^*$, with $\alpha^*$ large enough (but explicit and given during the proof), 
the spectral gap $\Sigma(\Lambda) \cap \Delta_{a^{**}}  = \{ \lambda \}$ holds in $X := L^1_\alpha$. 

%Theorem~\ref{theo:Frag} holds in $L^1_\alpha$ for any $a > a^{**}$, which in particular gives an  explicit bound
%on the spectral gap $\lambda - \sup \Re e \Sigma(\Lambda) \backslash \{ \lambda \}$.  
% .  For any $a > 0$,   there exists a constant $C_a$ such that 
% for any $f_0 \in L^1_\alpha$ the following estimate holds
% \bean
% \| e^{t\Lambda} f - \langle f \rangle f_\infty \, e^{\lambda t}  \|_{L^1_\alpha} \le C_a \, e^{at} \, 
% \|  f - \langle f \rangle f_\infty \|_{L^1_\alpha}
% \eean
%  \begin{itemize}
%\item blabla
%    \end{itemize}
\end{prop}

We use the following extension (shrinkage) of the functional space of the semigroup
  decay proved in \cite{MM*}. 
  
\begin{theo}[Extension of the functional space of the semigroup
  decay]   \label{theo:Extension}
  
  Let 
  
  \noindent
  $E$ and $\EE$ be two Banach spaces such that  $E \subset \EE$  with dense and 
  continuous embedding, and let $L$ be the generator of a semigroup
  $S_L(t) := e^{tL}$ on $E$, $\LL$ the generator of a semigroup $S_\LL(t) :=  e^{t\LL}$ on $\EE$ 
  with $\LL_{|E} = L$.

\smallskip\noindent
We assume that there exist two operators $\AA, \BB \in \mathscr{C}(\EE)$ such that 
$$
\LL = \AA + \BB, \,\,\,  L = A + B, \,\,\, A = \AA_{|E}, \,\,\, B = \BB_{|E},
$$
and a real number $a^{**} \in \R$ such that there holds: 

\begin{itemize}

\item[{\bf (i)}]   $B-a$ is  hypodissipative on $E$, $\BB-a$ is  hypodissipative on $\EE$ for any $a > a^{**}$;
  \item[{\bf (ii)}]  $A \in \BBB(E)$,  $\AA \in \BBB( \EE)$;
  \item[{\bf (iii)}]  there is $n \ge 1$ such that, for any $a > a^{**}$ and for some constant $C'_a \in (0,\infty)$, 
\[
\Bigl\|  (\AA S_\BB)^{(*n)}(t)  \Bigr\|_{\BBB(\EE,E)}  \le  C'_a \, e^{at}.
\]
   \end{itemize}
   
 \smallskip\noindent
 The following equivalence holds:
 
 \begin{itemize}
 \item[{\bf (1)}] There exists a finite rank projector $\Pi_{L} \in \BBB(E)$ which commutes with $L$ and satisfy $\Sigma(L_{|\Pi_{L}}) = \{ 0 \}$,  
 so that the semigroup $S_L =e^{tL}$ satisfies  the growth estimate
 \begin{equation}\label{eq:estimSGL} 
   \forall \, t \ge 0, \quad \left\|  S_L(t) -  \Pi_{L} \right\|_{\BBB(E)}
   \le C_{L,a} \, e^{a\, t}  
 \end{equation}
for any $a > a^{**}$ and some constant $C_{L, a}>0$;

\item[{\bf (2)}] There exists a  finite rank  projector  $\Pi_{\LL} \in \BBB(\EE)$  which commutes with $\LL$ and satisfy $\Sigma(\LL_{|\Pi_{\LL}}) = \{ 0 \}$,  
 so that the semigroup $S_\LL =e^{t\LL}$ satisfies  the growth estimate 
 \begin{equation}\label{eq:estimSGLL} 
   \forall \, t \ge 0, \quad \left\|  S_\LL(t) -  \Pi_{\LL} \right\|_{\BBB(\EE)}
   \le C_{\LL, a} \, e^{a \, t}  
 \end{equation}
for any $a > a^{**}$ and some  constant $C_{\LL, a} >0$.

    \end{itemize}
\end{theo}
\medskip

\noindent{\sl Proof of  Proposition~\ref{theo:K=Cste}. } {\sl Step 1. } We recall some facts presented in \cite{MR2114128,MR2536450}.
We introduce the rescaled function $g(t,x) := f(t,x) \, e^{-\lambda t}$ and the associated rescaled  equation
\beqn\label{eq:RescGFcst}
\partial_t g + \partial_x g + n_F K_0 g  = K_0 \int_x^\infty \kappa(y,x) \, g(y) \, dy
\eeqn
with vanishing boundary condition \eqref{eq:BoundaryCondition} and initial condition $g(0) = f_0$. We observe that the number of particles 
$$
\int_0^\infty g(t,x) \, dx
$$
is conserved. One can then show using the Tikhonov's infinite dimensional version of the Brouwer fixed point Theorem that there exists a steady state $f_\infty$ by proceeding 
exactly as for the self-similar fragmentation equation \cite[Section 3]{EMRR} (see also  \cite{MR2114128,MR2536450} where other arguments are
presented). Existence of the steady state $f_\infty$ is also given by the Krein-Rutman Theorem presented in section~\ref{subsec:KRGrFr}.
This steady sate corresponds to the first eigenfunction associated to the first eigenvalue $(n_F-1) \, K_0$ of the cell-division equation \eqref{eq:GrowthFragCst}.  

Anyway, under assumption \eqref{eq:intkappa<0}, it has been shown during the proof of  \cite[Theorem 1.1]{MR2114128} 
and \cite[Theorem~1.1]{MR2536450} that the solution $g$ to \eqref{eq:RescGFcst} satisfies
$$
\| g(t) - \langle f_0 \rangle \, f_\infty \|_{-1,1} \le e^{- \lambda t} 
\| f_0 - \langle f_0 \rangle \, f_\infty \|_{-1,1} 
$$
where for any $f \in L^1_1$ with mean $0$ we have defined
$$
\| f \|_{-1,1} := \int_0^\infty \Bigl| \int_0^x f(y) \, dy \Bigr| \, dx .
$$

\medskip\noindent
{\sl Step 2.  } For the mitosis equation, we introduce the splitting $\Lambda = \AA + \BB$ where
$$
\AA := \FF^{+}_R, \quad \BB := - \partial_x - n_F \, K_0 + \FF^{+,c}_R,
$$
with the notation of section~\ref{subsec:CellDiv}.  We define $E := L^1_\alpha$, $\alpha > 1$, and for any $f \in L^1_\alpha$ with mean $0$ we define 
$$
\| f \|_{-1,\alpha} := \int_0^\infty \Bigl| \int_0^x f(y) \, dy \Bigr| \, \langle x \rangle^{\alpha - 1} \, dx,
$$
as well as  the Banach space  $\EE$ obtained by completion of $L^1_\alpha$ with respect to the norm $\| \cdot \|_{-1,\alpha}$. 
Let us explain now why the hypothesis {\bf (i)}, {\bf (ii)} and {\bf (iii)} in Theorem~\ref{theo:Extension} are fulfilled. 
We clearly have that $\AA$ satisfies {\bf (ii)} in both spaces  $E$ and $\EE$, as well as  that  $\BB$ satisfies {\bf (i)}  in the space $E$ for any   $a > a^*$, 
$a^* := - K_0 (1 - \wp_\alpha)$ thanks to Lemma~\ref{lem:BorneBB}. We claim that  $\BB$ also satisfies {\bf (i)}  in the space $\EE$ for any $\alpha > \log_2 3$ and any $a > a' := 2K_0 \, (3 \times 2^{-\alpha} -1) $.  In order to prove that last claim we proceed along the line of the proof of \cite[Theorem~1.1]{MR2114128}. 
For any $g \in \EE$, we  introduce the notation
$$
g(t) := e^{\BB t} g, \quad G(t,x) := \int_0^x g(t,u) \, du,
$$
and we compute  
\bean
&&\partial_t G + \partial_x G +2  \, K_0 G 
= - 2 K_0 \int_x^\infty \chi^c_R(2y) \, \partial_y G(t,2y) \, dy
\\
&&\qquad =  2 K_0  \chi^c_R(2x) \, G(t,2x)  + 4 K_0\int_x^\infty (\chi^c_R)' (2y) \, G(t,2y) \, dy.
\eean
Similarly as in the proof of Lemma~\ref{lem:BorneBB}, we define 
$$
\phi(x) :=  {\bf 1}_{[0,x_2]} + { x^{\alpha-1} \over x_2^{\alpha-1} } \, {\bf 1}_{[x_2,\infty)}
$$
and we compute for $R \ge x_2$ 
\bean
&&\partial_t \int_0^\infty |G| \, \phi + \int_0^\infty |G| \, \phi \,  { \alpha -1 \over x} \, {\bf 1}_{x \ge x_2}  +2  \, K_0 \int_0^\infty |G| \,  \phi 
\\
&&\qquad \le  {K_0 \over 2^{\alpha-1}}  \int_0^\infty |G| \, \phi \, \chi^c_R   + {K_0 \over 2^{\alpha-1}} \int_0^\infty |(\chi^c_R)'| \, x \, |G| \, \phi \, dx  
\\
&&\qquad \le  3 \,  {K_0 \over 2^{\alpha-1}}   \int_0^\infty |G| \, \phi ,
\eean
where in the last line we have used that $| x \, (\chi^c_R)'| \le 2$ by definition of $\chi^c_R$. We conclude by taking $x_2 = R$ large enough. 
 
Moreover, we claim that   
\beqn\label{eq:U2explicit}
\Bigl\|  (\AA S_\BB)^{(*2)}(t)  \Bigr\|_{\BBB(\EE,E)}  \le  C' \,  (1+t) \, e^{- \mu t},
\eeqn
with $\mu := n_F K_0 = 2K_0$. 
In order to prove estimate \eqref{eq:U2explicit}, as in the proof of Lemma~\ref{lem:regMitose} and with the same notation, we compute starting from \eqref{eq:U*2-2}
\bean
U_0^{(*2)} (t) g (x) &=& 16 K_0^2 \chi_R(2x) \, e^{-\mu t} \int_{u_0}^{u_1} \chi_R(2u-4x+2t)  g(u) \, du 
\\
&=& 16 K_0^2 \chi_R(2x) \, e^{-\mu t} \,\Bigl\{ \Bigr[   \chi_R(2u-4x+2t)  G(u) \Bigr]_{u=u_0}^{u=u_1} 
\\&&\qquad- 
{2 \over R} \int_{u_0}^{u_1} \chi'_R(2u-4x+2t)  G(u) \, du \Bigr\} ,
\eean
and we get then for any $\beta \ge 0$
\bean
\| U_0^{(*2)} (t) g \|_{L^1_\beta} 
\le C \, (1+t)   \, e^{-\mu t} \, \| G \|_{L^1}. 
\eean
We have a similar estimate for $U_0 * U_0^c$ and we then obtain \eqref{eq:U2explicit} thanks to formula \eqref{eq:U*2-1}. 
We conclude by using the shrinkage of functional space result stated in Theorem~\ref{theo:Extension}, and we get for possible definition of  $\alpha^*$ the unique real number such that $3 \wp_{\alpha^*} - n_F= [a^{**}- (n_F-1) K_0]/K_0$ for any fixed $a^{**} \in (-  K_0, (n_F-1) K_0)$.

\medskip\noindent
{\sl Step 3.  } For the cell division equation with smooth offspring distribution $\wp$ we can proceed along the line of \cite[Theorem~1.1]{MR2536450} and of step 2. 
We introduce the same splitting as for the mitosis equation and we work in the same spaces.  
We clearly have  again that $\AA$ satisfies {\bf (ii)} in both spaces  $E$ and $\EE$ and $\BB$ satisfies {\bf (i)} in  $E$. 
We claim that  $\BB$ also satisfies {\bf (i)}  in the space $\EE$ for any  $a > a'  \in (-n_F K_0, 0)$ and any  $\alpha > \alpha'$ where 
$\alpha' > 1$ is such that $a' / K_0= - n_F + {4 \over \alpha'} n_F + \wp_{\alpha'-1}$.   Indeed, we have 
\bean
\partial_t G + \partial_x G + n_F  \, K_0 G 
&=& -   K_0 \int_x^\infty  \int_y^\infty  \kappa^c_R(z,y )  \, \partial_z G(z) \, dz dy
\\
&=&   - K_0  \int_x^\infty \beta^c_R(z,x) \, G(z)  \, dz, 
\eean
 with 
\bean
\beta^c_R(z,x) 
&:=& -   { \partial \over \partial z}  \int_x^\infty \kappa^c_R(z,y) \, dy
\\
&=& - n_F(x/z) \,   ( \chi^c_R(z) )' + \chi^c_R(z) \, \beta(z,x) ,
\eean
where 
\bean
\beta(z,x) 
:= -   { \partial \over \partial z}  \int_x^\infty \kappa(z,y) \, dy, \quad
n_F(u) := \int_u^1 \wp(u') \, du',
\eean
and where $\kappa^c_R(x,y) = \chi^c_R(x) \, \kappa (x,y)$ is defined on $\R^2_+$ by extended it  to $0$ outside of the set $\{ (x,y) \in \R^2; \, 0 < y < x \}$. 
On the one hand, we have 
\bean
\Phi_1(z) 
&:=&   \int_0^z  n_F(x/z) \,   |( \chi^c_R(z) )' |  \, \phi(x) \, dx 
\\
&\le&  n_F \, |( \chi^c_R(z) )' |  \, (x_2 + z^\alpha/\alpha) \le {4 \over \alpha} \, n_F \,  \phi(z)
\eean
for any $z \ge 0$ if $x_2 \le R^\alpha/\alpha$. On the other hand, we have 
\bean
\Phi_2(z) 
&:=&   \int_0^z  \chi^c_R(z) \, \beta(z,x) \, \phi(x) \, dx 
\\
&\le&  \chi^c_R(z) \Bigl\{ \eta(x_2/R) + \wp_{\alpha-1} \, \phi(z) \Bigr\}  , \quad \eta(u) := \int_0^u \wp(u') \, du', 
\eean
for any $z \ge 0$.  
Next, we compute 
\bean
&&\partial_t \int_0^\infty |G| \, \phi + \int_0^\infty |G| \, \phi \,  { \alpha -1 \over x} \, {\bf 1}_{x \ge x_2}  + n_F \, K_0 \int_0^\infty |G| \,  \phi 
\\
&&\qquad \le K_0   \int_0^\infty  \Bigl\{ \int_0^z  | \beta^c_R(z,x) | \, \phi(x) \, dx \Bigr\} \, |G(z)| \, dz
\\
&&\qquad \le K_0   \int_0^\infty  \Bigl\{ \Phi_1(z) + \Phi_2(z) \Bigr\} \, |G(z)| \, dz
\\
&&\qquad \le  K_0 \Bigl\{ {4 \over \alpha} \, n_F +  \eta(x_2/R) + \wp_{\alpha-1}   \Bigr\}  \int_0^\infty |G| \, \phi  \, dx ,
\eean
and we take $\alpha > 1$ large enough and next $R/x_2$ large enough.   

\smallskip
We next claim that $U_0 := \AA S_{\BB_0} $ with $\BB_0 := - \partial_x - \mu$, $\mu := n_F \, K_0$, satisfies
\beqn\label{eq:U0WL}
\| U_0 (t) g \|_{L^1_\beta} \le  C \, \bigl( 1 + {1 \over t } \bigr)  \,e^{-\mu t}  \| g \|_{-1,1}.
\eeqn
Starting from the definition 
\bean
(U_0 (t) )(x)  &:= & e^{-\mu t} \, K_0 \int_x^\infty \kappa_R(y,x) \, \partial_y G(y-t) \, dy
\\
&=& - e^{-\mu t} \,  K_0 \, \kappa_R(x,x) \, G(x-t) 
 - e^{-\mu t} \, K_0 \int_x^\infty \partial_y [\kappa_R(y,x)]  \,  G(y-t) \, dy,
\eean
we compute 
\bean
\| U_0 (t) g \|_{L^1_\beta} &\le  & K_0 \, e^{-\mu t}   \int_0^\infty  \chi_R(x) \, \langle x \rangle^\beta \, \wp (1) \, {|G(x-t)| \over x - t + t} \, dx 
\\
&& \!\!\!\!\!\!\!\!\!\!\!\!
+ K_0 \, e^{-\mu t}  \int_0^\infty \int_0^y \Bigl\{ |(\chi_R)'| \, {1 \over y} \wp({x \over y}) + \chi_R(y) \, {x \over y^3} \, |\wp'( {x \over y})|  \Bigr\} \langle x \rangle^\beta \, G(y-t) \, dy
\\
&\le  & K_0 \, e^{-\mu t} \, C_{R,\beta}  \, \bigl( 1 + \frac1t \bigr)  \int_0^\infty   \, |G(z) | \, dz,
\eean
and that ends the proof of \eqref{eq:U0WL}. We introduce the notation $\EE_1 := E$, $\EE_0 := \EE$ and $\EE_{1/2}$ as the $1/2$ complex interpolation between 
the spaces $\EE_0$ and $\EE_1$. From the above estimates, we have for any $a > - \mu$ that 
$$
 \| U_0(t) \|_{\EE_j,\EE_{j+1/2}} \le C \, t^{-1/2} \, e^{at}, \quad \hbox{for}\ j=0,1/2. 
$$
Thanks to \eqref{eq:U*2-1} it is not difficult now to prove that \eqref{eq:U2explicit} holds also in the present case.  We conclude again by using the shrinkage of functional space result stated in Theorem~\ref{theo:Extension}. 
\qed

%---------------------------------------------------------------------------------------------------------------------------------------------------------------------

\subsection{Quantified spectral gap for the self-similar fragmentation equation with positive kernel}

We present a second situation where a very accurate and quantitative description of the spectrum  (as formulated in point (ii) of Theorem~\ref{theo:Frag}) is possible.

\begin{prop}\label{prop:KRquant} Consider the self-similar fragmentation equation \eqref{eq:fragSSV} and assume that the fragmentation kernel
satisfies \eqref{eq:hypSSfrag}, \eqref{eq:wpsmooth} and \eqref{eq:wp>0}. Then, there exists a computable constant $a^{**} \in (a^*,0)$
such that the spectral gap $\Sigma(\Lambda) \cap \Delta_{a^{**}} = \{ 0 \}$ holds in the functional space $X = \dot L^1_\alpha \cap \dot L^1_\beta$, $0 \le \alpha < 1 < \beta$. 
 
\end{prop}

   \noindent
   {\sl Proof of Proposition~\ref{prop:KRquant}. } We split the proof into four steps. 
   
 \medskip\noindent
 {\sl Step 1. A priori bounds. }  We fix $a = a^*/2$, where $a^*$ is the same as in  Proposition~\ref{prop:FragSpectralMapping}. From Proposition~\ref{prop:FragSpectralMapping} and then Theorem~\ref{theo:QuantWeyl}, there exists a constant $R_a$ such that for any  eigenfunction $f$ associated to an eigenvalue $\xi \in \Sigma(\Lambda) \cap \Delta_a$,  which satisfies the normalization condition
 \beqn\label{eq:normalizationJordanFSS}
 \int_0^\infty |f(y)| \, y \, dy = 1, 
\eeqn
% in a Jordan basis associated to the eigenspace $R\Pi_{\Lambda,-1}$
there holds 
\beqn\label{eq:BdJordanFSS1}
\int_0^\infty |f(y)| \, \langle y \rangle^{2+\gamma} \, dy  \le R_a, \quad |\xi| \le R_a .
\eeqn
Together with the eigenvalue problem rewritten as  
$$
\partial_x( x^2 f) = x^2 \, [ x^\gamma \, f + \xi \, f - \FF^+ \!f],
$$
we deduce that $x^2 f \in BV(\R_+) \subset L^\infty(\R_+)$ and then,  iterating the argument, that for any $\delta \in (0,1)$ there exists $C_\delta$ such that 
\beqn\label{eq:BdJordanFSS2}
\| f \|_{W^{1,\infty}(\delta,1/\delta)} \le C_\delta. 
\eeqn

 \medskip\noindent
 {\sl Step 2. Positivity. } From \eqref{eq:normalizationJordanFSS} and \eqref{eq:BdJordanFSS1} we clearly have 
  \beqn\label{eq:positivityJordanFSS1}
 \int_\delta^{1/\delta}   |f(y)| \, y \, dy \ge  1 - 2 \, R \, \delta.
\eeqn
Taking $\delta_1 := 1/(4R)$, we see that there exists at least one point $x_1 \in (2\delta_1,1/(2\delta_1))$ such that 
  \beqn\label{eq:positivityJordanFSS2}
 |f(x_1) | \ge \delta_1^2/2.  
\eeqn 
Introducing the logarithm function $\theta$ defined by $f(y) = |f(y)| \, e^{i \theta(y)}$ which is locally well defined for any $y \in (0,\infty)$ such that $f(y) \not = 0$, we see from \eqref{eq:positivityJordanFSS2} and \eqref{eq:BdJordanFSS2} that there exists an interval $I_1 \subset (\delta_1,1/\delta_1)$, $x_1 \in I_1$,  and a computable real number $\eps_1 := \eps(\delta_1) > 0$ such that 
\beqn\label{eq:positivityFSS3}
|I_1| \ge \eps_1 \quad\hbox{and}\quad \Re e (f(x) \, e^{- i \theta(x_1)} ) \ge \delta_1^2/4 \quad \forall \, x \in I_1, 
\eeqn
and better (since $\cos$ and $|f(x)|$ are Lipschitz functions)
\beqn\label{eq:positivityFSS4}
|I_1| \ge \eps_1 \quad\hbox{and}\quad |f(x)| \ge  \delta_1^2/4, \quad |\theta(x_1) - \theta(x)| \le \pi/4 \quad \forall \, x \in I_1. 
\eeqn
 
On the other hand, we know that $\Pi_{\Lambda,0} f = 0$ which yields  
$$
\int_0^\infty \Re e \Bigl\{ f(y) \, e^{- i \theta(x_1)} \Bigr\} \, y \, dy = 0.
$$
From \eqref{eq:BdJordanFSS1} again, we deduce that 
$$
\int_\delta^{1/\delta}  \Re e \Bigl\{ f(y) \, e^{- i \theta(x_1)} \Bigr\} \, y \, dy \le 2 \, R \, \delta. 
$$
 But together with  the positivity property \eqref{eq:positivityFSS3}, there exists $\delta_2 < \delta_1$ such that $2R\delta_2 < \kappa_2 := \eps_1 \, \delta^2_1 / 8 > 0$ and %we must have for any $\delta_2 < .. $  
 for any $\delta \in (0,\delta_2)$ 
 $$
\int_\delta^{1/\delta}  \Bigl( \Re e \Bigl\{ f(y) \, e^{- i \theta(x_1)} \Bigr\} \Bigr)_- \, y \, dy \ge \kappa_2.
$$
Using the same arguments as above, there exists $x_2 \in (\delta_2,1/\delta_2)$ such that 
$$
|f(x_2)| \, \cos (\theta(x_2) - \theta(x_1)) \le - \kappa_2 \, \delta^2_2/2,
$$
and then there exist an interval  $ I_2 \subset (\delta_2,1/\delta_2)$ and some constants $\eps_2,\kappa_2 > 0$ such that 
%
% $$
%C \, \int_\delta^{1/\delta}  \Bigl( \cos (\theta(y) -  \theta(x_1))\Bigr)_- \, y\, dy \ge \eps_1 \, \delta^2_1/4 ,
%$$
%and then using the same argument as above, there exists $ I_2 \subset (\delta,1/\delta)$ such that 
$$
|I_2| \ge \eps_2 \quad\hbox{and}\quad  |f(y) |  \ge  \kappa_2 \quad \forall \, y \in I_2,
$$
%for some computable positive constants $\eps_2$ and $\kappa_2$.  And then also 
% $$
%|I_2| \ge \eps_2 \quad\hbox{and}\quad |f(y)| \ge  \kappa_2, \quad |\theta(x_1) - \pi \theta(y)| \le \pi/4 \quad \forall \, y \in I_2. 
%$$
as well as 
$$
\cos[\theta(y) - \theta (x)] \le 0   \quad \forall \, x \in I_1, \,\, \forall \, y \in I_2. 
$$
 We may assume without loss of generality that $x_1 > x_2$.

 \medskip\noindent
 {\sl Step 3}. By definition of the growth-fragmentation operator $\Lambda$ and denoting $\overline{\rm sign} f :=\bar f/|f|^2$ where $\bar f$ stands for the complex conjugate of $f$, we clearly have 
\bean
\int_0^\infty (\Lambda |f| - (\Lambda f) \, \overline{\rm sign} f) \, \phi =
\int_0^\infty\!\!\int_0^y k(y,x) \, |f(y)| (1 - {\rm sign} f(y) \,  \overline{\rm sign} f(x)) \, \phi(x) \, dxdy.
\eean
Since 
%$$
%\Re e \{1 - {\rm sign} f(y) \,  \overline{\rm sign} f(x) \}   \ge 0 \quad \forall \, x, y \ge 0
%$$
%and more accurately 
$$
\Re e \{1 - {\rm sign} f(y) \,  \overline{\rm sign} f(x) \}  = 1 - \cos [\theta(y) - \theta(x)] \ge 1  \quad \forall x \in I_1, \,\, \forall \, y \in I_2, 
$$
we deduce that 
\bean
&&\Re e \int_0^\infty (\Lambda |f| - (\Lambda f) \, \overline{\rm sign} f) \, \phi   
\\
&&\quad\ge \wp_* \,  \delta_2^{\gamma-1}  \int_{I_2} |f(y)| \int_{I_1} \Bigl\{ 1 - \cos [\theta(y) - \theta(x)]  \Bigr\}  \, x  \, dx \Bigr\} \, dy 
\\
&&\quad \ge \wp_* \,  \delta_2^{\gamma-1} \, \delta_1 \, \eps_1 \, \eps_2 \, \kappa_2 =:  - a^{**}.
 \eean

 \medskip\noindent
 {\sl Step 4. } Conclusion. For any mass normalized eigenvector $f \in D(\Lambda)$ associated to an eigenvalue $\xi \in \Delta_a \cap \Sigma(\Lambda) \backslash \{ 0 \}$, there holds thanks to step 3
\bean
\Re e \xi \, \langle |f|, \phi \rangle 
&=&
\Re e  \langle \xi f \, \overline{\rm sign} f , \phi \rangle  
\\
&=&
\Re e  \langle \Lambda f \, \overline{\rm sign} f , \phi \rangle  
\\
&\le& \langle \Lambda |f| , \phi \rangle + a^{**}
\eean
and then $\Re e  \xi \le a^{**}$. As a consequence, $\Delta_{a^{**}} \cap \Sigma(\Lambda) = \{ 0 \}$ and we conclude thanks to Theorem~\ref{theo:GalWeyl}. 
\qed

%\bigskip
%\bibliographystyle{acm}
%%\bibliography{biblio,mybiblio}
%\bibliography{./SGandCS}
%\end{document}

\bigskip%\newpage
\bibliographystyle{acm}
%\bibliography{biblio,mybiblio}

%\begin{thebibliography}{100}

\bigskip
\bigskip
 \signsm \signjs 

\end{document}